\documentclass[11pt,reqno]{amsart}

\usepackage[T1]{fontenc}
\usepackage[utf8]{inputenc}
\usepackage{lmodern}
\usepackage{amsmath,amssymb,amsthm,mathtools}
\usepackage{microtype}
\usepackage{enumitem}
\usepackage{tikz-cd}
\usepackage{xcolor}
\usepackage[colorlinks=true,linkcolor=blue!45!black,citecolor=blue!45!black,urlcolor=blue!45!black]{hyperref}
\usepackage[capitalise,noabbrev]{cleveref}

\setlist[itemize]{leftmargin=2em,itemsep=.25em,topsep=.35em}
\setlist[enumerate]{leftmargin=2.3em,itemsep=.25em,topsep=.35em}
\allowdisplaybreaks
\newtheorem{theorem}{Theorem}[section]
\newtheorem{proposition}[theorem]{Proposition}
\newtheorem{lemma}[theorem]{Lemma}
\newtheorem{corollary}[theorem]{Corollary}

\theoremstyle{definition}
\newtheorem{definition}[theorem]{Definition}
\newtheorem{example}[theorem]{Example}

\theoremstyle{remark}
\newtheorem{remark}[theorem]{Remark}

\crefname{criterion}{criterion}{criteria}
\Crefname{criterion}{Criterion}{Criteria}

\DeclareMathOperator{\Ob}{Ob}
\DeclareMathOperator{\Mor}{Mor}

\DeclareMathOperator{\defi}{def}

\DeclareMathOperator{\colim}{colim}
\DeclareMathOperator*{\hocolim}{hocolim}

\DeclareMathOperator{\im}{im}

\newcommand{\Set}{\mathsf{Set}}
\newcommand{\Cat}{\mathsf{Cat}}
\newcommand{\CAT}{\mathsf{CAT}}
\newcommand{\Gpd}{\mathsf{Gpd}}
\newcommand{\Grph}{\mathsf{Grph}}
\newcommand{\Top}{\mathsf{Top}}
\newcommand{\Comp}{\mathsf{Comp}}

\newcommand{\F}{\mathcal F}
\newcommand{\Pp}{\mathcal P}
\newcommand{\D}{\mathcal D}

\newcommand{\A}{\mathcal A}

\newcommand{\K}{\mathcal K}

\newcommand{\N}{\mathbb N}
\newcommand{\Z}{\mathbb Z}

\newcommand{\1}{\mathbf 1}
\newcommand{\2}{\mathbf 2}
\newcommand{\walking}[1]{\mathbf 2_{#1}}
\newcommand{\copow}{\mathbin{\boldsymbol\cdot}}

\newcommand{\op}{\mathrm{op}}
\newcommand{\freecat}{\F_1}
\newcommand{\freegpd}{\F_{\mathrm{gpd}}}
\newcommand{\FTopOne}{\F_{\Top_1}}
\newcommand{\CTopOne}{\mathcal C_{\Top_1}}
\newcommand{\ITopOne}{\mathcal I_{\Top_1}}
\newcommand{\DTop}{D_{\Top}}
\newcommand{\gpcomp}{\mathsf L}
\newcommand{\PiOne}{\Pi_1}

\newcommand{\abs}[1]{\lvert #1\rvert}

\newcommand{\defeq}{\coloneqq}

\title[Freely generated $n$-categories]{Freely generated $n$-categories, coinserters and presentations of low dimensional categories}
\author{Fernando Lucatelli Nunes}
\address{CMUC, Department of Mathematics, University of Coimbra, 3000-143 Coimbra, Portugal}
\email{fln@uc.pt}
\date{}
\hypersetup{
	pdftitle={Freely generated n-categories, coinserters and presentations of low dimensional categories},
	pdfauthor={Fernando Lucatelli Nunes},
	pdfsubject={Computads, coinserters, presentation complexes, and deficiency},
	pdfkeywords={computad, polygraph, coinserter, icon, categorical presentation, deficiency, fundamental groupoid, crossed module, rewriting systems}
}


\begin{document}
	
	\begin{abstract}
A presentation records not only a categorical structure, but also how that structure is assembled. This information is essential for rewriting, coherence, and questions of minimality: freely adjoining a cell with prescribed boundary is a fundamentally different operation from imposing an equation between cells already constructed. We show that these two stages are governed, respectively, by coinserters and coequifiers. The resulting framework treats computadic presentations uniformly, from ordinary categories to strict higher categories, and, in low dimensions, relates their generating and relation data to cellular attachments.
		
		The starting point is elementary.  If $G$ is a graph, the category
		$\freecat G$ freely generated by $G$ is the coinserter in $\Cat$ of its
		domain and codomain maps, regarded as functors between discrete categories.
		For every
		$n\geq1$, freely adjoining $n$-cells with prescribed parallel boundaries to
		a strict $(n-1)$-category is a coinserter in the $2$-category of strict
		$n$-categories, strict $n$-functors, and $n$-icons; here the two functors
		agree below dimension $n-1$, and the icon has only $n$-dimensional
		components.  The construction is left adjoint to the underlying
		derivation-scheme functor.
		When the lower skeleton is freely generated by a computad, it gives the
		usual free strict $n$-category, while coequifiers impose equations between
		the freely generated cells.  In dimension two, the same object satisfies a
		bicategorical universal property for normal pseudofunctors and icons.
		
		Replacing the walking arrow by the unit interval gives the topological
		coinserter of $G$.  The topological comparison then starts from a groupoidal
		$2$-computad $C$.
		Attaching one disk to $\abs{C_{\leq1}}$ for every relation gives a
		presentation complex $X_C$, and the groupoid presented by $C$ is canonically
		isomorphic to $\PiOne(X_C,C_0)$.  It follows that the rank-finite deficiency
		of a connected groupoid agrees with the classical deficiency of any of its
		isotropy groups; the finite version holds when the object set is finite.
		Homology then gives a sharp lower bound for the number of relations needed
		to present a thin groupoid on a fixed connected rank-finite graph.  Crossed
		modules carry the comparison one dimension further and detect relations
		among relations.  For the descent computad, the identity and associativity
		confluences form a homotopy basis and attain the resulting homological lower
		bound relative to each of the two fixed $2$-skeleta considered here.
	\end{abstract}
	
	\maketitle
	\clearpage
	\tableofcontents
	\clearpage
	
	\section{Introduction}
	\label{sec:introduction}

	A presentation answers a question which the presented object no longer
	answers: how was it built?  This question matters whenever the construction,
	and not only its quotient, enters the argument.  It is what permits one to
	perform calculations by rewriting, to formulate coherence as a problem about higher
	relations, and to ask whether a chosen family of relations is irredundant.
	The point of the present paper is that, for categorical structures, this
	information is intrinsically dimensional and can be organized by familiar
	universal constructions.
	
	The idea is simple.  First one freely adjoins cells with specified
	boundaries.  Only then does one impose equations between the cells so
	obtained.  Coinserters describe the first operation; coequifiers describe
	the second.  The distinction is already visible for ordinary categories,
	but it becomes indispensable in dimension two and above: a generating
	$2$-cell has two parallel paths as its boundary, whereas a relation between
	$2$-cells is one dimension higher.
	
	Let us begin with the smallest example in which both operations occur.  By a
	\emph{graph} we mean a set $G_0$ of objects, a set $G_1$ of arrows, and
	domain and codomain maps $s,t\colon G_1\rightrightarrows G_0$.  Consider the
	graph
	\begin{equation}
		\label{eq:o-triangulo-que-explica-tudo}
		\begin{tikzcd}[column sep=large]
			x \ar[r,"a"] \ar[rr,bend right=22,"c"']&
			y \ar[r,"b"]&
			z
		\end{tikzcd}.
	\end{equation}
	In the free category, the path $ba$ and the arrow $c$ are distinct.  We may
	either adjoin a $2$-cell between them,
	\begin{equation}
		\label{eq:a-celula-do-triangulo}
		\begin{tikzcd}[column sep=huge]
			x
			\ar[r,bend left=38,"ba",""{name=upper,below}]
			\ar[r,bend right=38,"c"',""{name=lower,above}]&
			z
			\arrow[Rightarrow,from=upper,to=lower,"\alpha" description]
		\end{tikzcd},
	\end{equation}
	or impose the equation $ba=c$.  The first construction remembers
	$\alpha$ as a genuine $2$-cell; the second identifies its boundary paths.
	Topologically, the graph in \eqref{eq:o-triangulo-que-explica-tudo} is a
	circle subdivided into three intervals.  The relation $ba=c$ attaches a disk
	along that circle.  Thus the categorical equation and the cellular
	attachment encode the same boundary datum, although they live in different
	categories.  This example is the model for the whole paper.
	
	We now give the universal construction behind it.  Let $\mathbb G$ be the
	category with objects $\mathsf 0,\mathsf 1$ and two arrows
	$d^0,d^1\colon\mathsf 0\rightrightarrows\mathsf 1$.  Equivalently, a graph
	is a functor
	\[
	G\colon\mathbb G^{\op}\longrightarrow\Set.
	\]
	We write $G_0=G(\mathsf 0)$ and $G_1=G(\mathsf 1)$ for its sets of objects
	and arrows.  The category $\freecat G$ freely generated by $G$ has object
	set $G_0$, and its arrows are the finite composable paths in $G$.
	
	Let $\D\colon\Set\to\Cat$ be the discrete-category functor.  Recall that a
	\emph{coinserter} of parallel functors
	$f,g\colon A\rightrightarrows B$ consists of a functor $q\colon B\to Q$ and
	a natural transformation $\theta\colon qf\Rightarrow qg$ with the following
	universal property: a functor $Q\to X$ is equivalently a functor
	$h\colon B\to X$ together with a natural transformation
	$hf\Rightarrow hg$, and the same correspondence holds for transformations
	between such functors.  The category $\freecat G$ has precisely this
	universal property:
	\begin{equation}
		\label{eq:the-first-coinserter}
		\begin{tikzcd}[column sep=large]
			\D G_1
			\ar[r,shift left=.65ex,"{\D s}"]
			\ar[r,shift right=.65ex,swap,"{\D t}"]&
			\D G_0 \ar[r,"q"]&
			\freecat G,
		\end{tikzcd}
		\qquad
		q\D s\xRightarrow{\ \theta\ }q\D t .
	\end{equation}
	More precisely, for every arrow $e\in G_1$, the component
	\[
	\theta_e\colon q(s(e))\longrightarrow q(t(e))
	\]
	is the generating arrow corresponding to $e$.  Thus the coinserter adjoins
	the arrows with their prescribed domains and codomains and closes them under
	composition.  Relations are imposed only afterwards, by a coequifier.
	
	There is an analogous weighted colimit in topology.  More precisely, the
	categorical weight is the walking arrow, while the topological weight is the
	interval with its two distinguished endpoints.  Write $D^1=[0,1]$ and
	$S^0=\{0,1\}=\partial D^1$.  If $S$ is a set and $X$ is a space,
	$S\copow X$ denotes the disjoint union of one copy of $X$ for each element of
	$S$.  Regarding $G_0$ and $G_1$ as discrete spaces, the realization
	$\abs G$ is the pushout
	\begin{equation}
		\label{eq:graph-in-the-introduction}
		\begin{tikzcd}
			G_1\copow S^0 \ar[r,"{s\sqcup t}"] \ar[d,hook] &
			G_0 \ar[d] \\
			G_1\copow D^1 \ar[r] & \abs{G}.
		\end{tikzcd}
	\end{equation}
	In \cref{sec:the-big-picture} we express this pushout as an
	interval-weighted colimit in the category $\Top$ of compactly generated weak
	Hausdorff spaces.  This is the topological coinserter of $G$.  The two
	coinserters lie in different categories, but both realize the prescribed
	domain and codomain data.  A categorical relation between parallel paths is
	then realized by attaching a disk along the loop obtained from those paths.
	This gives the usual presentation complex.
	
	The example suggests the two questions treated here.  First, does the
	coinserter description of the free category on a graph persist when one
	adjoins cells to a strict higher category?  Second, what can the topology of
	the resulting presentation complexes tell us about the number of relations
	and of relations among relations?  Part~I answers the first question;
	Part~II answers the second.  The two answers fit together because they begin
	with the same boundary data.
	
	\subsection*{Main results}
	
	\paragraph{Free attachment and equations.}
	
	The coinserter calculation does not require the lower-dimensional structure
	to be free.  Let $n\geq1$, let $B$ be a strict $(n-1)$-category, and let
	$S$ be a set.  Write $\walking{n-1}$ for the strict
	$(n-1)$-category freely generated by one $(n-1)$-cell.  A pair of functors
	\[
	b^-,b^+\colon S\copow\walking{n-1}\rightrightarrows B
	\]
	which agree below dimension $n-1$ selects two parallel
	$(n-1)$-cells for every $s\in S$.  We call this datum a
	\emph{derivation $n$-scheme} $D=(B,S,b^-,b^+)$.  If
	$\kappa_{n-1,s}$ denotes the distinguished $(n-1)$-cell in the $s$th
	copower summand, the free extension by $n$-cells adjoins one
	$n$-cell
	\[
	b^-(\kappa_{n-1,s})\Longrightarrow b^+(\kappa_{n-1,s})
	\]
	for every $s$, subject only to the strict $n$-categorical axioms.
	
	For $n=1$ the appropriate ambient $2$-category is $\Cat$, with ordinary
	natural transformations as $2$-cells.  For $n\geq2$, the appropriate
	ambient $2$-category has strict $n$-categories as objects, strict
	$n$-functors as $1$-cells, and \emph{$n$-icons} as $2$-cells.  An $n$-icon
	between strict $n$-functors which agree through dimension $n-2$ assigns an
	$n$-cell to every $(n-1)$-cell, compatibly with identities, compositions,
	and the $n$-cells already present.  If $d_nB$ denotes $B$ regarded as a
	strict $n$-category with only identity $n$-cells, the free extension is the
	coinserter
	\begin{equation}
		\label{eq:higher-picture-introduction}
		\begin{tikzcd}[column sep=large]
			d_n(S\copow\walking{n-1})
			\ar[r,shift left=.65ex,"{d_nb^-}"]
			\ar[r,shift right=.65ex,swap,"{d_nb^+}"]&
			d_nB \ar[r]&
			\operatorname{Ext}_n(D).
		\end{tikzcd}
	\end{equation}
	The components of the universal $n$-icon are exactly the newly adjoined
	$n$-cells.  Moreover, free extension by $n$-cells is left adjoint to the
	underlying derivation-scheme functor, which records a strict
	$(n-1)$-category together with all the $n$-cells and their boundaries
	(\cref{thm:derivation-scheme-coinserter,thm:as-derivacoes-sabem-o-caminho}).
	
	An \emph{$n$-computad} is the recursive case in which the lower category is
	itself freely generated.  It consists of a lower-dimensional computad
	$C_{<n}$, a set $C_n$ of generating $n$-cells, and two parallel
	$(n-1)$-cells in $\F_{n-1}C_{<n}$ as the boundary of every generator.
	Taking $B=\F_{n-1}C_{<n}$ in
	\eqref{eq:higher-picture-introduction} gives the strict $n$-category
	$\F_nC$ freely generated by $C$.  The one-step attachment theorem is
	therefore naturally stated for derivation schemes.  Computads give its
	recursive specialization, in which every skeleton is freely generated by
	the cells of the preceding dimensions.
	
	Once the $n$-cells have been freely adjoined, an equation between parallel
	$n$-cells is imposed by a \emph{coequifier}: the universal quotient which
	makes the two cells equal.  Consequently, a strict $n$-category admits a
	presentation by an $(n+1)$-computad if and only if its
	$(n-1)$-truncation is free on an $(n-1)$-computad
	(\cref{prop:higher-presentation-obstruction}).
	
	In dimension two, the same free strict $2$-category has the corresponding
	bicategorical universal property for normal pseudofunctors, meaning
	pseudofunctors whose unit constraints are identities, and for icons
	(\cref{prop:bicategory-coinserter}).
	
	To compare recursive computadic presentations with ordinary monadic
	presentations, we use a second adjunction on generator-preserving computad
	maps.  Its right adjoint records every cell together with formal
	representatives of its boundary
	(\cref{thm:free-underlying-computad-adjunction}).  This identifies higher
	computadic relation data with a distinguished class of relation data for
	the monad induced by this adjunction
	(\cref{thm:higher-computadic-monadic-comparison}).  A computadic relation
	is dimensionally typed: its two sides agree below the dimension in which
	the equation is imposed.  An arbitrary presentation for the free-category
	monad may also identify objects.  We show that these object equations may
	be imposed first, after which the remaining relation data are computadic
	(\cref{thm:computadic-monadic-comparison,thm:primeiro-os-objetos}).
	
	\paragraph{Presentation complexes and deficiency.}
	
	The categorical presentation and its topological realization can already be
	compared in dimension two.
	A \emph{$2$-computad} $C$ consists of a graph $G$, a set $C_2$ of
	generating $2$-cells, and two parallel paths in $\freecat G$ giving the
	boundary of every such generator.  Imposing each boundary pair as an
	equation gives the presented category $\Pp_1C$.  Write $\gpcomp A$ for the
	\emph{groupoid completion} of a category $A$, obtained by formally
	inverting every arrow, and write $\freegpd G$ for the groupoid freely
	generated by $G$.
	
	The \emph{groupoidal} variant allows signed paths in $\freegpd G$, rather
	than paths in $\freecat G$, as the two boundaries of a
	relation; it presents the corresponding quotient of $\freegpd G$.
	
	The \emph{presentation complex} $X_C$ is obtained from $\abs G$ by attaching
	one disk for every element of $C_2$, along the loop formed by its two
	boundary paths.  If $\PiOne(X,A)$ denotes the groupoid whose objects are
	the points of $A$ and whose arrows are endpoint-preserving homotopy classes
	of paths, let
	\[
	q_C\colon\freegpd G\longrightarrow\gpcomp(\Pp_1C)
	\quad\text{and}\quad
	\iota_C\colon\abs G\longrightarrow X_C
	\]
	denote respectively the quotient functor and the inclusion of the
	$1$-skeleton.  The comparison is the commutative square
	\begin{equation}
		\label{eq:computad-two-views}
		\begin{tikzcd}[column sep=large,row sep=large]
			\freegpd G \ar[r,two heads,"q_C"] \ar[d,"\cong"'] &
			\gpcomp(\Pp_1C) \ar[d,"\cong"]\\
			\PiOne(\abs G,G_0) \ar[r,"{(\iota_C)_*}"'] &
			\PiOne(X_C,G_0).
		\end{tikzcd}
	\end{equation}
	The fundamental groupoid is taken on the whole object set of the graph,
	rather than at one chosen basepoint.
	
	The square identifies the groupoidal quotient with the
	fundamental-groupoid calculation used in Part~II.  Brown's groupoid form of
	van Kampen's theorem computes this cellular pushout on the chosen object set
	\cite{Brown1967}.  Farjoun's theorem identifies the homotopy-invariant
	content of the same calculation: the fundamental groupoid functor preserves
	homotopy colimits \cite{Farjoun2004}.  Consequently, categorical
	quotient and cellular attachment are two realizations of the same boundary
	and relation data.  The precise comparison is
	\cref{thm:presentation-complex}.
	
	A given category or groupoid has many presentations.  Deficiency compares
	the number of generators and relations in such presentations.  For a finite
	connected groupoidal
	$2$-computad, put
	\[
	\defi(C)
	=|C_1|-|C_0|+1-|C_2|.
	\]
	The correction by the number of objects is forced by the geometry.  A
	\emph{maximal tree} is a spanning connected subgraph with no unoriented
	cycle.  It has $|C_0|-1$ arrows, used only to connect the objects.
	Collapsing the tree leaves a one-object presentation with
	$|C_1|-|C_0|+1$ loop generators and $|C_2|$ relations.  The same count is
	\[
	\defi(C)=1-\chi(X_C),
	\]
	where $\chi$ is Euler characteristic.  Thus the same normalization is
	obtained from the maximal-tree reduction and from the presentation complex.
	
	The argument also applies when the object set and the tree are infinite.  A
	connected graph is \emph{rank-finite} when only finitely many arrows lie
	outside some maximal tree; the maximal-tree normal form shows that their
	number is independent of the chosen tree, and we denote it by
	$\beta_1(G)$.  A rank-finite presentation is required in addition to have
	finitely many relations.  If a connected groupoid has finitely many objects
	and admits a finite computadic presentation, its finite deficiency is the
	classical deficiency of any isotropy group, namely the automorphism group
	at a chosen object.  For an arbitrary object set, the analogous equality
	holds for rank-finite deficiency whenever a rank-finite presentation exists
	(\cref{thm:connected-groupoid-isotropy-deficiency}).
	
	The Hopf exact sequence sharpens the Euler-characteristic calculation by
	detecting the contribution of second homology.  It yields a lower bound on
	the number of relations in terms of the first two Betti numbers of the
	presented groupoid
	(\cref{thm:homological-deficiency-bound}).  A groupoid is \emph{thin} when
	there is at most one arrow between each ordered pair of objects.  On a fixed
	connected rank-finite graph, thinness requires at least $\beta_1(G)$
	relations; choosing one relation for every arrow outside a maximal tree
	attains that bound
	(\cref{thm:cada-corda-paga-uma-relacao}).  The resulting minimality is
	relative to the chosen generating graph.
	
	For ordinary categories, groupoid completion may identify distinct arrows.
	It reflects thinness when the canonical functor from the category to its
	groupoid completion is faithful.  Alternatively, if an orientation of the
	relations is terminating and locally confluent, then the presented category
	is thin precisely when there is at most one normal path between each ordered
	pair
	(\cref{thm:faithful-localization,thm:rewriting-thinness}).
	
	\paragraph{Relations among relations.}
	
	In dimension three, the phrase ``relations among relations'' becomes
	literal.  Two composites of generating $2$-cells may have the same
	$1$-dimensional boundary without being equal.  A generating $3$-cell fills
	the resulting $2$-sphere.  In the connected groupoidal case, after choosing
	a root and a maximal tree, this datum is described by a
	\emph{crossed module}: a group homomorphism
	$\partial\colon M\to P$ together with an action of $P$ on $M$, written
	${}^pm$, such that
	\[
	\partial({}^pm)=p\partial(m)p^{-1},
	\qquad
	{}^{\partial(m)}n=mnm^{-1}.
	\]
	The classifying space of the crossed module has fundamental group
	$\operatorname{coker}\partial$ and second homotopy group
	$\ker\partial$.  Attaching a $3$-cell kills the corresponding identity among
	relations, together with all its translates under the $P$-action.
	
	A \emph{strict $2$-groupoid} is a strict $2$-category in which the
	$1$- and $2$-cells are invertible.  A \emph{cellular groupoidal
		$3$-computad} consists of a connected groupoidal $2$-computad, parallel
	$2$-cell boundaries for its $3$-generators, and a specified cellular
	attaching map $S^2\to X_C$ representing the sphere determined by each
	boundary pair.  Attaching a $3$-disk along every such map gives its
	realization $X_D$.  We compare the presented strict $2$-groupoid with the
	homotopy $2$-type of $X_D$, meaning its homotopy information through
	dimensions one and two
	(\cref{thm:crossed-module-comparison}).  If the presented strict
	$2$-groupoid is locally thin, its hom-groupoids have at most one arrow
	between each fixed pair of objects, and the second homotopy group vanishes.
	If $D$ is a finite connected cellular groupoidal $3$-computad and its
	presented strict $2$-groupoid is locally thin, the cellular count and the
	Hopf sequence give
	\[
	|D_3|\geq |D_0|-|D_1|+|D_2|-1
	+b_1(G;k)-b_2(G;k),
	\]
	where $G=\pi_1(X_D,x)$, $k$ is a field, and
	$b_i(G;k)=\dim_kH_i(G;k)$ is the $i$th Betti number
	(\cref{cor:three-cell-bound}).  It is therefore a quantitative obstruction
	to presenting a locally thin strict $2$-groupoid with too few $3$-cells.
	
	The descent computad exhibits the two arguments on the same finite
	generating data.  First apply local groupoid completion, which inverts the
	$2$-cells but not the $1$-cells.  The identity and associativity overlaps are
	then the two critical branchings of the resulting rewriting system.  Their
	confluences form a \emph{homotopy basis}: after the corresponding $3$-cells
	are imposed, every pair of parallel $2$-cells is equal.
	
	Applying groupoid completion also in dimension one gives a strict
	$2$-groupoid.  The same confluences now determine cellular $3$-cells, and
	their number attains the crossed-module and homological lower bound relative
	to each of the two fixed descent $2$-skeleta considered in
	\cref{prop:efficient-descent-presentations}.  The full subcomputad on the
	objects $1,2,3$ has no critical branching and its free track $2$-category is
	already locally thin
	(\cref{thm:intermediate-strict-descent-shape}).
	
	\subsection*{How the paper is organized}
	
	The paper has two parts.  Part~I begins with the interval-weighted
	realization of a graph and then isolates the categorical construction:
	coinserters freely adjoin cells, coequifiers impose relations, and
	derivation schemes separate one-step attachment from the recursive freeness
	of a computadic skeleton.  It also compares computadic and monadic
	presentations.  The path-length arguments in Part~I classify the nonempty
	totally ordered free categories as the finite ordinals, $\N$, $\N^{\op}$,
	and $\Z$.
	
	Part~II studies the corresponding cellular realizations.  After recalling
	the required topological background and constructing presentation complexes,
	we treat groupoid deficiency, thinness, crossed modules, and the descent
	computad.  Rewriting in its free track $2$-category and the subsequent
	groupoidal completion are treated separately.  The historical position of
	the work, including
	developments after the first 2017 preprint, is discussed only in the final
	section.
	\part{Coinserters, computads, and free higher categories}
	\section{The topological coinserter of a graph}
	\label{sec:the-big-picture}\label{panoramageral}
	
	A graph has a set of objects and a set of arrows, with specified domains and codomains: but no composition.  The free-category
	construction freely adds composites.  The topological construction
	considered below instead replaces every arrow by an interval and attaches
	its endpoints according to the domain and codomain maps.
	
	The interval-weighted formulation below is the construction used in the
	2017 preprint \cite{Nunes2017Preprint}.  After the $2$-generators are
	attached as disks, it gives the presentation complex used in Part~II.
	
	Recall the category $\mathbb G$ and the notation $G_0,G_1,s,t$ from the
	introduction.  We write $\Grph$ for the functor category
	$[\mathbb G^{\op},\Set]$.  The categorical coinserter is
	\eqref{eq:the-first-coinserter}; its value is $\freecat G$.
	
	We now make the parallel topological construction.  As above, $\Top$ denotes
	the category of compactly generated weak Hausdorff spaces.
	Let $\DTop\colon\Set\to\Top$ be the discrete-space functor, and define the
	weight
	\[
	\ITopOne\colon\mathbb G\longrightarrow\Top
	\]
	by
	\[
	\begin{aligned}
		\ITopOne(\mathsf 0)&=\{*\},&
		\ITopOne(\mathsf 1)&=[0,1],\\
		\ITopOne(d^0)(*)&=0,&
		\ITopOne(d^1)(*)&=1.
	\end{aligned}
	\]
	We regard $\mathbb G$ as a $\Top$-enriched category by giving each hom-set
	the discrete topology.  If $G\colon\mathbb G^{\op}\to\Set$, then
	$\DTop G$ denotes the objectwise composite
	\[
	\mathbb G^{\op}\xrightarrow{G}\Set\xrightarrow{\DTop}\Top .
	\]
	
	\begin{definition}[Topological realization of a graph]
		\label{def:topological-graph}\label{topographvolta}
		For a graph $G$, we define
		\begin{equation}
			\label{eq:ocoinsertertopologico}
			\FTopOne G
			\defeq \ITopOne\star\DTop G,
		\end{equation}
		where, for $W\colon\mathbb G\to\Top$ and
		$D\colon\mathbb G^{\op}\to\Top$,
		\[
		W\star D=\int^{j\in\mathbb G}W(j)\times D(j)
		\]
		denotes the $\Top$-enriched weighted colimit in the sense of
		\cite[Chapter~3]{Kelly1982}.  We call
		\eqref{eq:ocoinsertertopologico} the \emph{topological coinserter}, or the
		\emph{interval-weighted topological coinserter}, of $G$.
	\end{definition}
	
	In the 2017 preprint this interval-weighted colimit was called a
	\emph{topological isocoinserter}.  Since \emph{isocoinserter} also has the
	standard $2$-categorical meaning of a coinserter with invertible universal
	comparison $2$-cell, we use the unambiguous term
	\emph{interval-weighted topological coinserter} in the present version.
	
	That is to say, $\FTopOne G$ is obtained by adjoining one interval for every
	edge and identifying its endpoints with the prescribed vertices.  Unwinding
	the weighted colimit gives
	\begin{equation}
		\label{eq:the-interval-does-the-gluing}
		\FTopOne G\cong
		\bigl(G_0\sqcup(G_1\times[0,1])\bigr)/
		\bigl((e,0)\sim s(e),\ (e,1)\sim t(e)\bigr).
	\end{equation}
	In particular, it is the usual geometric realization of the graph as a
	$1$-dimensional CW complex.  The point of
	\cref{def:topological-graph} is that this familiar space is characterized
	by the same kind of enriched colimit that characterizes the free category.
	In this and the following topological coproducts we suppress $\DTop$ on the
	sets of objects and arrows.
	Equivalently, it is computed by the pushout
	\begin{equation}
		\label{eq:um-intervalo-para-cada-seta}
		\begin{tikzcd}[column sep=large,row sep=large]
			G_1\copow S^0 \ar[r,"{s\sqcup t}"] \ar[d,hook] &
			G_0 \ar[d] \\
			G_1\copow D^1 \ar[r] & \FTopOne G .
		\end{tikzcd}
	\end{equation}
	
	\begin{example}
		\label{ex:intervalos-antes-das-relacoes}
		Let $G$ have one object $x$ and one arrow $u\colon x\to x$.  Then
		$\freecat G$ is the one-object category whose endomorphism monoid is the
		free monoid on $u$, while $\FTopOne G$ is the circle.  For the triangular
		graph of \eqref{eq:o-triangulo-que-explica-tudo}, the realization is again a
		circle: its three arrows give a subdivision into three intervals.  The
		relation $ba=c$ does not belong to the realization of the graph.  It enters
		one dimension later, by attaching a $2$-cell, and the resulting presentation
		complex is a disk.  This is the elementary instance of the construction in
		\cref{sec:realization}.
	\end{example}
	
	For a space $X$, let $\CTopOne X$ be its \emph{path graph}: its objects
	are the points of $X$, its arrows are the continuous maps
	$[0,1]\to X$, and its source and target maps are evaluation at $0$ and
	$1$, respectively.
	
	\begin{proposition}
		\label{prop:topological-graph-adjunction}\label{oadjuntotopologico}
		The construction above is the left adjoint in an adjunction
		\[
		\FTopOne:\Grph\rightleftarrows\Top:\CTopOne.
		\]
		More precisely, there is a natural bijection
		\[
		\Top(\FTopOne G,X)\cong\Grph(G,\CTopOne X).
		\]
	\end{proposition}
	
	\begin{proof}
		A continuous map $\FTopOne G\to X$ consists of a point of $X$ for every
		object of $G$ and a path in $X$ for every arrow of $G$, whose endpoints are
		the selected points.  These are precisely the data of a graph morphism
		$G\to\CTopOne X$.  The two assignments are inverse and natural in $G$ and
		$X$.  This is also the universal property of the weighted colimit in
		\eqref{eq:ocoinsertertopologico}.
	\end{proof}
	
	\section{Conventions and notation}
	\label{sec:conventions}
	
	This section fixes terminology that otherwise changes meaning between the
	categorical and homotopical literature.
	
	\subsection{Size and categorical conventions}
	
	We work relative to a fixed Grothendieck universe.  This is only a size
	convention; class-based alternatives, including their interaction with
	ordinary and higher categorical constructions, are discussed by Levy
	\cite{Levy2018Classes}.  The symbols $\Set$,
	$\Grph$, $\Cat$, and $\Gpd$ denote the categories of small sets, graphs,
	categories, and groupoids.  The symbol $\CAT$ is used only for a
	larger ambient $2$-category.  All hom-categories appearing in a stated
	universal property are locally small.  This convention avoids treating
	$\Cat$ as an internal category in a still larger category.
	
	When a $2$-categorical universal property is asserted, we give $\Cat$ its
	usual $2$-category structure and regard $\Gpd$ as its locally full
	$2$-subcategory; the $2$-cells are natural transformations.  We use the same
	symbols for the underlying $1$-categories.
	
	Composition is written from right to left: if $f\colon x\to y$ and
	$g\colon y\to z$, then $gf$ means $g\circ f$.  Identity arrows are denoted
	by $1_x$.  A ``parallel pair'' always consists of two cells with exactly the
	same source and target.
	
	A \emph{categorical congruence} on a category $A$ is a family of
	equivalence relations on its hom-sets which is preserved by composition on
	both sides.  Its quotient has the same objects as $A$ and the corresponding
	equivalence classes as arrows.  A congruence on a groupoid is understood in
	this sense; its quotient is again a groupoid.
	
	The discrete-category functor is
	\[
	\D\colon\Set\longrightarrow\Cat.
	\]
	The fully faithful inclusion of groupoids is
	$U\colon\Gpd\to\Cat$, and its left adjoint
	\[
	\gpcomp\colon\Cat\longrightarrow\Gpd
	\]
	is \emph{groupoid completion}: $\gpcomp A$ is obtained from $A$ by formally
	inverting every arrow.  We write
	$\eta_A\colon A\to U\gpcomp A$ for the unit.  A category is
	\emph{faithfully localizable} when $\eta_A$ is faithful.  This is a property,
	not additional structure.
	
	For a set $S$ and an object $K$ in a category with coproducts, the
	\emph{copower}
	\[
	S\copow K\defeq\coprod_{s\in S}K
	\]
	records one independent copy of $K$ for each element of $S$.  In the
	cartesian categories of strict higher categories used below, it is
	canonically the product of $K$ with the discrete higher category on $S$.
	We use copower notation because it displays the indexing role of the set in
	the cell-attachment formulas.
	
	We distinguish three levels of sameness.  An
	\emph{isomorphism} has a strict inverse, an \emph{equivalence} is an
	equivalence of categories, and a \emph{biequivalence} is the corresponding
	notion for bicategories or $2$-categories.  None of these words is used as a
	synonym for another.
	
	\subsection{The two-dimensional language}
	
	We recall explicitly the part of $2$-category theory used in the paper.
	
	\begin{definition}
		A \emph{strict $2$-category} $\K$ consists of objects $A,B,\ldots$, a
		category $\K(A,B)$ for every ordered pair of objects, composition functors
		\[
		\K(B,C)\times\K(A,B)\longrightarrow\K(A,C),
		\]
		and identity objects $1_A\in\K(A,A)$, satisfying the associativity and unit
		axioms strictly.  The objects of $\K(A,B)$ are the \emph{$1$-cells}
		$A\to B$, and its arrows are the \emph{$2$-cells}.  A \emph{strict
			$2$-functor} preserves all this data and all these equations on the nose.
	\end{definition}
	
	If $\alpha\colon f\Rightarrow g$ and
	$\beta\colon g\Rightarrow h$ are $2$-cells in the same hom-category, their
	vertical composite is written $\beta\alpha$.  If
	$k\colon B\to C$ and $u\colon D\to A$ are $1$-cells, the whiskered
	$2$-cells are
	\[
	k\alpha\colon kf\Rightarrow kg,
	\qquad
	\alpha u\colon fu\Rightarrow gu.
	\]
	Horizontal composition is denoted by $\ast$.  The interchange equation
	between horizontal and vertical composition is what makes a rectangular
	pasting diagram independent of whether its rows or its columns are composed
	first.  Identity $2$-cells are denoted by $1_f$.
	
	The basic example is $\Cat$: its objects are small categories, its
	$1$-cells are functors, and its $2$-cells are natural transformations.  For
	functors $F,G\colon A\to B$, a natural transformation
	$\alpha\colon F\Rightarrow G$ is a family of arrows
	$\alpha_x\colon Fx\to Gx$ for which every arrow $u\colon x\to y$ of $A$
	gives a commutative square
	\[
	\begin{tikzcd}[column sep=large,row sep=large]
		Fx \ar[r,"Fu"] \ar[d,swap,"\alpha_x"]&
		Fy \ar[d,"\alpha_y"]\\
		Gx \ar[r,swap,"Gu"]&
		Gy.
	\end{tikzcd}
	\]
	Thus the naturality equation is
	$\alpha_yFu=Gu\alpha_x$.  Coinserters will freely adjoin components of such
	a transformation.  Icons and their higher analogues play the same role when
	only the top-dimensional cells are to be adjoined; they are defined in
	\cref{sec:higher-icons}.
	
	\subsection{Graphs, paths, and trees}
	
	Recall that a graph $G\colon\mathbb G^{\op}\to\Set$ is equivalently a
	diagram
	\[
	G_1\mathrel{\substack{\xrightarrow{s}\\[-.5ex]
			\xrightarrow[t]{}}}G_0.
	\]
	The elements of $G_0$ and $G_1$ are called the objects and arrows of $G$.
	When discussing $\abs G$ or the underlying unoriented multigraph, we also
	call them vertices and edges.  A path of length $k$ is a composable string
	$e_k\cdots e_1$ of $k$ arrows; the path of length zero at $x$ is denoted by
	$1_x$.  The free category on $G$ is denoted by $\freecat G$.
	
	We use non-reflexive graphs throughout.  A formulation with chosen reflexive
	arrows is equivalent after requiring the chosen arrows to represent identities,
	but it introduces extra structure that none of the results needs.
	
	The \emph{underlying unoriented multigraph} of $G$ has the same objects and
	arrows and associates with every arrow $e$ the unordered pair of endpoints
	$\{s(e),t(e)\}$.  A \emph{forest}
	$T\subseteq G$ means a spanning subgraph whose underlying unoriented
	multigraph has no cycle.  It is
	\emph{maximal} when its restriction to every connected component is a tree.
	Thus a maximal forest contains every object, is not generally unique, and a
	maximal forest in a finite graph has $\abs{G_0}-c(G)$ arrows, where $c(G)$ is
	the number of connected components.
	For any graph $G$, we say that $G$ is \emph{rank-finite} when
	$G_1\setminus T_1$ is finite for some maximal forest $T$.  The
	maximal-forest normal form of
	\cref{thm:maximal-forest-normal-form} shows that its cardinality is then
	independent of $T$; we put
	\[
	\beta_1(G)\defeq\abs{G_1\setminus T_1}.
	\]
	Every finite graph is rank-finite, and in that case
	\[
	\beta_1(G)=\abs{G_1}-\abs{G_0}+c(G).
	\]
	Thus $\beta_1(G)$ is a finite integer even when the sets of objects and
	tree arrows are infinite.
	
	\subsection{Thinness and connectedness}
	
	\begin{definition}
		A category $A$ is \emph{thin} when every ordered hom-set $A(x,y)$ has at most
		one element.  A $2$-category is \emph{locally thin} when every hom-category is
		thin, and it is \emph{locally groupoidal} when every hom-category is a
		groupoid.  A locally groupoidal strict $2$-category is also called a
		\emph{track $2$-category}.  A category is \emph{left and right cancellative}
		when $hf=hg$ implies $f=g$, and $fh=gh$ implies $f=g$, whenever the
		displayed composites are defined.  A $2$-category is
		\emph{locally cancellative} when every hom-category has these two
		cancellation properties.  A category is
		\emph{connected} when its groupoid completion is a connected groupoid.
	\end{definition}
	
	A thin groupoid can have one arrow $x\to y$ and one arrow $y\to x$; thinness
	does not mean that there is at most one arrow in the union of all hom-sets.
	A connected thin groupoid has exactly one arrow for every ordered pair of
	objects and is equivalent to the terminal groupoid.
	
	Throughout, a \emph{strict $2$-groupoid} is a strict $2$-category in which
	every $1$-cell is strictly invertible under $1$-cell composition and every
	$2$-cell is invertible under vertical composition.  This convention is
	stronger than being a track $2$-category: in a track $2$-category the
	$2$-cells are invertible, but the $1$-cells need not be.
	
	\subsection{Strict higher categories and walking cells}
	
	For $n\geq0$, $n\text{-}\Cat$ denotes the category of small strict globular
	$n$-categories and strict $n$-functors, with $0\text{-}\Cat=\Set$ and
	$1\text{-}\Cat=\Cat$.  We do not identify monoids with $0$-categories.
	Strict $n$-categories may equivalently be defined inductively as categories
	enriched in strict $(n-1)$-categories for the cartesian monoidal structure.
	
	For $k\geq1$, the strict $k$-category $\walking{k}$ is the
	\emph{walking $k$-cell}: the free strict $k$-category on one globular
	$k$-cell.  It contains a distinguished $k$-cell $\kappa_k$, together with
	the source and target cells forced by its iterated globular boundary, and
	no further generators.  Its boundary $\partial\walking{k}$ is the
	sub-$k$-category generated by those boundary cells, equivalently the result
	of deleting the top-dimensional generator $\kappa_k$.  Thus
	$\walking{1}=\2$ is the walking arrow category.  We put
	$\walking{0}=\1$, the singleton $0$-category with its distinguished object
	$\kappa_0$.  A $k$-cell of a strict $k$-category $A$ is equivalently a
	strict functor $\walking{k}\to A$.
	
	The locally discrete inclusion
	\[
	d_n\colon (n-1)\text{-}\Cat\longrightarrow n\text{-}\Cat
	\]
	adds only identity $n$-cells.  When a $2$-categorical universal property is
	used, the target is not the locally discrete $2$-category underlying
	$n\text{-}\Cat$; it is the iterated-icon $2$-category defined in
	\cref{sec:higher-icons}.
	
	\section{Coinserters and relations}
	\label{sec:coinserters}
	
	We recall only the $2$-colimits used in the sequel.  Our convention is the
	strict, Cat-enriched universal property.  Replacing the displayed
	isomorphisms of categories by equivalences gives the corresponding
	bicolimits.
	
	\begin{definition}[Coinserter]
		Let $f,g\colon A\rightrightarrows B$ be parallel $1$-cells in a $2$-category
		$\K$.  A \emph{coinserter} of $(f,g)$ consists of a $1$-cell
		$q\colon B\to Q$ and a $2$-cell
		\[
		\theta\colon qf\Rightarrow qg
		\]
		such that, for every object $X$, composition with $(q,\theta)$ induces an
		isomorphism of categories
		\begin{equation}
			\label{eq:coinserter-universal-property}
			\K(Q,X)\ \cong\
			\bigl\{(h,\alpha)\mid h\colon B\to X,
			\ \alpha\colon hf\Rightarrow hg\bigr\}.
		\end{equation}
		On the right, a morphism $(h,\alpha)\to(k,\beta)$ is a $2$-cell
		$\gamma\colon h\Rightarrow k$ satisfying
		$(\gamma g)\alpha=\beta(\gamma f)$.
	\end{definition}
	
	The comparison $2$-cell is best kept visible:
	\begin{equation}
		\label{eq:a-celula-universal-do-coinserter}
		\begin{tikzcd}[column sep=huge]
			A
			\ar[r,bend left=38,"qf",""{name=upper,below}]
			\ar[r,bend right=38,"qg"',""{name=lower,above}]&
			Q .
			\arrow[Rightarrow,from=upper,to=lower,"\theta" description]
		\end{tikzcd}
	\end{equation}
	For a cone $(h,\alpha)$ with vertex $X$, the corresponding picture is the
	same one with $(Q,q,\theta)$ replaced by $(X,h,\alpha)$.  The equation
	$(\gamma g)\alpha=\beta(\gamma f)$ says exactly that the two resulting
	pasted $2$-cells agree.
	
	The right-hand side of \eqref{eq:coinserter-universal-property} is the
	inserter category of the two precomposition functors together with its
	displayed comparison cell.  Thus the definition includes the action on
	$2$-cells; it is stronger than a bijection on objects.
	
	\begin{definition}[Isocoinserter and coequifier]
		An \emph{isocoinserter} has the same universal property as a coinserter, with
		$\alpha$ required to be invertible.  Given parallel $2$-cells
		$\alpha,\beta\colon f\Rightarrow g$, a \emph{coequifier} is a $1$-cell
		$q\colon B\to Q$, satisfying $q\alpha=q\beta$, such that for every $X$,
		precomposition with $q$ induces an isomorphism
		\[
		\K(Q,X)\cong
		\{\,h\colon B\to X\mid h\alpha=h\beta\,\},
		\]
		where the right-hand side is the full subcategory of $\K(B,X)$ on the
		displayed objects.
	\end{definition}
	
	Coinserters freely add a comparison cell; isocoinserters freely add an
	invertible comparison cell; coequifiers impose an equation between already
	existing cells.  The same order reappears in computadic presentations:
	an $n$-computad adjoins $n$-cells, while an $(n+1)$-computad presents an
	$n$-category by imposing equations between them.
	
	\begin{remark}
		Coinserters are weighted $2$-colimits in the sense of
		\cite{Street1976,Kelly1982}; for the elementary named $2$-categorical
		limits and colimits used here, see also \cite{Kelly1989}.  No general theory
		of weights is needed below: all universal properties are reduced explicitly
		to assignments of generators and cells.  Coinserters and coequifiers are
		also the elementary
		colimits from which Lack constructs lax codescent objects
		\cite[Proposition~2.1]{Lack2002Codescent}.  Here they occur in the cellular
		order used throughout the paper: the coinserter adjoins the generating
		cells, and the coequifier imposes equations between cells already
		constructed.
	\end{remark}
	
	\subsection{The free category as a coinserter}
	
	Let $G=(G_1\rightrightarrows G_0)$ be a graph.  Regard the source and target
	maps as functors between discrete categories:
	\begin{equation}
		\label{eq:graph-discrete-pair}
		\begin{tikzcd}[column sep=large]
			\D G_1
			\ar[r,shift left=.7ex,"\D s"]
			\ar[r,shift right=.7ex,swap,"\D t"]&
			\D G_0 .
		\end{tikzcd}
	\end{equation}
	
	\begin{theorem}[Free-category coinserter]
		\label{thm:free-category-coinserter}
		The coinserter of \eqref{eq:graph-discrete-pair} in $\Cat$ is the free category
		$\freecat G$.  Its universal $2$-cell has component at
		$e\in G_1$ equal to the length-one path
		$e\colon s(e)\to t(e)$.
	\end{theorem}
	
	Thus the universal cell may be drawn as
	\[
	\begin{tikzcd}[column sep=huge]
		\D G_1
		\ar[r,bend left=38,"q\D s",""{name=upper,below}]
		\ar[r,bend right=38,"q\D t"',""{name=lower,above}]&
		\freecat G .
		\arrow[Rightarrow,from=upper,to=lower,"\theta" description]
	\end{tikzcd}
	\qquad
	\theta_e=e\quad(e\in G_1).
	\]
	
	\begin{proof}
		Let $A$ be a category.  A functor $h\colon\D G_0\to A$ selects an object
		$h(x)$ for each object $x$ of $G$.  Since $\D G_1$ is discrete, a natural
		transformation
		\[
		\alpha\colon h\D s\Rightarrow h\D t
		\]
		is precisely a family of arrows
		$\alpha_e\colon h(s(e))\to h(t(e))$, one for every arrow of $G$.  By the path
		description of $\freecat G$, these data extend uniquely to a functor
		$\bar h\colon\freecat G\to A$, with
		\[
		\bar h(e_k\cdots e_1)=\alpha_{e_k}\cdots\alpha_{e_1}.
		\]
		Conversely, restriction of a functor along the objects and generating arrows
		recovers $(h,\alpha)$.
		
		A natural transformation $\gamma\colon\bar h\Rightarrow\bar k$ is determined
		by its components at the objects, and naturality on an arbitrary path follows
		by composition from naturality on each generating arrow.  Hence the two
		constructions give an isomorphism of categories, natural in $A$, exactly as
		required by \eqref{eq:coinserter-universal-property}.
	\end{proof}
	
	\begin{example}
		\label{ex:o-triangulo-como-apresentacao}
		For the graph in \eqref{eq:o-triangulo-que-explica-tudo}, the universal
		transformation of \cref{thm:free-category-coinserter} has components
		$a,b,c$.  The coequifier of the parallel arrows $ba,c\colon x\to z$
		imposes the single relation $ba=c$.  The resulting category has precisely
		one arrow between comparable objects
		\[
		x<y<z,
		\]
		and is therefore isomorphic to the ordinal category $[2]$.  The categorical
		presentation and the disk described in
		\cref{ex:intervalos-antes-das-relacoes} are the algebraic and topological
		forms of the same attachment.
	\end{example}
	
	A closely related two-dimensional universal property for the free category
	with finite products on a multigraph was given by Walters
	\cite{Walters1989}.
	
	\begin{theorem}[Free-groupoid isocoinserter]
		\label{thm:free-groupoid-coinserter}
		For a graph $G$, the following objects have the same universal property:
		\begin{enumerate}
			\item the coinserter of \eqref{eq:graph-discrete-pair} in $\Gpd$;
			\item the isocoinserter of \eqref{eq:graph-discrete-pair} in $\Cat$;
			\item the free groupoid $\freegpd G$ on $G$.
		\end{enumerate}
		In particular, there is a canonical isomorphism
		$\freegpd G\cong\gpcomp\freecat G$.
	\end{theorem}
	
	\begin{proof}
		For a groupoid $H$, every component of a natural transformation between
		functors into $H$ is invertible.  Thus the coinserter data in $\Gpd$ are an
		assignment of objects to the objects of $G$ and invertible arrows to its
		arrows.  This is the universal property of $\freegpd G$.  For an arbitrary
		category $A$, an isocone in $\Cat$ is given by the same assignments, so it
		extends uniquely to a functor $\freegpd G\to A$.
		The naturality equations at the arrows of $G$ identify morphisms of isocones with natural
		transformations between the extensions.  This proves the full
		hom-category-valued isocoinserter property.  Finally, maps
		$\gpcomp\freecat G\to H$ correspond under $\gpcomp\dashv U$ to maps
		$\freecat G\to UH$, which have the same generator data.
	\end{proof}
	
	\begin{remark}
		The word \emph{coinverter} is sometimes used for the universal operation that
		inverts the comparison cell of a previously formed coinserter.  We use
		\emph{isocoinserter} for the one-step construction above, so no unnamed
		intermediate object is needed.  Coinverters and their relation with
		categories of fractions are treated in \cite{KellyLackWalters1993}.
	\end{remark}
	
	\subsection{Coequifiers as relations}
	
	Suppose $B$ is a category and $u_i,v_i\colon x_i\to y_i$ are parallel arrows,
	indexed by a set $R$.  Let $Q$ be the quotient of $B$ by the least categorical
	congruence containing $u_i=v_i$ for every $i$.  This familiar construction is
	the one-dimensional shadow of a coequifier.
	
	\begin{proposition}
		\label{prop:quotient-coequifier}
		Let $\D R$ be the discrete category on the indexing set, and define functors
		$x,y\colon\D R\rightrightarrows B$ by
		$x(i)=x_i$ and $y(i)=y_i$.  The families $u=(u_i)$ and $v=(v_i)$ are parallel
		natural transformations $x\Rightarrow y$.  Their coequifier in $\Cat$ is the
		quotient functor $q\colon B\to Q$.  Equivalently, precomposition with $q$
		identifies $\Cat(Q,A)$ with the full subcategory of $\Cat(B,A)$ on the
		functors $h$ satisfying $hu_i=hv_i$ for every $i\in R$.
	\end{proposition}
	
	\begin{proof}
		Because $\D R$ is discrete, $u$ and $v$ are precisely the displayed
		families of parallel arrows.  The equality $hu=hv$ of natural transformations
		holds exactly when $hu_i=hv_i$ for every $i$.  Factoring $h$ through the
		quotient by the least categorical congruence generated by these equations is
		therefore the coequifier universal property, on both functors and natural
		transformations.
	\end{proof}
	
	\section{Free categories and free groupoids}
	\label{sec:graphs}
	
	The coinserter theorems determine free categories by universal properties,
	while the deficiency theory of Part~II requires concrete representatives.
	We therefore record path, reduced-word, and maximal-forest normal forms.  A
	reader interested first in the formal higher-dimensional theorem may pass to
	\cref{sec:presentations} and return here at the beginning of Part~II.
	
	\subsection{Paths in a free category}
	
	Every arrow of $\freecat G$ is a path
	$e_k\cdots e_1$, and equality is literal equality of strings after deleting
	identity paths.  In particular, path length is well defined and additive.
	
	\begin{proposition}
		\label{prop:free-category-cancellation}
		Every free category satisfies left and right cancellation.  Moreover,
		$\freecat G$ is thin if and only if, for each ordered pair of objects
		$(x,y)$, the graph $G$ has at most one path from $x$ to $y$, including
		the length-zero path when $x=y$.
	\end{proposition}
	
	\begin{proof}
		If $rp=rq$ as paths, deletion of the common terminal string $r$ gives $p=q$;
		the other cancellation law is analogous.  The arrows
		$x\to y$ in $\freecat G$ are exactly the paths from $x$ to $y$, so the second
		statement is the definition of thinness applied to the path normal form.
	\end{proof}
	
	Cancellation is a property of the word calculus of a category; by itself it
	does not imply that the unit into groupoid completion is faithful.
	
	\subsection{Totally ordered free categories}
	
	\begin{definition}
		A \emph{totally ordered category} is a thin category $A$ for which
		\[
		x\leq y\quad\Longleftrightarrow\quad A(x,y)\neq\varnothing
		\]
		is a total order on the objects.  Equivalently, it is the category associated
		with a totally ordered set: there is one arrow $x\to y$ when $x\leq y$, and no
		arrow otherwise.
	\end{definition}
	
	The path length in a free category makes total order much more restrictive
	than thinness alone.
	
	\begin{lemma}[Finite intervals]
		\label{lem:ordered-free-category-intervals}
		Let $\freecat G$ be totally ordered.
		
		\begin{enumerate}
			\item For an object $x$ and an integer $m\geq0$, there is at most one arrow of
			length $m$ with domain $x$, and at most one arrow of length $m$ with codomain
			$x$.
			\item Every object which is not greatest has a unique immediate successor,
			joined to it by an arrow of length one.  Dually, every object which is not
			least has a unique immediate predecessor.
			\item If $x\leq y$ and the arrow $x\to y$ has length $m$, then the interval
			\[
			[x,y]=\{z\mid x\leq z\leq y\}
			\]
			has exactly $m+1$ objects.
		\end{enumerate}
	\end{lemma}
	
	\begin{proof}
		Suppose that $a\colon x\to y$ and $b\colon x\to z$ both have length $m$.
		Totality allows us, after interchanging $a$ and $b$, to assume $y\leq z$.
		Let $c\colon y\to z$ be the unique arrow.  Thinness gives $ca=b$, and
		additivity of path length gives
		\[
		\ell(c)+m=\ell(ca)=\ell(b)=m.
		\]
		Thus $\ell(c)=0$, so $c$ is an identity, $y=z$, and $a=b$.  Applying the same
		argument in
		\[
		(\freecat G)^{\op}\cong\freecat(G^{\op})
		\]
		proves the assertion with fixed codomain.
		
		Now suppose that $x$ is not greatest.  Choose $y>x$ and write the unique arrow
		$x\to y$ as a nonempty path.  Its first arrow has domain $x$.  By the first
		part, it is the only length-one arrow with domain $x$; call its codomain
		$x^+$.  Every path from $x$ to an object strictly above $x$ begins with this
		arrow, again by uniqueness.  Hence $x^+$ is the least object strictly above
		$x$, and is therefore the unique immediate successor.  The assertion about
		predecessors follows by passage to the opposite category.
		
		It remains to identify an arbitrary interval.  Write the unique arrow
		$p\colon x\to y$ as
		\[
		x=x_0\xrightarrow{e_1}x_1\xrightarrow{e_2}\cdots
		\xrightarrow{e_m}x_m=y.
		\]
		The objects are strictly increasing: a loop arrow would give a nonidentity
		endomorphism parallel to an identity.  Let $z\in[x,y]$, and denote the unique
		arrows $x\to z$ and $z\to y$ by $a$ and $b$.  Since $ba=p$, equality of paths
		in a free category says that this factorization is obtained by cutting the
		displayed arrow word at one place.  Thus $z=x_k$ for a unique
		$0\leq k\leq m$.  Conversely, every $x_k$ lies in the interval.  Therefore
		$[x,y]=\{x_0,\ldots,x_m\}$.
	\end{proof}
	
	\begin{theorem}[Classification of totally ordered free categories]
		\label{thm:totally-ordered-free-categories}
		Let $G$ be a graph such that $\freecat G$ is nonempty and totally ordered.
		Then $\freecat G$ is isomorphic, as an ordered category, to exactly one of
		\[
		[m]=\{0<1<\cdots<m\}\quad(m\geq0),\qquad
		\N,\qquad \N^{\op},\qquad\text{or}\qquad\Z .
		\]
		Conversely, every ordered category in this list is a free category.
	\end{theorem}
	
	\begin{proof}
		There are four cases, according to the existence of a least or a greatest
		object.
		
		Suppose first that both exist, and denote them by $\bot$ and $\top$.  If the
		arrow $\bot\to\top$ has length $m$, then
		\cref{lem:ordered-free-category-intervals} identifies the whole object set
		$[\bot,\top]$ with a chain of $m+1$ objects.  Since a totally ordered category
		has its arrows determined by its order, this gives
		$\freecat G\cong[m]$.
		
		Suppose next that there is a least object but no greatest object.  Put
		$x_0=\bot$ and let $x_{n+1}$ be the immediate successor of $x_n$.  The
		successor exists at every stage, and the composite $x_0\to x_n$ has length
		$n$.  If $y$ is any object and the arrow $\bot\to y$ has length $m$, the
		finite-interval description shows that its successive objects are
		$x_0,x_1,\ldots,x_m=y$.  Hence $n\mapsto x_n$ is an order isomorphism
		$\N\cong\freecat G$.  The case of a greatest object and no least object is
		dual and gives $\N^{\op}$.
		
		Finally, suppose that neither endpoint exists.  Choose an object $x_0$ and
		iterate its unique successors and predecessors:
		\[
		\cdots<x_{-2}<x_{-1}<x_0<x_1<x_2<\cdots .
		\]
		For $y\geq x_0$, the finite interval $[x_0,y]$ shows that $y=x_m$, where $m$
		is the length of $x_0\to y$.  For $y\leq x_0$, the dual argument gives
		$y=x_{-m}$, where $m$ is the length of $y\to x_0$.  Thus
		$n\mapsto x_n$ is an order isomorphism $\Z\cong\freecat G$.
		
		For the converse, take the graph having one arrow between each pair of
		successive elements of the indicated order.  Its finite paths are in
		bijection with comparable pairs, so its free category is the corresponding
		ordered category.  Finally, the finite orders are distinguished by their
		cardinalities, while $\N$, $\N^{\op}$, and $\Z$ are distinguished by the
		existence of least and greatest objects.  Hence the order type in the
		statement is unique.
	\end{proof}
	
	\subsection{Reduced arrow words}
	
	For every arrow $e$ of $G$, adjoin a formal inverse $e^{-1}$ with
	$s(e^{-1})=t(e)$ and $t(e^{-1})=s(e)$.  A \emph{signed-arrow word} is a
	composable word in these letters.  It is \emph{reduced} if it has no adjacent
	subword $e^{-1}e$ or $ee^{-1}$.
	
	\begin{theorem}[Reduced-word normal form]
		\label{thm:reduced-word-normal-form}
		Every arrow of $\freegpd G$ has a unique reduced signed-arrow representative.
		Composition is concatenation followed by cancellation of the unique maximal
		inverse overlap.
	\end{theorem}
	
	\begin{proof}
		Let $W(G)$ have object set $G_0$, with arrows the signed-arrow words modulo
		the congruence generated by
		\[
		ee^{-1}\sim 1_{t(e)},\qquad e^{-1}e\sim 1_{s(e)}.
		\]
		Concatenation and reversal of signed words make $W(G)$ a groupoid.  An
		assignment of the arrows of $G$ to arrows of a groupoid extends uniquely to
		a functor from $W(G)$; hence $W(G)$ has the universal property of
		$\freegpd G$, and we identify the two.
		
		Repeated deletion of adjacent inverse pairs terminates because it strictly
		shortens a word.  For uniqueness, define the stack reduction by reading a
		word from right to left: cancel the top letter when it is inverse to the next
		input letter, and otherwise push the new letter.  This reduction is unchanged
		by either insertion or deletion of an adjacent inverse pair: the two
		successive stack operations are mutually inverse, whether the first operation
		pushes a letter or cancels the letter already on top.  Consequently it
		is constant on the congruence classes defining $W(G)$, and a reduced word is
		returned unchanged.  Every congruence class therefore contains exactly one
		reduced word.  If two reduced words are concatenated, all cancellations occur
		in their unique maximal inverse overlap, which gives the stated composition
		rule.
	\end{proof}
	
	\begin{remark}[Why there is no natural integer length]
		\label{rem:no-natural-groupoid-length}
		The reduced length $\ell_G\colon\Mor(\freegpd G)\to\N$ is useful for a fixed
		graph, but it is not functorial under arbitrary graph maps.  A graph map may
		identify two arrows, turning a reduced word into one with cancellation.  We
		use reduced words only objectwise, never as a natural transformation of
		monads.
	\end{remark}
	
	\subsection{Maximal forests}
	
	Let $T\subseteq G$ be a maximal forest.  Choose a root $r_j$ in each connected
	component and, for every object $x$ in that component, let
	$\tau_x\colon r_j\to x$ be the unique reduced signed-arrow word in $T$.
	For an arrow $e\colon x\to y$ not in $T$, define the loop
	\begin{equation}
		\label{eq:chord-loop}
		\lambda_e=\tau_y^{-1}e\tau_x\colon r_j\longrightarrow r_j.
	\end{equation}
	The three factors are visible in the triangle
	\begin{equation}
		\label{eq:o-laco-de-uma-corda}
		\begin{tikzcd}[column sep=large,row sep=large]
			& x \ar[dr,"e"]&\\
			r_j \ar[ur,"\tau_x"]&&
			y \ar[ll,bend left=18,"\tau_y^{-1}"]
		\end{tikzcd}
	\end{equation}
	We call an arrow of $G_1\setminus T_1$ a \emph{chord} of the maximal forest
	$T$.
	
	\begin{theorem}[Maximal-forest normal form]
		\label{thm:maximal-forest-normal-form}
		For each component $G^{(j)}$, the isotropy group
		$\freegpd G(r_j,r_j)$ is the free group on the set
		$G^{(j)}_1\setminus T^{(j)}_1$, with basis given by the loops
		\eqref{eq:chord-loop}.  Consequently,
		\[
		\freegpd G\simeq
		\coprod_j \Sigma F\bigl(G^{(j)}_1\setminus T^{(j)}_1\bigr),
		\]
		where $F(S)$ is the free group on $S$ and $\Sigma H$ denotes the one-object
		groupoid with automorphism group $H$.
	\end{theorem}
	
	\begin{proof}
		For an original chord $e\colon x\to y$, the definition gives
		\[
		e=\tau_y\lambda_e\tau_x^{-1},
		\qquad
		e^{-1}=\tau_x\lambda_e^{-1}\tau_y^{-1}.
		\]
		If $a\colon u\to v$ is a signed tree arrow, then
		$\tau_v^{-1}a\tau_u$ is the identity.  Indeed, the reduced form of
		$a\tau_u$ is $\tau_v$, since both are paths in the tree from the root to
		$v$.  Thus every signed
		letter $a\colon u\to v$ has an expression
		\[
		a=\tau_v\,\omega(a)\,\tau_u^{-1},
		\]
		where $\omega(a)$ is $1$ for a tree arrow, $\lambda_e$ for a positively
		oriented chord $e$, and $\lambda_e^{-1}$ for its inverse.  Substitution into
		a signed path from $x$ to $y$ makes the intervening tree paths telescope and
		gives $\tau_yw\tau_x^{-1}$, with $w$ a word in the chord loops.  Conversely,
		expanding a reduced word in the chord loops gives a signed-arrow word in $G$.
		
		To see directly that no relation has been introduced, collapse every tree
		arrow to an identity and send each chord $e$ to a distinct generator of the
		free group $F(G^{(j)}_1\setminus T^{(j)}_1)$.  This gives a functor from the
		free groupoid of the component to the corresponding one-object groupoid, and
		it sends $\lambda_e$ to the generator indexed by $e$.  The homomorphism in
		the other direction sends that generator to $\lambda_e$.  The preceding
		rewriting shows that these homomorphisms are inverse on the isotropy group.
		The chosen arrows $\tau_x$
		then exhibit the full connected groupoid as equivalent to the one-object
		isotropy groupoid.  Taking the coproduct over components proves the result.
	\end{proof}
	
	For finite $G$, the rank of the displayed free group is
	$\abs{G_1}-\abs{G_0}+c(G)=\beta_1(G)$.  The formulation by chords also works
	without cardinal subtraction when $G$ is rank-finite.
	
	By a subgroupoid $H\subseteq K$ we mean a subcategory closed under the
	inverses in $K$; it need not contain every object of $K$.
	
	\begin{corollary}
		\label{cor:subgroupoid-free}
		Every subgroupoid of a free groupoid is a free groupoid.
	\end{corollary}
	
	\begin{proof}
		Each connected component of a subgroupoid is equivalent to one of its
		isotropy groups.  That isotropy group is a subgroup of an isotropy group of
		the ambient free groupoid, hence a subgroup of a free group.  It is free by
		the Nielsen--Schreier theorem
		\cite[Chapter~I, Proposition~3.10]{LyndonSchupp1977}.
		
		To pass from this isotropy calculation to an isomorphism of groupoids, fix a
		root $r$ in the component and choose an arrow
		$\tau_x\colon r\to x$ for every object $x$, with $\tau_r=1_r$.  Take the
		star-shaped tree with one generating arrow from $r$ to each $x\neq r$, and
		adjoin one loop at $r$ for every element of a free basis of the isotropy
		group.  Send these generators to the corresponding $\tau_x$ and basis
		elements.  Every arrow $f\colon x\to y$ has the unique normal form
		\[
		f=\tau_y h\tau_x^{-1},
		\qquad h\in H(r,r),
		\]
		so the induced functor from the free groupoid is an isomorphism.  Repeating
		the construction in every component proves the result.
	\end{proof}
	
	\subsection{Faithful groupoid completion}
	
	The next proposition isolates the faithfulness hypothesis used later.
	
	\begin{proposition}
		\label{prop:localization-reflects-thin-under-faithfulness}
		If $A$ is faithfully localizable and $\gpcomp A$ is thin, then $A$ is thin.
	\end{proposition}
	
	\begin{proof}
		For parallel arrows $f,g\colon x\to y$ in $A$, thinness gives
		$\eta_A(f)=\eta_A(g)$ in $\gpcomp A$.  Faithfulness of $\eta_A$ gives $f=g$.
	\end{proof}
	
	\begin{remark}
		Faithful localizability cannot be replaced by two-sided cancellation.  For a
		one-object category, the former says that its endomorphism monoid embeds in a
		group, while the latter says only that the monoid is left and right
		cancellative.  Johnstone gives categorical criteria for groupoid
		embeddability \cite{Johnstone2008}.  Explicit cancellative monoids that do not
		embed in groups are exhibited in \cite{EdwardesHeath2025}.  Thus faithfulness
		is the hypothesis that permits thinness to be reflected from the completion.
	\end{remark}
	
	\section{Monadic and computadic presentations}
	\label{sec:presentations}
	
	We consider presentations of individual algebras for a fixed monad.  This
	is distinct from a presentation of the monad itself, as in the theory of
	finitary enriched monads of Kelly and Power, and from the comparison between
	monads and theories
	\cite{KellyPower1993,BourkeGarner2019}.
	
	We first describe presentations for an arbitrary monad by parallel arrows
	into a free algebra.  We then identify the data supplied by a $2$-computad
	for the free-category monad.  The comparison records the dimension of each
	relation and makes explicit the requirement that the boundary maps preserve
	objects.  We finally show that this requirement causes no loss: an arbitrary
	free-category-monad presentation can first be quotiented on objects and then
	rewritten as an object-agreeing computadic presentation.
	
	\subsection{Presentations for a monad}
	
	\begin{definition}
		A \emph{monad} on a category $\A$ is an endofunctor
		$T\colon\A\to\A$ equipped with natural transformations
		\[
		\eta\colon1_{\A}\Rightarrow T,
		\qquad \mu\colon TT\Rightarrow T,
		\]
		satisfying the unit and associativity equations
		\[
		\mu\,T\eta=1_T=\mu\,\eta T,
		\qquad
		\mu\,T\mu=\mu\,\mu T.
		\]
		A \emph{$T$-algebra} is an object $A$ together with an arrow
		$a\colon TA\to A$ such that
		$a\eta_A=1_A$ and $aTa=a\mu_A$.  An algebra homomorphism
		$h\colon(A,a)\to(B,b)$ is an arrow $h\colon A\to B$ satisfying
		$ha=bTh$.  The resulting category is denoted by $\A^T$.
	\end{definition}
	
	The free $T$-algebra on $X$ is $(TX,\mu_X)$.  It gives the
	Eilenberg--Moore adjunction
	\[
	F^T\colon\A\rightleftarrows\A^T\mathbin{:}U^T
	\]
	whose induced monad is $T$.  An arrow $r\colon R\to TX$ in $\A$ is a
	\emph{Kleisli arrow} $R\rightsquigarrow X$ for $T$.  It corresponds to the
	algebra homomorphism
	\begin{equation}
		\label{eq:kleisli-extension}
		r^\sharp=\mu_X\,T r\colon TR\longrightarrow TX.
	\end{equation}
	
	\begin{definition}[Monadic presentation]
		\label{def:monadic-presentation}
		A \emph{$T$-relation datum} is a quadruple
		$(R,X,r_0,r_1)$ with parallel arrows
		\[
		R\mathrel{\substack{\xrightarrow{r_0}\\[-.5ex]
				\xrightarrow[r_1]{}}}TX.
		\]
		If the coequalizer of
		$r_0^\sharp,r_1^\sharp\colon F^T R\rightrightarrows F^T X$ exists in
		$\A^T$, it is the \emph{algebra presented by the datum}.  The object $X$ is
		the object of generators and $R$ is the object indexing relations.
	\end{definition}
	
	Relation data therefore contain two parallel maps, one for each side of an
	equation.  They form a category $\operatorname{Rel}_T$: morphisms are pairs
	$(u,v)\colon(R,X,r_0,r_1)\to(R',X',r'_0,r'_1)$ satisfying
	$Tv\,r_i=r'_i u$ for $i=0,1$, with identities and composition defined
	componentwise.  On the full subcategory of $\operatorname{Rel}_T$ on which
	the indicated coequalizers exist, choosing those coequalizers makes the
	construction in
	\cref{def:monadic-presentation} a functor to $\A^T$.  Without chosen
	coequalizers it is canonical up to unique isomorphism.
	
	A parallel pair $f,g\colon A\rightrightarrows B$ is
	\emph{reflexive} when it has a common section: an arrow
	$s\colon B\to A$ such that $fs=1_B=gs$.
	
	\begin{theorem}[Canonical reflexive presentation]
		\label{thm:canonical-monadic-presentation}
		Every $T$-algebra $a\colon TA\to A$ is the coequalizer in $\A^T$ of the
		reflexive pair
		\begin{equation}
			\label{eq:canonical-em-coequalizer}
			F^T(TA)\mathrel{\substack{\xrightarrow{\mu_A}\\[-.5ex]
					\xrightarrow[T a]{}}}F^T A\xrightarrow{a}(A,a).
		\end{equation}
		Equivalently, it is presented by the Kleisli pair with relation object $TA$,
		generator object $A$, and maps
		\[
		1_{TA},\ \eta_A a\colon TA\rightrightarrows TA.
		\]
	\end{theorem}
	
	\begin{proof}
		Monad associativity says that
		$\mu_A\colon F^T(TA)\to F^TA$ is an algebra homomorphism, and naturality of
		$\mu$ at $a$ says the same for
		$Ta\colon F^T(TA)\to F^TA$.  The algebra associativity axiom makes
		$a\colon F^TA\to(A,a)$ an algebra homomorphism.  Finally, naturality of
		$\mu$ at $\eta_A$ shows that
		$T\eta_A\colon F^TA\to F^T(TA)$ is an algebra homomorphism.
		
		Under \eqref{eq:kleisli-extension}, $1_{TA}$ corresponds to $\mu_A$ and
		$\eta_Aa$ corresponds to
		\[
		\mu_A T(\eta_Aa)=\mu_A T\eta_A\,T a=T a.
		\]
		The algebra map $a$ coequalizes $\mu_A$ and $T a$ by the algebra axiom.
		The pair is reflexive, with common section
		$T\eta_A\colon TA\to TTA$, since
		$\mu_A T\eta_A=1$ and $T a\,T\eta_A=T(a\eta_A)=1$.
		
		If an algebra homomorphism $h\colon F^T A\to(B,b)$ coequalizes the pair, put
		$\bar h=h\eta_A\colon A\to B$.  Then
		$\bar h a=h$: indeed, naturality of $\eta$ and the coequalizing equation give
		\[
		h\eta_Aa=h\,T a\,\eta_{TA}=h\mu_A\eta_{TA}=h.
		\]
		Moreover, $bT\bar h=\bar h a$ follows from the fact that $h$ is an algebra
		map.  Uniqueness among algebra maps satisfying $\bar h a=h$ is immediate:
		if $k$ is another, then
		$k=ka\eta_A=h\eta_A=\bar h$.  This is the required coequalizer property.
	\end{proof}
	
	\begin{remark}
		Numerical finiteness enters only after generator and relation objects have
		been equipped with specified size functions.  In the computadic cases below,
		these sizes are cardinalities or ranks of the relevant cell sets.  No
		finiteness hypothesis is implicit in
		\cref{thm:canonical-monadic-presentation}.
	\end{remark}
	
	\subsection{Computads in dimension two}
	
	Let $G$ be a graph.  The source and target of a prospective generating
	$2$-cell are arbitrary parallel arrows in $\freecat G$; in particular, they
	need not be generating arrows of $G$.
	We write
	\[
	U_{\Grph}\colon\Cat\longrightarrow\Grph
	\]
	for the underlying-graph functor; thus every arrow of a category, including
	an identity, becomes an arrow of its underlying graph.
	
	We use Street's two-dimensional computad, with the two boundary maps written
	explicitly \cite{Street1976}.
	
	\begin{definition}[$2$-computad]
		\label{def:two-computad}
		A \emph{$2$-computad} $C$ consists of a graph
		$C_{\leq1}=(C_1\rightrightarrows C_0)$, a set $C_2$, and two functions
		\[
		\partial^-,\partial^+\colon C_2\longrightarrow
		\Mor(\freecat C_{\leq1})
		\]
		whose values at each $\alpha\in C_2$ are parallel paths.  We write
		\[
		\alpha\colon\partial^-\alpha\Rightarrow\partial^+\alpha.
		\]
		The category presented by $C$ is
		\[
		\Pp_1 C\defeq
		\freecat C_{\leq1}/
		\bigl(\partial^-\alpha=\partial^+\alpha\mid\alpha\in C_2\bigr).
		\]
	\end{definition}
	
	A \emph{map of $2$-computads} $f\colon C\to D$ consists of a function on
	objects, a path $f(e)$ in $\freecat D_{\leq1}$ with the transported
	endpoints for every generating arrow $e\in C_1$, and, for every
	$\alpha\in C_2$, a formal
	$2$-cell in the free strict $2$-category on $D$ from
	$f(\partial^-\alpha)$ to $f(\partial^+\alpha)$.  Here ``formal'' means
	generated by the elements of $D_2$ using identities and horizontal and
	vertical composition, subject only to the strict $2$-category axioms; the
	construction is given formally in
	\cref{def:ascelulassaolivres}.
	Equivalently, its one-dimensional part is a Kleisli graph map
	\[
	C_{\leq1}\longrightarrow U_{\Grph}\freecat D_{\leq1}.
	\]
	It is \emph{generator-preserving} when its one-dimensional part is a graph
	map $f_{\leq1}\colon C_{\leq1}\to D_{\leq1}$ and there is a function
	$f_2\colon C_2\to D_2$ satisfying
	\[
	(\freecat f_{\leq1})(\partial^\epsilon\alpha)
	=\partial^\epsilon(f_2\alpha)
	\qquad(\epsilon\in\{-,+\}).
	\]
	The associated general map sends $\alpha$ to the target generator
	$f_2(\alpha)$.  We write $2\text{-}\Comp$ for general maps and
	$2\text{-}\Comp^{\mathrm{gen}}$ for generator-preserving maps.
	
	We say that $C$ \emph{presents} a category $A$ when an isomorphism
	$\Pp_1C\cong A$ is specified or clear from context.  We say ``presents up to
	equivalence'' when only an equivalence is asserted.
	
	Let $\2_{\Grph}$ be the graph with two objects and one arrow.  The boundary
	functions assemble into graph maps
	\begin{equation}
		\label{eq:computad-kleisli-pair}
		C_2\copow\2_{\Grph}
		\mathrel{\substack{\xrightarrow{\partial^-}\\[-.5ex]
				\xrightarrow[\partial^+]{}}}
		U_{\Grph}\freecat C_{\leq1},
	\end{equation}
	where the two maps agree on objects.
	
	\begin{proposition}
		\label{prop:computad-category-coequalizer}
		The category $\Pp_1 C$ is the coequalizer in $\Cat$ of
		\begin{equation}
			\label{eq:computad-category-pair}
			C_2\copow\2
			\mathrel{\substack{\xrightarrow{\partial^-}\\[-.5ex]
					\xrightarrow[\partial^+]{}}}
			\freecat C_{\leq1}.
		\end{equation}
	\end{proposition}
	
	\begin{proof}
		A functor $H\colon\freecat C_{\leq1}\to A$ coequalizes the two functors in
		\eqref{eq:computad-category-pair} exactly when
		$H(\partial^-\alpha)=H(\partial^+\alpha)$ for every $\alpha$.  This is the
		universal property of the quotient by the least categorical congruence
		generated by those equations.
	\end{proof}
	
	Let $T_1$ be the free-category monad on $\Grph$.  The maps in
	\eqref{eq:computad-kleisli-pair} are a parallel pair of Kleisli arrows for
	$T_1$.  The $T_1$-algebra axioms give precisely identities and associative
	composition on the underlying graph, so the Eilenberg--Moore category
	$\Grph^{T_1}$ is canonically isomorphic to $\Cat$.
	
	\begin{theorem}[Computadic relations as $T_1$-relation data]
		\label{thm:computadic-monadic-comparison}
		The $T_1$-algebra presented by the Kleisli pair
		\eqref{eq:computad-kleisli-pair} is $\Pp_1C$.  On objects, this identifies
		ordinary $2$-computads with the $T_1$-relation data whose
		relation graph is equipped with a chosen identification with a copower of
		$\2_{\Grph}$ and whose two boundary maps agree on objects.  With
		generator-preserving maps this identification is functorial.
		Every small category admits a presentation of this form.
	\end{theorem}
	
	\begin{proof}
		Under the free-category adjunction, the Kleisli pair induces exactly the two
		functors in \eqref{eq:computad-category-pair}, so the first assertion is
		\cref{prop:computad-category-coequalizer}.
		
		A generator-preserving computad map determines a morphism of these Kleisli
		pairs: its functions on objects and arrows give the generator-graph map, and
		its function on $2$-generators gives the map between the copowers of
		$\2_{\Grph}$.  Boundary preservation is exactly the pair of commutative
		squares required of a morphism of relation data.  The induced map between the
		coequalizers is the functor $\Pp_1f$, and this assignment preserves identities
		and composition.  This proves the functoriality assertion.
		
		For a small category $A$, take $C_0=\Ob A$ and take one $1$-generator
		$[f]\colon x\to y$ for every arrow $f\colon x\to y$ of $A$.  Add a relation
		$[1_x]=1_x$ for every object and a relation
		$[gf]=[g][f]$ for every composable pair.  Evaluation defines a functor
		$\Pp_1C\to A$.  The relations reduce every nonempty path to the single
		generator labeled by its composite and every identity-labeled generator to
		the empty path.  Evaluation is therefore bijective on objects and hom-sets.
	\end{proof}
	
	The last assertion is an existence statement, not a minimality claim.
	Choosing a smaller generating graph and a smaller relation set is the
	substance of presentation theory.
	
	The preceding theorem singles out the relation data which already have
	computadic form.  We call a $T_1$-relation datum \emph{object-agreeing}, or a
	\emph{$1$-cell relation datum}, when its two object functions are equal.  An
	arbitrary datum may instead identify objects.  The object equations can be
	imposed first, after which the remaining equations are again boundaries of
	generating $2$-cells.
	
	\begin{theorem}[Object normalization of $T_1$-presentations]
		\label{thm:primeiro-os-objetos}
		Let $R$ and $G$ be small graphs and let
		\begin{equation}
			\label{eq:datum-antes-dos-objetos}
			r_0,r_1\colon R\rightrightarrows
			T_1G=U_{\Grph}\freecat G
		\end{equation}
		be a $T_1$-relation datum.  Write
		\[
		r_{i,0}\colon R_0\longrightarrow G_0,
		\qquad
		r_{i,1}\colon R_1\longrightarrow\Mor(\freecat G)
		\qquad(i=0,1)
		\]
		for its object and arrow functions.  Form the coequalizer
		\begin{equation}
			\label{eq:primeiro-coequalizamos-os-objetos}
			\begin{tikzcd}[column sep=large]
				R_0
				\ar[r,shift left=.65ex,"r_{0,0}"]
				\ar[r,shift right=.65ex,swap,"r_{1,0}"]&
				G_0 \ar[r,two heads,"q_0"]&
				\overline G_0
			\end{tikzcd}
		\end{equation}
		in $\Set$.  Define a graph $\overline G$ by
		\[
		\overline G_1=G_1,\qquad
		s_{\overline G}=q_0s_G,\qquad
		t_{\overline G}=q_0t_G,
		\]
		and let $q=(q_0,1_{G_1})\colon G\to\overline G$ be the resulting graph map.
		For every $\rho\in R_1$, put
		\begin{equation}
			\label{eq:as-relacoes-depois-dos-objetos}
			\overline\partial^{-}\rho
			=(\freecat q)\bigl(r_{0,1}(\rho)\bigr),
			\qquad
			\overline\partial^{+}\rho
			=(\freecat q)\bigl(r_{1,1}(\rho)\bigr).
		\end{equation}
		These are parallel paths in $\freecat\overline G$.  Hence they define a
		$2$-computad $\overline C$ with
		\[
		\overline C_{\leq1}=\overline G,\qquad
		\overline C_2=R_1.
		\]
		The composite
		\begin{equation}
			\label{eq:quociente-normalizado}
			\freecat G\xrightarrow{\freecat q}
			\freecat\overline G\longrightarrow\Pp_1\overline C
		\end{equation}
		is a coequalizer in $\Cat$ of the Kleisli extensions
		\[
		r_0^\sharp,r_1^\sharp\colon\freecat R\rightrightarrows\freecat G.
		\]
		Consequently, the original datum and the object-agreeing $1$-cell
		presentation associated with $\overline C$ present canonically isomorphic
		categories.
		
		With the quotient in \eqref{eq:primeiro-coequalizamos-os-objetos} taken to be
		the set of equivalence classes, this normalization is functorial: a morphism
		of $T_1$-relation data induces a generator-preserving map between the
		normalized $2$-computads.
	\end{theorem}
	
	\begin{proof}
		Since $r_i$ is a graph map, for $\rho\in R_1$ we have
		\[
		\begin{aligned}
			\operatorname{source}\bigl(r_{i,1}(\rho)\bigr)&=r_{i,0}(s_R\rho),&
			\operatorname{target}\bigl(r_{i,1}(\rho)\bigr)&=r_{i,0}(t_R\rho).
		\end{aligned}
		\]
		The coequalizer equation $q_0r_{0,0}=q_0r_{1,0}$ therefore identifies the
		sources of the two paths in
		\eqref{eq:as-relacoes-depois-dos-objetos}, and it also identifies their
		targets.  Thus $\overline\partial^{-}\rho$ and
		$\overline\partial^{+}\rho$ are parallel, so $\overline C$ is well defined.
		
		Let
		\[
		p\colon\freecat\overline G\longrightarrow\Pp_1\overline C
		\]
		be the quotient functor.  The composite $p\,\freecat q$ coequalizes
		$r_0^\sharp$ and $r_1^\sharp$.  Indeed, on an object $v\in R_0$ this is the
		equation
		\[
		q_0r_{0,0}(v)=q_0r_{1,0}(v),
		\]
		and on a generating arrow $\rho\in R_1$ it is the defining relation
		\[
		p\bigl(\overline\partial^{-}\rho\bigr)
		=p\bigl(\overline\partial^{+}\rho\bigr).
		\]
		Since $\freecat R$ is free on $R$, equality on its objects and generating
		arrows gives equality of the two functors.
		
		It remains to prove universality.  Let $A$ be a category and let
		$H\colon\freecat G\to A$ be a functor.  The equality
		\[
		Hr_0^\sharp=Hr_1^\sharp
		\]
		holds if and only if the two functors agree on the generators of
		$\freecat R$.  Explicitly, this means
		\begin{align}
			H_0r_{0,0}(v)&=H_0r_{1,0}(v)
			&&(v\in R_0),                         \label{eq:normalizacao-nos-objetos}\\
			H\bigl(r_{0,1}(\rho)\bigr)&=H\bigl(r_{1,1}(\rho)\bigr)
			&&(\rho\in R_1).                      \label{eq:normalizacao-nas-setas}
		\end{align}
		By \eqref{eq:primeiro-coequalizamos-os-objetos},
		\eqref{eq:normalizacao-nos-objetos} is equivalent to a unique function
		\[
		\overline H_0\colon\overline G_0\longrightarrow\Ob A
		\quad\text{with}\quad
		\overline H_0q_0=H_0.
		\]
		For an arrow $e\in\overline G_1=G_1$, set
		$\overline H(e)=H(e)$.  Its source and target are correct because
		$\overline H_0q_0=H_0$.  The universal property of the free category now
		gives a unique functor
		\[
		\overline H\colon\freecat\overline G\longrightarrow A
		\quad\text{such that}\quad
		\overline H\,\freecat q=H.
		\]
		Under this factorization, \eqref{eq:normalizacao-nas-setas} says precisely
		\[
		\overline H\bigl(\overline\partial^{-}\rho\bigr)
		=\overline H\bigl(\overline\partial^{+}\rho\bigr)
		\qquad(\rho\in R_1).
		\]
		By \cref{prop:computad-category-coequalizer}, this is equivalent to a unique
		factorization of $\overline H$ through $p$.  We have therefore obtained,
		naturally in $A$, a bijection between functors
		$\Pp_1\overline C\to A$ and functors
		$H\colon\freecat G\to A$ which coequalize
		$r_0^\sharp,r_1^\sharp$.  This proves the coequalizer assertion.
		
		Finally, let
		$(u,v)\colon(R,G,r_0,r_1)\to(R',G',r'_0,r'_1)$ be a morphism of
		$T_1$-relation data, and write
		$q'\colon G'\to\overline G'$ for the quotient map in the normalization of
		the target datum.  Thus
		\[
		T_1v\,r_i=r'_iu\qquad(i=0,1).
		\]
		Its object equation sends every generating equivalence
		$r_{0,0}(a)\sim r_{1,0}(a)$ to
		$r'_{0,0}(u_0a)\sim r'_{1,0}(u_0a)$.  It induces a function
		\[
		\overline v_0\colon\overline G_0\longrightarrow\overline G'_0,
		\qquad
		\overline v_0([x])=[v_0(x)].
		\]
		Together with $v_1\colon G_1\to G'_1$, this is a graph map
		$\overline v\colon\overline G\to\overline G'$.  On $2$-generators use
		$u_1\colon R_1\to R'_1$.  From
		$\overline vq=q'v$ and $T_1v\,r_i=r'_iu$ we obtain
		\[
		(\freecat\overline v)(\overline\partial^\pm\rho)
		=\overline\partial'{}^\pm(u_1\rho).
		\]
		Thus $(\overline v,u_1)$ is a generator-preserving map of $2$-computads.
		The construction respects identities and composition, which proves the final
		assertion.
	\end{proof}
	
	\begin{remark}
		\label{rem:novos-caminhos-depois-do-quociente}
		The functor $\freecat q$ need not be full.  Identifying two objects can make
		arrows composable which were not composable in $G$, and
		$\freecat\overline G$ then contains the resulting new paths.  This is why the
		normalization uses the free category on the quotient graph, rather than an
		arrowwise quotient of $\freecat G$.  The two presentations are equivalent in
		the precise sense that their coequalizers are canonically isomorphic; their
		relation data need not be isomorphic.
	\end{remark}
	
	\subsection{Groupoidal computads and classical group presentations}
	
	\begin{definition}[Groupoidal $2$-computad]
		\label{def:groupoidal-computad}
		A \emph{groupoidal $2$-computad} has a graph $G$, a set $R$, and parallel
		boundary words
		\[
		\partial^-,\partial^+\colon R\longrightarrow\Mor(\freegpd G).
		\]
		It presents the quotient groupoid
		\[
		\Pp_{\mathrm{gpd}}C
		=\freegpd G/
		(\partial^-\rho=\partial^+\rho\mid\rho\in R).
		\]
		An ordinary $2$-computad determines a groupoidal one by applying
		$\freecat G\to\freegpd G$ to its boundaries, and the resulting groupoid is
		canonically $\gpcomp\Pp_1C$.
	\end{definition}
	
	Since $\gpcomp$ preserves coproducts and coequalizers, applying it to
	\eqref{eq:computad-category-pair} replaces the walking arrow by the walking
	isomorphism and yields exactly the displayed groupoidal quotient.
	
	Presentations of groupoids and their use in the study of group
	presentations go back to Higgins \cite{Higgins1964}.  The formulation above
	keeps the computadic boundary functions explicit.
	
	A map $f\colon C\to D$ of groupoidal $2$-computads consists of a function
	on objects, a signed path in $\freegpd D_{\leq1}$ with the transported
	endpoints for every generating arrow of $C$, and, for every
	$\alpha\in C_2$, a $2$-cell in the free strict $2$-groupoid on $D$ from
	$f(\partial^-\alpha)$ to $f(\partial^+\alpha)$.  This free strict
	$2$-groupoid is characterized by the following universal property:
	assigning its $1$- and $2$-cells into a strict $2$-groupoid is equivalent
	to assigning the graph arrows and generating $2$-cells of $D$ compatibly
	with their signed boundaries.  We write
	$2\text{-}\Comp_{\mathrm{gpd}}$ for the resulting category and
	$2\text{-}\Comp_{\mathrm{gpd}}^{\mathrm{gen}}$ for its subcategory of
	generator-preserving maps.
	
	\begin{definition}[Connected computad]
		\label{def:connected-computad}
		An ordinary or groupoidal $2$-computad is \emph{connected} when the underlying
		unoriented multigraph of its generating graph is connected.  Equivalently,
		its free category or free groupoid is connected; equivalently again, the
		category or groupoid which it presents is connected.
	\end{definition}
	
	Indeed, the free object has an arrow between two objects precisely when
	there is a zigzag between them in the generating graph.  Imposing relations
	identifies parallel arrows but neither creates arrows between graph
	components nor identifies objects.
	
	The following right adjoint makes the relation object explicit: it has one
	generator for each ordered pair of parallel arrows.
	
	\begin{proposition}[Cofree groupoidal computad]
		\label{prop:cofree-groupoidal-computad}
		For a graph $G$, let $R(G)$ have underlying graph $G$ and one $2$-generator
		$[u,v]\colon u\Rightarrow v$ for every ordered pair of parallel arrows
		$u,v$ in $\freegpd G$.  Then $R$ is right adjoint to the functor
		\[
		U_{\mathrm{gen}}\colon
		2\text{-}\Comp_{\mathrm{gpd}}^{\mathrm{gen}}\longrightarrow\Grph
		\]
		that forgets the $2$-generators of a groupoidal computad with
		generator-preserving maps.
	\end{proposition}
	
	\begin{proof}
		Given a graph map $f\colon C_{\leq1}\to G$, the image of a generating
		$2$-cell $\rho\colon u\Rightarrow v$ of $C$ is forced to be the generator
		$[\freegpd(f)(u),\freegpd(f)(v)]$ of $R(G)$.  This gives a unique computad map
		over $f$, naturally in $C$ and $G$.
	\end{proof}
	
	Write $\Sigma H$ for the one-object groupoid associated with a group $H$.
	
	\begin{theorem}[Comparison with group presentations]
		\label{thm:group-presentation-comparison}
		Let $C$ be a connected groupoidal $2$-computad and $T$ a maximal tree in its
		generating graph.  Collapsing $T$ produces a one-object group presentation
		with generating set $C_1\setminus T_1$ and one relator for every element of
		$C_2$.  It presents the isotropy group of $\Pp_{\mathrm{gpd}}C$ at the chosen
		root.  If the generating graph is finite or rank-finite, then
		\[
		\abs{C_1\setminus T_1}=\beta_1(C_{\leq1}).
		\]
		Conversely, every group presentation is the one-object case of a groupoidal
		computad.
	\end{theorem}
	
	\begin{proof}
		Use the tree words $\tau_x$ of \cref{thm:maximal-forest-normal-form}.  Replace
		each occurrence of an arrow $e\colon x\to y$ in a boundary word by the rooted
		loop $\tau_y^{-1}e\tau_x$, deleting the tree letters after collapse.  Each
		equation $u=v$ becomes the relator $v^{-1}u=1$ in the free group on the
		chords.  The maximal-forest normal form identifies the resulting quotient
		group with the chosen isotropy group.  A classical presentation
		$\langle S\mid R\rangle$ gives the converse computad with one object, loop
		arrows $S$, and one $2$-generator $r\Rightarrow1$ for each relator.
	\end{proof}
	
	This comparison will identify groupoid deficiency with classical group
	deficiency in
	\cref{thm:connected-groupoid-isotropy-deficiency}.
	
	\subsection{The simplicial category}
	\label{sec:simplicial-presentation}
	
	The simplicial category provides a familiar example in which the typing
	conventions matter.  Let
	$[n]=\{0<\cdots<n\}$.  For $n\geq1$ and $0\leq i\leq n$, let
	$\delta_i^n\colon[n-1]\to[n]$ be the injection omitting $i$.  For $n\geq0$
	and $0\leq i\leq n$, let
	$\sigma_i^n\colon[n+1]\to[n]$ be the surjection identifying $i$ and $i+1$.
	
	\begin{theorem}[Simplicial presentation]
		\label{thm:simplicial-presentation}
		The category $\Delta$ of finite nonempty ordinals and monotone maps is
		presented by the objects $[n]$, the generators $\delta_i^n,\sigma_i^n$, and
		the following typed relations.  In the mixed relation,
		$0\leq j\leq n$ and $0\leq i\leq n+1$.
		\begin{align}
			\delta_j^{n+1}\delta_i^n
			&=\delta_i^{n+1}\delta_{j-1}^n
			&&(n\geq1,\ 0\leq i<j\leq n+1),                              \label{eq:coface}\\
			\sigma_j^n\sigma_i^{n+1}
			&=\sigma_i^n\sigma_{j+1}^{n+1}
			&&(n\geq0,\ 0\leq i\leq j\leq n),                            \label{eq:codegeneracy}\\
			\sigma_j^n\delta_i^{n+1}
			&=
			\begin{cases}
				\delta_i^n\sigma_{j-1}^{n-1},&i<j,\\
				1_{[n]},&i=j\text{ or }i=j+1,\\
				\delta_{i-1}^n\sigma_j^{n-1},&i>j+1.
			\end{cases}
			&&(n\geq0).                                                   \label{eq:mixed}
		\end{align}
		Terms with superscript $n-1$ occur only when their displayed indices make
		$n\geq1$ automatic.
	\end{theorem}
	
	\begin{proof}
		The displayed identities follow directly from the definitions of the
		injections and surjections.  Orient the mixed identities from a subword
		$\sigma\delta$ to a subword $\delta\sigma$, or to the identity in the two
		cancellation cases.  Let $\operatorname{inv}(w)$ be the number of pairs
		consisting of an occurrence of a codegeneracy to the left of an occurrence
		of a coface in $w$.  A noncancelling mixed rewrite decreases
		$\operatorname{inv}(w)$, while a cancelling rewrite decreases the word
		length.  Hence the lexicographically ordered pair
		\[
		\bigl(\operatorname{length}(w),\operatorname{inv}(w)\bigr)
		\]
		strictly decreases at every mixed rewrite.  The mixed rewriting terminates
		and leaves a coface block followed by a codegeneracy block.
		
		Within the coface block, the coface identities put the omitted values into
		their canonical order; within the codegeneracy block, the codegeneracy
		identities put the nontrivial fibres into their canonical order.  These are
		the standard normal forms for monotone injections and surjections,
		respectively
		\cite[Chapter~I, \S2, equation~(3)]{May1967}.  Thus every word is equal,
		using the displayed relations, to an epi--mono word.  The monotone map
		represented by that word determines its epi--mono factorization uniquely:
		the omitted values determine the coface block, and the nontrivial fibres
		determine the codegeneracy block.  Hence two normal words representing the
		same monotone map coincide, and the induced functor from the presented
		category to $\Delta$ is bijective on every hom-set.
	\end{proof}
	
	\section{Computads and iterated icons}
	\label{sec:higher-icons}
	
	The higher-dimensional theorem rests on two choices: the unital free
	enriched-category monad and the class of $2$-cells used between strict higher
	functors.
	
	\subsection{The free enriched-category monad}
	
	Let $\mathcal V$ be a cocomplete cartesian monoidal category in which products
	preserve coproducts in each variable.  A $\mathcal V$-graph $A$ consists of a
	set $A_0$ and objects $A(x,y)$ of $\mathcal V$.  The free
	$\mathcal V$-category monad has hom-objects
	\begin{equation}
		\label{eq:free-enriched-category-monad}
		(T_{\mathcal V}A)(x,y)=
		\coprod_{k\geq0}\ 
		\coprod_{x=x_0,\ldots,x_k=y}
		A(x_{k-1},x_k)\times\cdots\times A(x_0,x_1).
	\end{equation}
	For $k=0$, the product is the terminal monoidal unit when $x=y$; there is no
	such summand when $x\neq y$.
	
	On a morphism of $\mathcal V$-graphs the endofunctor $T_{\mathcal V}$ applies
	the given hom-morphisms to every factor and the given object function to every
	indexing string.  Its unit includes $A(x,y)$ as the length-one summand.  Its
	multiplication flattens a string of strings; the hypothesis that products
	preserve coproducts supplies the canonical distributivity isomorphisms needed
	to identify the resulting iterated coproduct with
	\eqref{eq:free-enriched-category-monad}.  On every coproduct summand the two
	possible iterated multiplications are the same flattening map, and the unit
	maps insert empty or one-letter strings.  The monad laws therefore follow from
	the universal property of the coproducts, so these maps define a monad on the
	category of small $\mathcal V$-graphs.
	
	\begin{proposition}
		\label{prop:free-enriched-category}
		The algebras of $T_{\mathcal V}$ are small $\mathcal V$-categories.  The
		$k=0$ summand in \eqref{eq:free-enriched-category-monad} supplies identity
		arrows, while multiplication concatenates strings.
	\end{proposition}
	
	\begin{proof}
		An algebra map evaluates every composable string.  Its restriction to the
		$k=0$ and $k=2$ summands gives identities and composition.  The algebra unit
		and multiplication axioms give the unit and associativity axioms.  Conversely,
		iterated composition in a $\mathcal V$-category evaluates every summand, and
		the category axioms give a $T_{\mathcal V}$-algebra.  The two constructions
		are inverse.  An algebra homomorphism preserves the length-zero and
		length-two evaluations exactly when its hom-morphisms preserve identities
		and composition; hence the algebra homomorphisms are precisely the enriched
		functors.  This is the usual free enriched-category construction
		\cite{Kelly1982}.
	\end{proof}
	
	For every $m\geq0$, the category $m\text{-}\Cat$ is cocomplete and cartesian
	closed.  Inductively, cocompleteness follows from the local-presentability
	theorem for enriched categories \cite{KellyLack2001}, while the cartesian
	closed structure is the corresponding enriched functor-category
	construction \cite[Sections~2.2--2.3]{Kelly1982}.  Hence cartesian product
	preserves coproducts in each variable.  Taking
	$\mathcal V=(n-1)\text{-}\Cat$ constructs the free strict $n$-category on an
	$(n-1)\text{-}\Cat$-enriched graph.  Applied after freely constructing the
	lower computadic skeleton, this is the recursive free strict $n$-category
	used below.  The $k=0$ summand supplies units.
	
	\subsection{\texorpdfstring{$n$-icons}{n-icons}}
	
	For the strict categories used here, the relevant part of the iterated-icon
	construction admits a concise description.
	
	In dimension two, an icon is defined between strict $2$-functors that agree
	on objects.  It assigns a $2$-cell to every $1$-cell, compatibly with
	identities, composition, and naturality with respect to $2$-cells
	\cite{Lack2010Icons}.  For $n\geq2$, the $n$-icons used here are the
	corresponding top-dimensional components of the iterated-icon construction:
	the functors agree in dimensions at most $n-2$, and only the components in
	dimension $n$ may be nonidentity.
	
	More precisely, when $n=2$, a $1$-cell $u\colon x\to y$ has the component
	\[
	\begin{tikzcd}[column sep=huge]
		F(x)=G(x)
		\ar[r,bend left=38,"F(u)",""{name=upper,below}]
		\ar[r,bend right=38,"G(u)"',""{name=lower,above}]&
		F(y)=G(y).
		\arrow[Rightarrow,from=upper,to=lower,"\alpha_u" description]
	\end{tikzcd}
	\]
	The equations below say that these components compose as the $1$-cells do
	and are natural with respect to the $2$-cells.
	
	For $n=1$ we put
	$1\text{-}\Cat_{\mathrm{icon}}=\Cat$, with ordinary natural transformations as
	$2$-cells.  For higher cells, $v\circ_j u$ denotes composition along the
	common $j$-dimensional boundary, written in the paper's right-to-left order.
	The definition is as follows.
	
	\begin{definition}[$n$-icon]
		\label{def:n-icon}\label{def:os-icones-nao-mexem-nos-objetos}
		Let $n\geq2$, and let $F,G\colon A\to B$ be strict $n$-functors that agree on
		all cells of dimension at most $n-2$.  An \emph{$n$-icon}
		$\alpha\colon F\Rightarrow G$ assigns to every $(n-1)$-cell $u$ of $A$ an
		$n$-cell
		\[
		\alpha_u\colon F(u)\Rightarrow G(u).
		\]
		If $w$ is an $(n-2)$-cell, if $u,v$ are $j$-composable
		$(n-1)$-cells, and if $\Gamma\colon u\Rightarrow v$ is an $n$-cell, the
		following equations are required:
		\begin{align}
			\alpha_{1_w}&=1_{1_{F(w)}},                                      \label{eq:n-icon-unit}\\
			\alpha_{v\circ_j u}&=\alpha_v\circ_j\alpha_u
			&&(0\leq j<n-1),                                          \label{eq:n-icon-composition}\\
			G(\Gamma)\circ_{n-1}\alpha_u
			&=\alpha_v\circ_{n-1}F(\Gamma).                               \label{eq:n-icon-naturality}
		\end{align}
		For $n=2$ these are exactly the axioms for an icon between strict
		$2$-functors.
	\end{definition}
	
	In dimension two the naturality equation says that the following square in
	the hom-category of $B$ commutes:
	\begin{equation}
		\label{eq:o-quadrado-do-icon}
		\begin{tikzcd}[column sep=huge,row sep=large]
			F(u) \ar[r,Rightarrow,"\alpha_u"]
			\ar[d,Rightarrow,swap,"F(\Gamma)"]&
			G(u) \ar[d,Rightarrow,"G(\Gamma)"]\\
			F(v) \ar[r,Rightarrow,swap,"\alpha_v"]&
			G(v).
		\end{tikzcd}
	\end{equation}
	Thus an icon has no component at an object.  Its first nontrivial components
	are the $2$-cells attached to $1$-cells, which is exactly why icons detect the
	operation of adjoining generating $2$-cells without changing the objects.
	
	\begin{proposition}
		\label{prop:n-icon-two-category}
		Strict $n$-categories, strict $n$-functors, and $n$-icons form a strict
		$2$-category, denoted by $n\text{-}\Cat_{\mathrm{icon}}$.
	\end{proposition}
	
	\begin{proof}
		Vertical composition and identities are pointwise in the top direction:
		\[
		(\beta\alpha)_u=\beta_u\circ_{n-1}\alpha_u,
		\qquad (1_F)_u=1_{F(u)}.
		\]
		Whiskering is given by
		$(H\alpha)_u=H(\alpha_u)$ and $(\alpha L)_v=\alpha_{L(v)}$.
		For $n$-icons
		\[
		\alpha\colon F\Rightarrow G\colon A\to B,
		\qquad
		\beta\colon H\Rightarrow K\colon B\to C,
		\]
		their horizontal composite has component
		\begin{equation}
			\label{eq:n-icon-horizontal-composition}
			(\beta\ast\alpha)_u
			=\beta_{G(u)}\circ_{n-1}H(\alpha_u)
			=K(\alpha_u)\circ_{n-1}\beta_{F(u)}.
		\end{equation}
		The two expressions are equal by
		\eqref{eq:n-icon-naturality} for $\beta$, applied to the $n$-cell
		$\alpha_u$.  Equations
		\eqref{eq:n-icon-unit}--\eqref{eq:n-icon-naturality} are preserved by
		these operations.  Strict associativity and the unit laws follow from strict
		functoriality, and middle-four interchange follows from interchange in $C$.
		Thus the displayed data form a strict $2$-category.
	\end{proof}
	
	This is the strict top-dimensional fragment of the iterated enriched-icon
	construction of \cite{ChengGurski2014}.  It is exactly the fragment required
	by free attachment: the lower-dimensional data remain fixed and the first
	nonidentity components occur in dimension $n$.  The agreement condition is
	important: an arbitrary $n$-natural transformation has lower-dimensional
	components and does not encode the operation of adjoining only top cells.
	
	\begin{lemma}[The walking-cell calculation]
		\label{lem:walking-cell-icons}
		Let $n\geq2$, let $A=S\copow\walking{n-1}$, and let
		$f,g\colon A\rightrightarrows B$ be strict $(n-1)$-functors that agree below
		dimension $n-1$.  For a strict $n$-functor
		$h\colon d_nB\to X$, $n$-icons
		\[
		h d_n f\Rightarrow h d_n g
		\]
		are in bijection with families of $n$-cells
		\[
		\xi_s\colon h(f\kappa_{n-1,s})\Rightarrow
		h(g\kappa_{n-1,s})\qquad(s\in S),
		\]
		where $\kappa_{n-1,s}$ is the distinguished cell in the $s$th copower
		summand.
	\end{lemma}
	
	\begin{proof}
		The unit axiom fixes the components on identity $(n-1)$-cells, and the only
		nonidentity top-dimensional generator of each globe is
		$\kappa_{n-1,s}$.  Compatibility with composition therefore determines the
		icon on that summand from its component at $\kappa_{n-1,s}$.  Different
		copower summands have no nonidentity cells between them, so the choices are
		independent.  Naturality at the identity $n$-cells of $d_nA$ is automatic.
	\end{proof}
	
	\subsection{Derivation schemes and free extensions}
	
	Computads are obtained by adjoining top-dimensional generators to a
	\emph{free} lower skeleton.  The coinserter calculation itself does not
	require that freeness.  We isolate the more general datum first.
	
	\begin{definition}[Derivation $n$-scheme]
		\label{def:derivation-n-scheme}
		Let $n\geq1$.  A \emph{derivation $n$-scheme} is a quadruple
		\[
		D=(B,S,b^-,b^+)
		\]
		consisting of a strict $(n-1)$-category $B$, a set $S$, and strict
		$(n-1)$-functors
		\[
		b^-,b^+\colon S\copow\walking{n-1}\rightrightarrows B
		\]
		which agree on every cell of dimension below $n-1$.  Thus, for each
		$s\in S$, the distinguished cells
		$b^-(\kappa_{n-1,s})$ and $b^+(\kappa_{n-1,s})$ are parallel
		$(n-1)$-cells of $B$.
		
		The \emph{free extension by $n$-cells} $\operatorname{Ext}_n(D)$ has $B$ as its
		$(n-1)$-truncation and has one generating $n$-cell
		\[
		\xi_s\colon b^-(\kappa_{n-1,s})
		\Rightarrow b^+(\kappa_{n-1,s})
		\qquad(s\in S).
		\]
		Its $n$-cells are the well-typed expressions obtained from the $\xi_s$ and
		identity $n$-cells by all strict compositions, modulo only the structural
		equations of a strict globular $n$-category.  These include associativity
		and units for every composition, interchange, and the compatibility of
		identity $n$-cells with lower-dimensional composition; for instance, for
		$j$-composable $(n-1)$-cells $u,v$,
		\[
		1_{v\circ_j u}=1_v\circ_j1_u
		\qquad(0\leq j<n-1).
		\]
	\end{definition}
	
	\begin{definition}[Maps of derivation schemes]
		\label{def:maps-of-derivation-schemes}
		Let
		\[
		D=(B,S,b^-,b^+),\qquad
		D'=(B',S',b'^-,b'^+)
		\]
		be derivation $n$-schemes.  A \emph{generator-preserving morphism}
		\[
		(F,\varphi)\colon D\longrightarrow D'
		\]
		consists of a strict $(n-1)$-functor $F\colon B\to B'$ and a function
		$\varphi\colon S\to S'$ such that
		\begin{equation}
			\label{eq:morphism-of-derivation-schemes}
			F b^\epsilon
			=b'^\epsilon(\varphi\copow1_{\walking{n-1}})
			\qquad(\epsilon\in\{-,+\}).
		\end{equation}
		Composition is componentwise.  We denote the resulting category by
		$n\text{-}\mathsf{Der}$.
	\end{definition}
	
	The assignment $D\mapsto\operatorname{Ext}_n(D)$ is functorial on these
	morphisms.  Indeed, $(F,\varphi)$ extends uniquely to the strict $n$-functor
	\[
	\operatorname{Ext}_n(F,\varphi)\colon
	\operatorname{Ext}_n(D)\longrightarrow\operatorname{Ext}_n(D')
	\]
	whose lower truncation is $F$ and which sends
	$\xi_s$ to $\xi'_{\varphi(s)}$.  We thus obtain a functor
	\[
	\operatorname{Ext}_n\colon n\text{-}\mathsf{Der}
	\longrightarrow n\text{-}\Cat .
	\]
	
	When $n=1$, the category $B$ is a set, the two maps give the endpoints of
	the arrows indexed by $S$, and $\operatorname{Ext}_1(D)$ is the free category
	on that graph.  For $n\geq2$, this construction freely adjoins
	top-dimensional cells to an already constructed strict lower category.
	
	\begin{theorem}[Derivation-scheme coinserter]
		\label{thm:derivation-scheme-coinserter}
		Let $D=(B,S,b^-,b^+)$ be a derivation $n$-scheme.  In
		$n\text{-}\Cat_{\mathrm{icon}}$, the diagram
		\[
		\begin{tikzcd}[column sep=large]
			d_n(S\copow\walking{n-1})
			\ar[r,shift left=.7ex,"d_nb^-"]
			\ar[r,shift right=.7ex,swap,"d_nb^+"]&
			d_nB \ar[r,"j"]&
			\operatorname{Ext}_n(D)
		\end{tikzcd}
		\]
		is a coinserter.  The map $j$ is the identity through dimension $n-1$, and
		the component of the universal $n$-icon at the distinguished cell in the
		$s$th summand is $\xi_s$.
	\end{theorem}
	
	\begin{proof}
		Let $X$ be a strict $n$-category.  A cone with vertex $X$ consists of a
		strict $n$-functor $h\colon d_nB\to X$ and an $n$-icon
		\[
		\zeta\colon hd_nb^-\Rightarrow hd_nb^+.
		\]
		For $n\geq2$, \cref{lem:walking-cell-icons} identifies $\zeta$ with a family
		of $n$-cells
		\[
		h\bigl(b^-(\kappa_{n-1,s})\bigr)
		\xRightarrow{\ \zeta_s\ }
		h\bigl(b^+(\kappa_{n-1,s})\bigr)
		\qquad(s\in S).
		\]
		For $n=1$, the same statement says that a natural transformation between
		two functors out of the discrete category on $S$ is a family of arrows with
		the prescribed endpoints.
		
		By the term construction in
		\cref{def:derivation-n-scheme}, this family extends uniquely to a strict
		$n$-functor
		\[
		\bar h\colon\operatorname{Ext}_n(D)\longrightarrow X
		\]
		which restricts to $h$ and sends $\xi_s$ to $\zeta_s$.  If
		$(h,\zeta)$ and $(k,\omega)$ are two cones, a morphism between them is an
		$n$-icon $\gamma\colon h\Rightarrow k$ satisfying
		\[
		(\gamma d_nb^+)\zeta
		=\omega(\gamma d_nb^-).
		\]
		The equation says that the following square of $n$-icons
		commutes:
		\[
		\begin{tikzcd}[column sep=huge,row sep=huge]
			hd_nb^- \ar[r,Rightarrow,"\zeta"]
			\ar[d,Rightarrow,swap,"\gamma d_nb^-"]&
			hd_nb^+ \ar[d,Rightarrow,"\gamma d_nb^+"]\\
			kd_nb^- \ar[r,Rightarrow,swap,"\omega"]&
			kd_nb^+ .
		\end{tikzcd}
		\]
		Its component at the $s$th distinguished cell is precisely the naturality
		equation required at the generator $\xi_s$.  Induction over formal
		composites, using the icon composition equations and interchange, extends
		$\gamma$ uniquely to an $n$-icon $\bar h\Rightarrow\bar k$.  Restriction and
		extension are inverse on objects and morphisms and commute with
		postcomposition in $X$.  Hence they give a natural isomorphism between the
		mapping category out of $\operatorname{Ext}_n(D)$ and the category of cones,
		which is the coinserter universal property.
	\end{proof}
	
	\begin{example}[The triangle revisited]
		\label{ex:o-triangulo-volta}
		Let $G$ be the graph of
		\eqref{eq:o-triangulo-que-explica-tudo}, put $B=\freecat G$, and let
		$S=\{\alpha\}$.  The equations
		\[
		b^-(\kappa_{1,\alpha})=ba,
		\qquad
		b^+(\kappa_{1,\alpha})=c
		\]
		define a derivation $2$-scheme.  Its free extension is the strict
		$2$-category obtained from $\freecat G$ by adjoining the single $2$-cell
		\[
		\alpha\colon ba\Rightarrow c.
		\]
		The coinserter remembers $\alpha$ as a cell.  If the same boundary pair is
		instead used as relation data for the underlying category, its coequifier
		imposes $ba=c$ and gives the ordinal category $[2]$, as in
		\cref{ex:o-triangulo-como-apresentacao}.  Thus the two operations in the
		opening example are precisely the coinserter and coequifier constructions.
	\end{example}
	
	\begin{definition}[The derivation scheme underlying a strict category]
		\label{def:underlying-derivation-scheme}
		Let $X$ be a strict $n$-category.  Write $X_n$ for its set of $n$-cells and
		$X_{\leq n-1}$ for its lower truncation.  The source and target of every
		$\xi\in X_n$ are parallel $(n-1)$-cells, and therefore determine strict
		$(n-1)$-functors
		\[
		c_X^-,c_X^+\colon
		X_n\copow\walking{n-1}\rightrightarrows X_{\leq n-1}
		\]
		which agree below dimension $n-1$.  Define
		\[
		\mathcal C_{n\text{-}\mathrm{Der}}(X)
		=\bigl(X_{\leq n-1},X_n,c_X^-,c_X^+\bigr).
		\]
		For a strict $n$-functor $H\colon X\to Y$, its lower truncation and its
		function $H_n\colon X_n\to Y_n$ form a morphism of derivation schemes.
		This defines a functor
		\[
		\mathcal C_{n\text{-}\mathrm{Der}}\colon
		n\text{-}\Cat\longrightarrow n\text{-}\mathsf{Der}.
		\]
		For $n=1$, it is the ordinary underlying-graph functor.  Notice that this
		functor includes the lower category itself; it is therefore different from
		the recursively underlying \emph{computad} functor constructed below.
	\end{definition}
	
	\begin{theorem}[The derivation-scheme adjunction]
		\label{thm:as-derivacoes-sabem-o-caminho}
		For every $n\geq1$, there is an adjunction
		\[
		\operatorname{Ext}_n\colon n\text{-}\mathsf{Der}
		\rightleftarrows n\text{-}\Cat
		\mathbin{:}\mathcal C_{n\text{-}\mathrm{Der}} .
		\]
		Its counit
		\[
		\operatorname{Ext}_n\mathcal C_{n\text{-}\mathrm{Der}}(X)
		\longrightarrow X
		\]
		is the identity through dimension $n-1$ and sends the generating copy of an
		$n$-cell $\xi$ to $\xi$ itself.
	\end{theorem}
	
	\begin{proof}
		Let $D=(B,S,b^-,b^+)$ be a derivation $n$-scheme and let $X$ be a strict
		$n$-category.  A morphism
		\[
		(h,\varphi)\colon
		D\longrightarrow\mathcal C_{n\text{-}\mathrm{Der}}(X)
		\]
		consists of a strict $(n-1)$-functor
		$h\colon B\to X_{\leq n-1}$ and, for every $s\in S$, an $n$-cell
		\[
		\varphi(s)\colon
		h\bigl(b^-(\kappa_{n-1,s})\bigr)
		\Rightarrow
		h\bigl(b^+(\kappa_{n-1,s})\bigr).
		\]
		The boundary equations in
		\eqref{eq:morphism-of-derivation-schemes} say exactly that these cells have
		the displayed sources and targets.  By the term construction in
		\cref{def:derivation-n-scheme}, equivalently by
		\cref{thm:derivation-scheme-coinserter}, these data extend uniquely to a
		strict $n$-functor
		\[
		\operatorname{Ext}_n(D)\longrightarrow X.
		\]
		Restriction to $B$ and evaluation on the generators $\xi_s$ give the inverse
		assignment.  The two constructions commute with morphisms of derivation
		schemes and with strict $n$-functors in $X$, and hence give the natural
		bijection required for the adjunction.
		
		When $n=1$, the same argument says that a functor from the free category on
		a graph is determined uniquely by its values on the objects and generating
		arrows.  Finally, applying the bijection to the identity of
		$\mathcal C_{n\text{-}\mathrm{Der}}(X)$ gives the stated counit.
	\end{proof}
	
	\subsection{Recursive computads}
	
	\begin{definition}[$n$-computad]
		\label{def:n-computad}
		A $0$-computad is a set, and a $1$-computad is a graph.  Inductively, for
		$n\geq2$, an $n$-computad $C$ consists of an $(n-1)$-computad $C_{<n}$, a set
		$C_n$ of generating $n$-cells, and parallel boundary functions
		\[
		\partial^-,\partial^+\colon C_n\longrightarrow
		\bigl(\F_{n-1}C_{<n}\bigr)_{n-1}.
		\]
		Here $\F_{n-1}C_{<n}$ is the free strict $(n-1)$-category constructed at the
		previous stage.  The boundary functions determine strict functors
		\begin{equation}
			\label{eq:higher-boundary-pair}
			b^-,b^+\colon
			C_n\copow\walking{n-1}\rightrightarrows
			\F_{n-1}C_{<n}
		\end{equation}
		that agree on all cells below dimension $n-1$.
	\end{definition}
	
	For a single generator $\alpha\in C_n$, write $x_\alpha$ and $y_\alpha$ for
	the common $(n-2)$-source and $(n-2)$-target of its two boundaries.  The
	cell being adjoined has the schematic form
	\[
	\begin{tikzcd}[column sep=huge]
		x_\alpha
		\ar[r,bend left=38,"\partial^-\alpha",""{name=upper,below}]
		\ar[r,bend right=38,"\partial^+\alpha"',""{name=lower,above}]&
		y_\alpha .
		\arrow[Rightarrow,from=upper,to=lower,"\alpha" description]
	\end{tikzcd}
	\]
	The single arrows in this picture represent $(n-1)$-cells and the double
	arrow represents the generating $n$-cell.  This is the usual globular
	computad, or polygraphic, datum
	\cite{Street1976,Street1995Higher,Burroni1993}.
	
	\begin{definition}[The free strict $n$-category on a computad]
		\label{def:ascelulassaolivres}
		The free strict $n$-category $\F_nC$ is defined together with computads by
		induction on $n$.  For $n=0$ it is the set $C_0$, regarded as a strict
		$0$-category, and for $n=1$ it is the free category on the graph $C$.
		Suppose that $n\geq2$ and put
		\[
		B=\F_{n-1}C_{<n}.
		\]
		The cells of $\F_nC$ through dimension $n-1$ are those of $B$.  Its
		$n$-cells are the well-typed formal expressions generated by
		\[
		1_u\quad(u\in B_{n-1}),\qquad
		\alpha\colon\partial^-\alpha\Rightarrow\partial^+\alpha
		\quad(\alpha\in C_n),
		\]
		under all compositions $\circ_j$, $0\leq j<n$, modulo the least congruence
		compatible with all globular source and target maps and containing all the
		structural equations of a strict globular $n$-category, including the
		compatibility of identity $n$-cells with lower-dimensional composition.  No
		further equations between expressions containing the generators in $C_n$
		are imposed.
	\end{definition}
	
	This term construction is the globular form of the iterated free
	enriched-category construction described above.  In particular, it does not
	depend upon the coinserter calculation that follows.
	
	\begin{proposition}[Assignment of generators]
		\label{prop:mandaosgeradores}
		Let $n\geq1$, let $C$ be an $n$-computad, and let $X$ be a strict
		$n$-category.  Write $X_{\leq n-1}$ for its underlying strict
		$(n-1)$-category.  To give a strict
		$n$-functor
		\[
		H\colon\F_nC\longrightarrow X
		\]
		is equivalent to giving a strict $(n-1)$-functor
		$h\colon\F_{n-1}C_{<n}\to X_{\leq n-1}$ and, for each
		$\alpha\in C_n$, an $n$-cell
		\[
		h(\partial^-\alpha)\ \xRightarrow{\ \xi_\alpha\ }\
		h(\partial^+\alpha)
		\]
		in $X$.  Under this correspondence $H$ restricts to $h$ and sends
		$\alpha$ to $\xi_\alpha$.
	\end{proposition}
	
	\begin{proof}
		Extend $h$ to formal $n$-cells by sending identities to identities, each
		generator $\alpha$ to $\xi_\alpha$, and each formal composite to the
		corresponding composite in $X$.  The strict $n$-category axioms in $X$ make the
		extension constant on the defining congruence of
		\cref{def:ascelulassaolivres}; hence it defines a strict $n$-functor
		$H\colon\F_nC\to X$.  Every value of $H$ on a formal cell is forced by these
		requirements, which proves uniqueness.  Restriction and evaluation on the
		generators give the inverse construction.
	\end{proof}
	
	Maps are defined by the same recursion.  A map of $0$-computads is a
	function.  Given $n\geq1$, a general map $f\colon C\to D$ consists of a map
	$f_{<n}\colon C_{<n}\to D_{<n}$ and, for each $\alpha\in C_n$, an $n$-cell
	of $\F_nD$ with source and target
	\[
	(\F_{n-1}f_{<n})(\partial^-\alpha)
	\quad\text{and}\quad
	(\F_{n-1}f_{<n})(\partial^+\alpha).
	\]
	At $n=1$ this says that a generating arrow may be sent to a path.  A map is
	\emph{generator-preserving} when its lower part is generator-preserving and
	a function $f_n\colon C_n\to D_n$ supplies the displayed cells as target
	generators; equivalently,
	\[
	(\F_{n-1}f_{<n})(\partial^\epsilon\alpha)
	=\partial^\epsilon(f_n\alpha)
	\qquad(\epsilon\in\{-,+\}).
	\]
	Composition is substitution of the chosen formal cells.  These
	maps form the categories $n\text{-}\Comp$ and
	$n\text{-}\Comp^{\mathrm{gen}}$, and evaluation on formal cells makes
	$\F_n\colon n\text{-}\Comp\to n\text{-}\Cat$ a functor.  The definition in
	\cref{def:two-computad} is the case $n=2$.
	
	This is the usual globular computad/polygraph data, expressed so that the
	top boundary is a genuine parallel pair.  In particular, the symbol
	$C_n\copow\walking{n-1}$ in \eqref{eq:higher-boundary-pair} emphasizes that
	the construction uses one independent copy of the walking cell for each
	generator.  Equivalently, in strict $(n-1)$-categories it is the product of
	$\walking{n-1}$ with the discrete $(n-1)$-category on $C_n$.
	
	\subsection{The free--underlying adjunction}
	\label{sec:free-underlying-computad}
	
	For completeness, we give the finite-dimensional truncation of the
	free--underlying polygraph adjunction used by M\'etayer
	\cite[Section~1]{Metayer2016}; the explicit recursion fixes the boundary
	conventions needed below.  For this adjunction we fix the
	generator-preserving maps.  This convention matters: the larger category
	$n\text{-}\Comp$ is useful for substitution, but a map there may send a
	generator to an arbitrary formal composite.  The right adjoint below instead
	records an actual cell together with chosen formal representatives of its
	boundary.
	
	We construct that right adjoint recursively.  Put
	$\mathcal U_0=1_{\Set}$.  Suppose that
	\[
	\F_{n-1}\colon (n-1)\text{-}\Comp^{\mathrm{gen}}
	\rightleftarrows (n-1)\text{-}\Cat\mathbin{:}\mathcal U_{n-1}
	\]
	has been constructed, and denote its counit by
	$\varepsilon_{n-1}$.  Let $X$ be a strict $n$-category and put
	$B=X_{\leq n-1}$.  The lower skeleton of $\mathcal U_nX$ is
	$\mathcal U_{n-1}B$.  Its generating $n$-cells are the triples
	\begin{equation}
		\label{eq:underlying-computad-generators}
		(u,v,\xi),
	\end{equation}
	where $u,v$ are parallel $(n-1)$-cells of
	$\F_{n-1}\mathcal U_{n-1}B$ and
	\[
	\xi\colon\varepsilon_{n-1,B}(u)
	\Rightarrow\varepsilon_{n-1,B}(v)
	\]
	is an $n$-cell of $X$.  The two boundaries of
	\eqref{eq:underlying-computad-generators} are $u$ and $v$.  If
	$H\colon X\to Y$ is a strict $n$-functor, $\mathcal U_nH$ applies
	$\mathcal U_{n-1}(H_{\leq n-1})$ to the lower skeleton and sends
	\[
	(u,v,\xi)\longmapsto
	\bigl(\F_{n-1}\mathcal U_{n-1}(H_{\leq n-1})(u),
	\F_{n-1}\mathcal U_{n-1}(H_{\leq n-1})(v),H(\xi)\bigr).
	\]
	Naturality of $\varepsilon_{n-1}$ makes the last triple well typed.
	
	For $n=1$, this construction is the ordinary underlying graph: its
	generating arrows are all the arrows of the category.  In higher dimensions,
	the formal representatives $u,v$ in
	\eqref{eq:underlying-computad-generators} are essential.  Replacing them
	only by their values in $B$ would lose the boundary data needed for a
	computad map.  This is the finite-dimensional form of the right adjoint
	constructed for polygraphs in \cite[Section~1]{Metayer2016}.
	
	\begin{theorem}[Free--underlying computad adjunction]
		\label{thm:free-underlying-computad-adjunction}
		For every $n\geq0$, the preceding construction gives an adjunction
		\[
		\F_n\colon n\text{-}\Comp^{\mathrm{gen}}
		\rightleftarrows n\text{-}\Cat\mathbin{:}\mathcal U_n .
		\]
		Its counit
		\[
		\varepsilon_{n,X}\colon\F_n\mathcal U_nX\longrightarrow X
		\]
		evaluates every formal cell in $X$.
	\end{theorem}
	
	\begin{proof}
		The assertion is immediate for $n=0$.  Assume it in dimension $n-1$.
		Let $C$ be an $n$-computad and $X$ a strict $n$-category.  A
		generator-preserving map
		\[
		f\colon C\longrightarrow\mathcal U_nX
		\]
		has a lower part
		$f_{<n}\colon C_{<n}\to\mathcal U_{n-1}(X_{\leq n-1})$.  By the inductive
		adjunction, this is equivalently a strict $(n-1)$-functor
		\[
		h\colon\F_{n-1}C_{<n}\longrightarrow X_{\leq n-1}.
		\]
		For each $\alpha\in C_n$, preservation of the two boundaries forces
		$f(\alpha)$ to have the form
		\[
		\bigl(\F_{n-1}f_{<n}(\partial^-\alpha),
		\F_{n-1}f_{<n}(\partial^+\alpha),\xi_\alpha\bigr).
		\]
		The inductive adjunction and its counit identify the evaluations of the first
		two entries with
		$h(\partial^-\alpha)$ and $h(\partial^+\alpha)$.  Thus the remaining entry is
		exactly an $n$-cell
		\[
		\xi_\alpha\colon h(\partial^-\alpha)
		\Rightarrow h(\partial^+\alpha)
		\]
		in $X$.  By \cref{prop:mandaosgeradores}, these data are equivalent to a
		unique strict $n$-functor $\F_nC\to X$.
		
		Conversely, such a strict functor supplies $h$ and the cells $\xi_\alpha$.
		Take the inductive adjunct $f_{<n}$ of $h$ and send $\alpha$ to the displayed
		triple.  The two constructions are inverse and natural in $C$ and $X$.
		This proves the adjunction.  Taking $C=\mathcal U_nX$ and the identity map
		under the resulting bijection gives the stated evaluation counit.
	\end{proof}
	
	We write
	\[
	T_n=\mathcal U_n\F_n
	\]
	for the induced monad on
	$n\text{-}\Comp^{\mathrm{gen}}$ and call it the computadic free strict
	$n$-category monad.
	
	\begin{theorem}[Computad coinserter]
		\label{thm:computad-coinserter}
		For a $2$-computad $C$, its free strict $2$-category $\F_2C$ is the
		coinserter in $2\text{-}\Cat_{\mathrm{icon}}$ of
		\[
		d_2(C_2\copow\2)
		\mathrel{\substack{\xrightarrow{d_2b^-}\\[-.5ex]
				\xrightarrow[d_2b^+]{}}}
		d_2(\freecat C_{\leq1}).
		\]
		Its underlying category in dimensions zero and one is
		$\freecat C_{\leq1}$, and its universal icon has component at
		$\alpha\in C_2$ equal to the generating $2$-cell $\alpha$.
	\end{theorem}
	
	\begin{proof}
		For a strict $2$-category $X$, a strict $2$-functor $\F_2C\to X$ is an
		assignment of objects and $1$-cells to the generating graph, together with a
		$2$-cell between the images of the two boundary paths for each
		$\alpha\in C_2$, by \cref{prop:mandaosgeradores}.  The first part is a strict
		$2$-functor
		$d_2\freecat C_{\leq1}\to X$; by \cref{lem:walking-cell-icons}, the second
		part is an icon between its two composites with the boundary functors.
		Transformations between such assignments satisfy exactly the compatibility
		condition in \eqref{eq:coinserter-universal-property}.  Indeed, an icon
		between the lower-skeleton functors extends componentwise to the only
		possible icon between the resulting functors out of $\F_2C$; its naturality
		on a generating $2$-cell is precisely that compatibility condition, and
		induction over formal vertical and horizontal composites proves naturality
		on every $2$-cell.  Restriction and extension are inverse and commute with
		postcomposition.  They therefore give the required natural isomorphism of
		mapping categories.
	\end{proof}
	
	The same argument is dimension-independent.
	
	\begin{theorem}[Free strict $n$-categories as coinserters]
		\label{thm:higher-coinserter}
		Let $n\geq2$ and let $C$ be an $n$-computad.  In
		$n\text{-}\Cat_{\mathrm{icon}}$, the free strict $n$-category $\F_nC$ is the
		coinserter of
		\begin{equation}
			\label{eq:n-coinserter}
			\begin{tikzcd}[column sep=large]
				d_n(C_n\copow\walking{n-1})
				\ar[r,shift left=.7ex,"d_nb^-"]
				\ar[r,shift right=.7ex,swap,"d_nb^+"]&
				d_n(\F_{n-1}C_{<n}) \ar[r,"j"]&
				\F_nC .
			\end{tikzcd}
		\end{equation}
		The coinserter map is the identity through dimension $n-1$, and the component
		of its universal $n$-icon at $\alpha\in C_n$ is the generating $n$-cell
		$\alpha\colon\partial^-\alpha\Rightarrow\partial^+\alpha$.
	\end{theorem}
	
	In particular, the universal $n$-icon in
	\eqref{eq:n-coinserter} is the displayed $2$-cell in
	$n\text{-}\Cat_{\mathrm{icon}}$:
	\[
	\begin{tikzcd}[column sep=huge]
		d_n(C_n\copow\walking{n-1})
		\ar[r,bend left=38,"jd_nb^-",""{name=upper,below}]
		\ar[r,bend right=38,"jd_nb^+"',""{name=lower,above}]&
		\F_nC .
		\arrow[Rightarrow,from=upper,to=lower,"\theta" description]
	\end{tikzcd}
	\qquad \theta_\alpha=\alpha .
	\]
	
	\begin{proof}
		Put
		\[
		A=C_n\copow\walking{n-1},\qquad
		B=\F_{n-1}C_{<n},
		\]
		and let $j\colon d_nB\to\F_nC$ be the inclusion of the lower skeleton.  The
		universal $n$-icon
		$\theta\colon jd_nb^-\Rightarrow jd_nb^+$ has component at the distinguished
		cell of the $\alpha$th summand equal to the generating $n$-cell $\alpha$.
		
		For a strict $n$-category $X$, let
		$\mathsf{Cone}_X(b^-,b^+)$ be the category whose objects are pairs
		\[
		(h,\xi),\qquad h\colon d_nB\to X,\qquad
		\xi\colon hd_nb^-\Rightarrow hd_nb^+,
		\]
		and whose morphisms $\gamma\colon(h,\xi)\to(k,\zeta)$ are $n$-icons
		$\gamma\colon h\Rightarrow k$ satisfying
		\begin{equation}
			\label{eq:higher-cone-morphism}
			(\gamma d_nb^+)\xi=\zeta(\gamma d_nb^-).
		\end{equation}
		Restriction defines a functor
		\begin{equation}
			\label{eq:higher-restriction-functor}
			\Phi_X\colon
			n\text{-}\Cat_{\mathrm{icon}}(\F_nC,X)
			\longrightarrow \mathsf{Cone}_X(b^-,b^+),
			\qquad H\longmapsto(Hj,H\theta).
		\end{equation}
		
		We construct its inverse.  Given $(h,\xi)$, the walking-cell calculation
		identifies $\xi$ with a family
		\[
		\xi_\alpha\colon h(\partial^-\alpha)
		\Rightarrow h(\partial^+\alpha)\qquad(\alpha\in C_n).
		\]
		The generator-assignment property
		\cref{prop:mandaosgeradores} gives a unique
		strict $n$-functor
		$\bar h\colon\F_nC\to X$ that restricts to $h$ and sends the generating cell
		$\alpha$ to $\xi_\alpha$.
		
		Now let $\gamma\colon(h,\xi)\to(k,\zeta)$ be a cone morphism.  Since $j$ is
		the identity through dimension $n-1$, the components of $\gamma$ give the only
		possible $n$-icon $\bar\gamma\colon\bar h\Rightarrow\bar k$.  At a generating
		$n$-cell $\alpha$, its naturality equation is
		\begin{equation}
			\label{eq:higher-generator-naturality}
			\gamma_{\partial^+\alpha}\circ_{n-1}\xi_\alpha
			=\zeta_\alpha\circ_{n-1}\gamma_{\partial^-\alpha},
		\end{equation}
		which is precisely the $\alpha$-component of
		\eqref{eq:higher-cone-morphism}.  Every $n$-cell of $\F_nC$ is built from
		generating $n$-cells and identities by the operations $\circ_j$.  Induction on
		such an expression, using
		\eqref{eq:n-icon-unit}--\eqref{eq:n-icon-composition} and interchange,
		therefore proves naturality of $\bar\gamma$ at every $n$-cell.  The same
		induction proves uniqueness.  Thus $(h,\xi)\mapsto\bar h$ and
		$\gamma\mapsto\bar\gamma$ define a functor $\Psi_X$ inverse to $\Phi_X$.
		
		Indeed, both composites are identities: a strict functor out of $\F_nC$ is
		determined by its restriction to $B$ and its values on the generating
		$n$-cells, while an $n$-icon is determined by its components on the unchanged
		$(n-1)$-skeleton and naturality on those generators.  Postcomposition with a
		strict $n$-functor commutes with restriction and extension, so the isomorphism
		\[
		n\text{-}\Cat_{\mathrm{icon}}(\F_nC,X)
		\cong\mathsf{Cone}_X(b^-,b^+)
		\]
		is natural in $X$.  This is
		\eqref{eq:coinserter-universal-property}, and it also proves the assertions
		about the lower skeleton and the universal components.  For $n=1$ the same
		statement is \cref{thm:free-category-coinserter}.
	\end{proof}
	
	\subsection{Relations and monadic presentation data}
	
	\begin{definition}[Presentation by an $(n+1)$-computad]
		\label{def:higher-presentation}
		Let $R$ be an $(n+1)$-computad with $n$-skeleton $C$.  It presents the strict
		$n$-category
		\[
		\Pp_nR=\F_nC/
		(\partial^-\rho=\partial^+\rho\mid\rho\in R_{n+1}),
		\]
		the coequalizer in $n\text{-}\Cat$ of the two boundary functors determined by
		the $(n+1)$-generators.
	\end{definition}
	
	\medskip
	\noindent\emph{The higher monadic comparison.}
	
	We now compare \cref{def:higher-presentation} with the monadic presentations
	of \cref{def:monadic-presentation}.  Let
	$\widehat{\walking{n}}$ denote the canonical globular $n$-computad with one
	generating $n$-cell, two generating $k$-cells in every dimension
	$0\leq k<n$, and the source and target maps of the $n$-globe.  For $n=0$ it
	has one $0$-generator.  By construction,
	\[
	\F_n\widehat{\walking{n}}\cong\walking{n}.
	\]
	If $S$ is a set,
	$S\copow\widehat{\walking{n}}$ denotes the coproduct of one copy of this
	computad for each $s\in S$.
	
	The adjunction of
	\cref{thm:free-underlying-computad-adjunction} has a comparison functor
	\[
	K_n\colon n\text{-}\Cat\longrightarrow
	\bigl(n\text{-}\Comp^{\mathrm{gen}}\bigr)^{T_n}.
	\]
	It sends $X$ to the $T_n$-algebra on $\mathcal U_nX$ whose structure map is
	$\mathcal U_n(\varepsilon_{n,X})$.  This comparison functor exists for every
	adjunction; no monadicity assertion is needed for the result below.
	
	\begin{theorem}[Computadic relations as $T_n$-relation data]
		\label{thm:higher-computadic-monadic-comparison}
		Let $n\geq1$, and let $R$ be an $(n+1)$-computad with $n$-skeleton $C$ and
		set $S=R_{n+1}$ of generating relations.  Its two boundary assignments
		determine a pair of Kleisli arrows for $T_n$,
		\begin{equation}
			\label{eq:higher-computadic-kleisli-pair}
			r^-,r^+\colon
			S\copow\widehat{\walking{n}}\rightrightarrows T_nC.
		\end{equation}
		The Kleisli extensions of
		\eqref{eq:higher-computadic-kleisli-pair} form a parallel pair of free
		$T_n$-algebra homomorphisms.  Under the canonical identifications
		\[
		F^{T_n}\bigl(S\copow\widehat{\walking{n}}\bigr)
		\cong K_n(S\copow\walking{n}),
		\qquad
		F^{T_n}C\cong K_n(\F_nC),
		\]
		the two Kleisli extensions are
		$K_n(\bar\partial^-)$ and $K_n(\bar\partial^+)$, where
		\[
		\bar\partial^-,\bar\partial^+\colon
		S\copow\walking{n}\rightrightarrows\F_nC
		\]
		are the boundary functors.  The coequalizer of
		$\bar\partial^-,\bar\partial^+$ in $n\text{-}\Cat$ is $\Pp_nR$.
		
		On objects, this identifies $(n+1)$-computads with the $T_n$-relation data
		whose relation object is a copower of
		$\widehat{\walking{n}}$ and whose two Kleisli arrows agree on every
		generator below dimension $n$.  With generator-preserving maps, the
		identification is functorial.
	\end{theorem}
	
	\begin{proof}
		For each $\rho\in S$, the two $n$-cells
		$\partial^-\rho$ and $\partial^+\rho$ of $\F_nC$ are parallel.  Assignment
		of the distinguished cell in the $\rho$th copy therefore gives strict
		$n$-functors
		\[
		\bar\partial^-,\bar\partial^+\colon
		S\copow\walking{n}\rightrightarrows\F_nC
		\]
		which agree on the boundary of every globe.  Under the adjunction
		$\F_n\dashv\mathcal U_n$, these are precisely the Kleisli arrows
		$r^-,r^+$ in \eqref{eq:higher-computadic-kleisli-pair}.
		
		The free $T_n$-algebra on an $n$-computad $A$ is
		$K_n(\F_nA)$.  Under this identification, the algebra homomorphisms
		$(r^-)^\sharp,(r^+)^\sharp$ obtained by Kleisli extension are the images
		under $K_n$ of the boundary functors:
		\[
		(r^\pm)^\sharp=K_n(\bar\partial^\pm).
		\]
		By \cref{def:higher-presentation}, the coequalizer in $n\text{-}\Cat$ of
		the corresponding strict functors
		$\bar\partial^-,\bar\partial^+$ is exactly $\Pp_nR$.  This proves the first
		assertion.
		
		Conversely, suppose that a $T_n$-relation datum has relation object
		$S\copow\widehat{\walking{n}}$ and that its two Kleisli arrows agree below
		dimension $n$.  Their adjuncts send the distinguished top generator in each
		summand to two $n$-cells of $\F_nC$ with the same globular boundary.  Taking
		these pairs as the boundaries of the elements of $S$ defines an
		$(n+1)$-computad with $n$-skeleton $C$.  The two constructions are inverse.
		A generator-preserving map carries the distinguished globes and their two
		boundary maps to the corresponding globes and boundary maps, so the
		identification respects identities and composition.
	\end{proof}
	
	The theorem identifies the relation data and the corresponding pair of free
	algebra homomorphisms.  It does not claim that $K_n$ is an equivalence or
	that it creates the displayed coequalizer.  If one invokes a monadicity
	theorem for precisely the generator-preserving category of computads fixed
	above, then the Eilenberg--Moore algebra presented by
	\eqref{eq:higher-computadic-kleisli-pair} is $K_n(\Pp_nR)$.  The
	coinserter and coequifier results of this paper do not require that additional
	step.
	
	The preceding quotient is a coequifier.  This is the second stage of the
	coinserter--coequifier construction introduced above.
	
	\begin{proposition}[Computadic relations as coequifiers]
		\label{prop:higher-relations-coequifier}
		Let $n\geq2$ and let $R$ be an $(n+1)$-computad with $n$-skeleton $C$.  The
		boundary $n$-cells of its $(n+1)$-generators determine parallel $n$-icons
		\[
		\alpha^-,\alpha^+\colon f\Rightarrow g
		\colon d_n(R_{n+1}\copow\walking{n-1})\longrightarrow \F_nC.
		\]
		The quotient $q\colon\F_nC\to\Pp_nR$ is their coequifier in
		$n\text{-}\Cat_{\mathrm{icon}}$.
	\end{proposition}
	
	\begin{proof}
		For each $\rho\in R_{n+1}$, the parallel $n$-cells
		$\partial^-\rho$ and $\partial^+\rho$ have a common $(n-1)$-source and a
		common $(n-1)$-target.  These common boundaries define strict functors
		\[
		a,b\colon R_{n+1}\copow\walking{n-1}
		\rightrightarrows \F_{n-1}C_{<n}.
		\]
		Composing $d_na$ and $d_nb$ with the inclusion
		$j\colon d_n\F_{n-1}C_{<n}\to\F_nC$ gives
		$f=jd_na$ and $g=jd_nb$.  By
		\cref{lem:walking-cell-icons}, the two families
		$(\partial^-\rho)_\rho$ and $(\partial^+\rho)_\rho$ determine the parallel
		$n$-icons $\alpha^-$ and $\alpha^+$.
		
		For every strict $n$-category $X$, a strict functor
		$h\colon\F_nC\to X$ satisfies $h\alpha^-=h\alpha^+$ if and only if
		\[
		h(\partial^-\rho)=h(\partial^+\rho)
		\qquad(\rho\in R_{n+1}).
		\]
		This is equivalent to factoring uniquely through the quotient $q$.  The
		quotient changes no cell below dimension $n$, so an $n$-icon between two such
		functors has the same components after factorization.  Its naturality on a
		quotient $n$-cell follows from naturality on any representative, and is
		independent of that representative because both functors respect the
		relations.  Thus factorization gives an isomorphism of the relevant mapping
		categories, natural in $X$, which is the coequifier universal property.
	\end{proof}
	
	\begin{proposition}
		\label{prop:higher-presentation-obstruction}
		Let $n\geq1$.  If a strict $n$-category $A$ is presented by an
		$(n+1)$-computad, its
		$(n-1)$-truncation is free on an $(n-1)$-computad.  Conversely, if
		\[
		A_{\leq n-1}\cong\F_{n-1}C_{<n}
		\]
		for some $(n-1)$-computad $C_{<n}$, then $A$ admits a presentation by an
		$(n+1)$-computad.
	\end{proposition}
	
	\begin{proof}
		Relations in dimension $n+1$ identify only $n$-cells, so the lower truncation
		of $\Pp_nR$ is the lower truncation of $\F_nC$, namely
		$\F_{n-1}C_{<n}$.  Conversely, choose an isomorphism
		between the $(n-1)$-truncation of $A$ and $\F_{n-1}C_{<n}$.  Extend
		$C_{<n}$ to an $n$-computad $C$ by
		taking one generating $n$-cell $[a]$ for every $n$-cell $a$ of $A$, with the
		boundary transported across the chosen isomorphism.  Evaluation gives a
		strict $n$-functor
		\[
		e\colon\F_nC\longrightarrow A
		\]
		that is the chosen isomorphism below dimension $n$ and sends $[a]$ to $a$.
		
		For every $n$-cell $z$ of $\F_nC$, adjoin an $(n+1)$-generator with boundaries
		\[
		z\ ,\ [e(z)].
		\]
		These are parallel, and every resulting equation lies in the kernel
		congruence of $e$.  Conversely, if $e(z)=e(z')=a$, then the two added equations
		identify both $z$ and $z'$ with $[a]$.  They therefore generate the entire
		kernel congruence of $e$.  Since $e$ is surjective on $n$-cells and is an
		isomorphism in lower dimensions, its quotient by this congruence is
		isomorphic to $A$.  The resulting $(n+1)$-computad is the required
		presentation.
	\end{proof}
	
	The obstruction explains why $(n+1)$-computads are not presentations of
	all strict $n$-categories.  General monadic presentations from
	\cref{thm:canonical-monadic-presentation} do present all of them, but can also
	impose relations in the lower skeleton.
	
	\subsection{The bicategorical dimension-two statement}
	
	A \emph{bicategory} has objects and hom-categories as a strict $2$-category
	does, but its composition is associative and unital only up to specified
	invertible associator and unitor $2$-cells, subject to the pentagon and
	triangle equations.  A \emph{pseudofunctor} preserves composition and
	identities up to specified invertible compositor and unit $2$-cells satisfying
	the corresponding coherence equations.  It is \emph{normal} when its unit
	comparison cells are identities; its compositor cells may still be
	nonidentity.
	
	If two normal pseudofunctors agree on objects, an \emph{icon} between them
	assigns a $2$-cell to every $1$-cell.  These cells are natural with respect
	to $2$-cells and compatible with the compositor cells of the two
	pseudofunctors.  This is the bicategorical version of
	\cref{def:n-icon} for $n=2$ \cite{Lack2010Icons}.  The same
	dimension-two argument now has the following form.
	
	\begin{definition}
		\label{def:cones-pseudofuntoriais}
		For a $2$-computad $C$, put
		$B=\freecat C_{\leq1}$ and $A=C_2\copow\2$, and let
		$b^-,b^+\colon A\rightrightarrows B$ be its boundary pair.  If $X$ is a
		bicategory, let
		$\mathsf{Cone}^{\mathrm{ps}}_X(C)$ be the category whose objects are pairs
		\[
		h\colon d_2B\longrightarrow X,\qquad
		\xi\colon hd_2b^-\Rightarrow hd_2b^+,
		\]
		where $h$ is a normal pseudofunctor and $\xi$ is an icon.  An arrow
		$(h,\xi)\to(k,\zeta)$ is an icon $\gamma\colon h\Rightarrow k$ satisfying
		\[
		(\gamma d_2b^+)\xi
		=\zeta(\gamma d_2b^-).
		\]
	\end{definition}
	
	\begin{proposition}[The strict representative in bicategories]
		\label{prop:bicategory-coinserter}
		Let $\mathsf{NHom}$ be the $2$-category of bicategories, normal
		pseudofunctors, and icons, and let $C$ be a $2$-computad.  For every
		bicategory $X$, restriction along the lower-skeleton inclusion induces an
		isomorphism of categories
		\begin{equation}
			\label{eq:coinserter-pseudofuntorial}
			\mathsf{NHom}(\F_2C,X)
			\ \cong\ \mathsf{Cone}^{\mathrm{ps}}_X(C),
		\end{equation}
		natural in $X$.  Consequently, the strict free $2$-category $\F_2C$,
		regarded as a bicategory, is the coinserter in $\mathsf{NHom}$ of the
		boundary pair in \cref{thm:computad-coinserter}.
	\end{proposition}
	
	\begin{proof}
		Write $B=\freecat C_{\leq1}$ and $A=C_2\copow\2$.  The inclusion
		$j\colon d_2B\to\F_2C$ and the icon whose components are the generating
		$2$-cells define a cone on the boundary pair.  For every bicategory $X$,
		restriction gives a functor from
		$\mathsf{NHom}(\F_2C,X)$ to $\mathsf{Cone}^{\mathrm{ps}}_X(C)$.
		
		Conversely, let such a cone consist of a normal pseudofunctor
		$h\colon d_2B\to X$ and an icon
		$\xi\colon hd_2b^-\Rightarrow hd_2b^+$.  Keep the objects and the images of
		all paths chosen by $h$, as well as the unit and compositor constraints of
		$h$, and send a generating $2$-cell $\alpha$ to the component $\xi_\alpha$.
		Here $\xi_\alpha$ denotes the component at the distinguished arrow in the
		$\alpha$th copy of $\2$; the other components of $\xi$ are forced by its
		unit axioms.
		
		Write the compositor of $h$ as
		\[
		\phi_{v,u}\colon h(v)h(u)\Longrightarrow h(vu).
		\]
		We now define the value of $\bar h$ on every formal $2$-cell.  It preserves
		identity and vertical composition.  If
		$\alpha\colon u\Rightarrow u'$ and
		$\beta\colon v\Rightarrow v'$ are horizontally composable, put
		\begin{equation}
			\label{eq:a-extensao-pseudofuntorial}
			\bar h(\beta\ast\alpha)
			=\phi_{v',u'}\,
			\bigl(\bar h(\beta)\ast\bar h(\alpha)\bigr)\,
			\phi_{v,u}^{-1}.
		\end{equation}
		The source of the right-hand side is $h(vu)$ and its target is
		$h(v'u')$, as required.
		
		We check that this recursive evaluation descends to the defining congruence
		of $\F_2C$.  Vertical associativity and units hold in each hom-category of
		$X$.  Horizontal associativity follows from the pentagon equation for the
		compositors of $h$, and the horizontal unit equations follow from normality
		and the triangle equations.  Naturality of $\phi$ gives interchange:
		applying \eqref{eq:a-extensao-pseudofuntorial} before or after vertical
		composition gives the same pasted $2$-cell.  These are precisely the
		relations used in
		\cref{def:ascelulassaolivres}; hence the evaluation is independent of a
		formal representative.  It defines the hom-functors of a normal
		pseudofunctor
		\[
		\bar h\colon\F_2C\longrightarrow X.
		\]
		Its unit and compositor constraints are those of $h$.
		Equation \eqref{eq:a-extensao-pseudofuntorial} is exactly their naturality
		with respect to the newly adjoined $2$-cells.  Thus $\bar h$ restricts to
		$h$ and sends each generating $2$-cell to the prescribed component of
		$\xi$.
		
		An icon between two cones already has a component at every path of $B$.  The
		cone-morphism equation is precisely its naturality equation at each generating
		$2$-cell of $C$.  Induction over vertical composites is immediate.
		For a horizontal composite, the compositor axiom for an icon and naturality
		of the two compositors extend the equation by
		\eqref{eq:a-extensao-pseudofuntorial}.  Thus it holds at every formal
		$2$-cell and is unaffected by the defining congruence.  The icon therefore
		extends uniquely to the two normal pseudofunctors just constructed.
		
		Restriction and extension are inverse on objects and morphisms: the values
		on all $1$-cells and the compositor constraints are already part of $h$, and
		the values on the generating $2$-cells determine the rest.  The two
		constructions commute with postcomposition, proving the naturality of
		\eqref{eq:coinserter-pseudofuntorial}.  This is the coinserter universal
		property.  The ambient $2$-category and its icon calculus are those of
		\cite{Lack2010Icons}.
	\end{proof}
	
	Thus the strict free $2$-category represents the indicated coinserter when
	it is tested by normal pseudofunctors and icons.

	\part{Topology, deficiency, and low-dimensional presentations}
	\section{Topological background}
	\label{sec:topology-background}
	
	Part~I expressed generators and relations internally to higher category
	theory.  We now realize the same generating and relation cells
	topologically.  In
	dimensions up to three, this topology controls both relations and relations
	among relations.  We recall the four ingredients used below: vertex-based
	fundamental groupoids, Euler characteristic, group homology, and crossed
	modules.
	
	\subsection{Fundamental groupoids and van Kampen}
	
	For a space $X$ and a subset $A\subseteq X$, let $\PiOne(X,A)$ denote the full
	subgroupoid of the fundamental groupoid of $X$ on the points of $A$.  Its
	arrows are endpoint-preserving homotopy classes of paths.  If $A$ meets every
	path component, the inclusion
	$\PiOne(X,A)\to\PiOne(X,X)$ is an equivalence.
	We abbreviate the full fundamental groupoid by
	$\PiOne(X)\defeq\PiOne(X,X)$.  A continuous map
	$f\colon X\to Y$ sends the class of a path $\gamma$ to the class of
	$f\gamma$; this defines the functor $\PiOne\colon\Top\to\Gpd$.
	
	The appropriate general form of van Kampen is a statement about homotopy
	colimits.  For a small category $A$, write $BA$ for the geometric
	realization of its nerve.  We write $\hocolim_J X$ for the Bousfield--Kan
	homotopy colimit of a diagram $X\colon J\to\Top$.  If
	$D\colon J\to\Gpd$ is a diagram of groupoids, its
	\emph{Grothendieck construction} $\int_JD$ has objects $(j,x)$ with
	$j\in J$ and $x\in D(j)$.  An arrow $(j,x)\to(k,y)$ is a pair
	$(u,\varphi)$ consisting of an arrow $u\colon j\to k$ in $J$ and an arrow
	$\varphi\colon D(u)(x)\to y$ in $D(k)$; composition is
	$(v,\psi)(u,\varphi)=(vu,\psi D(v)(\varphi))$.  We take
	\[
	\hocolim\nolimits_J^{\Gpd}D
	\defeq\gpcomp\!\left(\int_JD\right)
	\]
	as a concrete model, natural up to equivalence, for its homotopy colimit in
	the folk model structure on groupoids.  In this model structure the weak
	equivalences are equivalences of groupoids, the cofibrations are the functors
	injective on objects, and the fibrations are the functors which lift
	isomorphisms.  Equivalently, the displayed groupoid is the
	fundamental groupoid of the topological homotopy colimit of the classifying
	spaces $BD(j)$, by the Grothendieck-construction theorem
	\cite{Thomason1979}.  It therefore models its homotopy $1$-type, namely the
	Postnikov truncation which retains components and fundamental groups.
	
	\begin{theorem}[Farjoun's homotopy-colimit theorem]
		\label{thm:farjoun-vankampen}
		\label{fundamentalgrupoidenaoseperde}
		For every small category $J$ and every diagram $X\colon J\to\Top$, there is a
		natural equivalence of groupoids
		\[
		\PiOne\!\left(\hocolim_J X\right)
		\simeq
		\hocolim\nolimits_J^{\Gpd}(\PiOne\circ X).
		\]
		That is to say, the fundamental groupoid functor preserves homotopy
		colimits.
	\end{theorem}
	
	\begin{proof}
		Farjoun's commutation formula identifies the fundamental groupoid of a
		homotopy colimit with the groupoid homotopy colimit of the induced diagram
		\cite{Farjoun2004}.  By Thomason's Grothendieck-construction theorem, the
		classifying space of $\int_J(\PiOne\circ X)$ is the topological homotopy
		colimit of the corresponding classifying spaces \cite{Thomason1979};
		applying groupoid completion gives precisely the homotopy $1$-type used in
		our definition of $\hocolim_J^{\Gpd}$.  Thus Farjoun's formula has the
		displayed form.  The assertion is naturally an equivalence of groupoids, not
		an equality of chosen strict models.
	\end{proof}
	
	For the cellular computations below we need a strict pushout on a specified
	set of vertices.  Brown's theorem gives precisely this refinement.
	
	\begin{theorem}[Groupoid van Kampen]
		\label{thm:groupoid-vankampen}
		Suppose $X=U\cup V$, the interiors of $U$ and $V$ cover $X$, and
		$A\subseteq X$ meets every path component of $U$, $V$, and $U\cap V$.  Then
		the square
		\[
		\begin{tikzcd}
			\PiOne(U\cap V,A\cap U\cap V) \ar[r] \ar[d] &
			\PiOne(U,A\cap U) \ar[d] \\
			\PiOne(V,A\cap V) \ar[r] & \PiOne(X,A)
		\end{tikzcd}
		\]
		is a pushout in $\Gpd$.
	\end{theorem}
	
	\begin{proof}
		This is Brown's van Kampen theorem for fundamental groupoids
		\cite[Theorem~3.4]{Brown1967}; the multiple-basepoint condition is precisely
		what makes the statement valid without imposing connectedness on the
		intersection.
	\end{proof}
	
	The cell-attachment form follows by applying the theorem to a collar of the
	attaching circle.
	
	\begin{corollary}[Attaching $2$-cells]
		\label{cor:attaching-two-cells}
		Let $Y$ be a CW complex whose components meet $A\subseteq Y$, and form
		\[
		X=Y\cup_{\coprod_{r\in R}f_r}
		(R\copow D^2)
		\]
		from cellular based loops $f_r\colon(S^1,*)\to(Y,A)$.  Then
		$\PiOne(X,A)$ is the
		quotient of $\PiOne(Y,A)$ by the least groupoid congruence forcing every
		arrow represented by $f_r$ to be an identity.  An unbased attaching loop may
		instead be transported to $A$ along a chosen path.
	\end{corollary}
	
	\begin{proof}
		For a finite set of cells, attach them successively and apply
		\cref{thm:groupoid-vankampen}.  At the chosen basepoint, the map
		\[
		\PiOne(S^1,\{*\})\longrightarrow\PiOne(D^2,\{*\})
		\]
		sends the generator to the identity.  The van Kampen pushout therefore
		quotients by the groupoid congruence generated by the attaching loop.
		
		For arbitrary $R$, let $X_S$ denote the subcomplex obtained by attaching the
		cells indexed by a finite subset $S\subseteq R$.  Then
		$X=\colim_{S\subseteq_{\mathrm{fin}}R}X_S$.  The image of a path or a
		homotopy from the compact CW complexes $I$ or $I\times I$ is contained in a
		finite subcomplex of $X$
		\cite[Appendix, Proposition~A.1]{Hatcher2002}.  Hence
		\[
		\PiOne(X,A)\cong
		\colim_{S\subseteq_{\mathrm{fin}}R}\PiOne(X_S,A),
		\]
		and the finite calculation gives exactly the congruence generated by all
		the attaching loops.
	\end{proof}
	
	\subsection{CW complexes and Euler characteristic}
	
	A CW complex is built inductively by attaching $i$-disks to its
	$(i-1)$-skeleton along maps from their boundary spheres.  We write $X^i$
	for the $i$-skeleton.  A map of CW complexes is \emph{cellular} when it
	sends $X^i$ into the $i$-skeleton of the target for every $i$.
	
	For a finite CW complex $X$, let $c_i(X)$ be the number of $i$-cells.  Its
	Euler characteristic is
	\[
	\chi(X)=\sum_i(-1)^ic_i(X)
	=\sum_i(-1)^i\dim_k H_i(X;k)
	\]
	for every field $k$.  The second equality is the Euler--Poincar\'e formula for
	the cellular chain complex.  We put
	$b_i(X;k)=\dim_kH_i(X;k)$.  In particular, a finite connected
	$2$-complex satisfies
	\begin{equation}
		\label{eq:euler-two-complex}
		1-\chi(X)=b_1(X;k)-b_2(X;k),
	\end{equation}
	and a finite connected simply connected $3$-complex satisfies
	\begin{equation}
		\label{eq:euler-three-complex}
		\chi(X)=1+b_2(X;k)-b_3(X;k).
	\end{equation}
	
	If $A\subseteq X$ is a contractible CW subcomplex, then the inclusion is a
	cofibration and the quotient map $X\to X/A$ is a homotopy equivalence.  This
	is the collapse result used below \cite[Proposition~0.17]{Hatcher2002}.
	We also use the Hurewicz theorem in the form
	$\pi_2(X)\cong H_2(X;\Z)$ for a simply connected CW complex
	\cite[Theorem~4.32]{Hatcher2002}.
	
	\begin{example}[Presentation complexes need not split as wedges]
		The presentation $\langle a,b\mid aba^{-1}b^{-1}\rangle$ has presentation
		complex the torus, not a wedge of circles, disks, and spheres.  Cell counts
		still give $\chi=1-2+1=0$.  Our arguments use cellular chains and van Kampen,
		not a general wedge decomposition.
	\end{example}
	
	\subsection{Classifying spaces and the Hopf bound}
	
	For a small groupoid $\mathcal G$, put
	$H_*(\mathcal G;k)=H_*(B\mathcal G;k)$.  If $\mathcal G$ is connected and
	$x$ is an object, the inclusion of the isotropy group
	$\mathcal G(x,x)$ is an equivalence of groupoids and induces a homotopy
	equivalence
	\[
	B\mathcal G(x,x)\simeq B\mathcal G.
	\]
	
	Let $X$ be a connected CW complex and let $G=\pi_1(X,x)$.  The classifying
	map $X\to BG$ induces an isomorphism on $H_1(-;k)$, and the Hopf exact
	sequence over the integers is described in
	\cite[Chapter~II, \S5]{Brown1982}.  In the natural universal-coefficient
	diagram, the map is surjective on the tensor term in degree two and an
	isomorphism on the $\operatorname{Tor}$ term coming from first homology.
	The universal coefficient theorem for homology
	\cite[Section~3.A]{Hatcher2002} therefore gives a surjection
	\begin{equation}
		\label{eq:hopf-surjection}
		H_2(X;k)\twoheadrightarrow H_2(G;k).
	\end{equation}
	Consequently,
	\begin{equation}
		\label{eq:betti-inequalities}
		b_1(X;k)=b_1(G;k),\qquad b_2(X;k)\geq b_2(G;k).
	\end{equation}
	The same statements apply to a connected groupoid after choosing an object.
	We use rational coefficients when writing an unqualified Betti number.
	
	\subsection{Crossed modules and homotopy two-types}
	
	\begin{definition}
		A \emph{crossed module} is a group homomorphism
		$\partial\colon M\to P$ together with a left action of $P$ on $M$, written
		${}^pm$, satisfying
		\[
		\partial({}^pm)=p\partial(m)p^{-1},\qquad
		{}^{\partial(m)}n=mnm^{-1}.
		\]
		Its homotopy groups are
		$\pi_1=\operatorname{coker}\partial$ and
		$\pi_2=\ker\partial$.
	\end{definition}
	
	For a presentation complex, $P$ records loops in the $1$-skeleton, $M$
	records relative $2$-cells, and $\partial$ records their attaching words.  The
	action of $P$ changes the whiskering path, while the Peiffer identity is the
	algebraic form of interchange.  A $3$-cell therefore imposes an identity
	among relations \cite{BrownSalleh1999}.
	
	Crossed modules of groups are algebraic models of connected pointed homotopy
	$2$-types.  Brown and Spencer identify them with categorical groups
	\cite{BrownSpencer1976}.  Janelidze internalized this correspondence in the
	semi-abelian setting \cite{Janelidze2003Internal}; only the classical
	group-valued case is used here.  The comparison with strict $2$-groupoids
	and homotopy $2$-types is developed in \cite{Noohi2007} and
	\cite[Chapters~6 and~7]{BrownHigginsSivera2011}.  A strict $2$-groupoid is
	\emph{connected} when its underlying $1$-groupoid, consisting of its objects
	and $1$-cells, is connected.  After choosing an object,
	every connected strict $2$-groupoid is biequivalent to a one-object model;
	crossed modules over groupoids give the literal multiobject formulation.  In
	a chosen connected component, with $x\in X^1$, the fundamental crossed
	module of a CW pair $(X^2,X^1)$ is
	\[
	\partial\colon\pi_2(X^2,X^1,x)\longrightarrow\pi_1(X^1,x).
	\]
	Whitehead's free crossed-module theorem says that it is freely generated by
	the attaching maps of the $2$-cells.  Attaching a $3$-cell quotients by the
	corresponding identity among relations.  Colimit theorems for crossed modules
	and crossed complexes make these constructions compatible with cellular
	pushouts \cite{Whitehead1949,BrownHiggins1981,BrownHigginsSivera2011}.
	
	Write $B_2\mathcal Q$ for the classifying space of the crossed module
	corresponding to a connected strict $2$-groupoid $\mathcal Q$; equivalently,
	one may use a standard geometric realization of its $2$-nerve.  Its first
	two homotopy groups are the cokernel and kernel of the crossed-module
	boundary \cite[Proposition~2.6]{BrownHiggins1991}.
	
	For a connected strict $2$-groupoid $\mathcal Q$, the following elementary
	observation will be used repeatedly:
	\begin{equation}
		\label{eq:local-thin-pi2}
		\mathcal Q\text{ is locally thin}
		\quad\Longleftrightarrow\quad
		\pi_2(B_2\mathcal Q,x)=0.
	\end{equation}
	Indeed, two parallel $2$-cells differ by an automorphism of a $1$-cell, and
	whiskering identifies every such automorphism group with the automorphism
	group of an identity $1$-cell.
	
	\section{Presentation complexes}
	\label{sec:realization}
	
	The topological coinserter of \cref{def:topological-graph} realizes every
	generating arrow as an interval.  We now attach one disk for each generating
	relation and prove that the resulting CW complex realizes the presented
	groupoid.  We write
	\[
	\abs{G}=\FTopOne G
	\]
	for the graph realization.  By
	\eqref{eq:um-intervalo-para-cada-seta}, it is a $1$-dimensional CW complex
	with one $0$-cell per vertex and one $1$-cell per edge.
	
	\begin{proposition}
		\label{prop:graph-fundamental-groupoid}\label{issoeomesmogrupoide}
		Sending a signed edge to the corresponding oriented edge path induces a
		canonical isomorphism
		\[
		\freegpd G\ \cong\ \PiOne(\abs{G},G_0).
		\]
		In particular, $G$ is connected if and only if $\abs{G}$ is path-connected.
	\end{proposition}
	
	\begin{proof}
		The unit $G\to\CTopOne\abs{G}$ of
		\cref{prop:topological-graph-adjunction} gives the canonical edge-path map.
		The compact image of a path or endpoint-preserving homotopy is contained in
		a finite subcomplex \cite[Appendix, Proposition~A.1]{Hatcher2002}.  The
		edge-path theorem for a finite CW graph then shows that every path class
		between vertices has a finite edge-path representative
		\cite[Section~1.A]{Hatcher2002}.  Deleting a backtracking pair
		$e^{-1}e$ or $ee^{-1}$ is an endpoint-preserving homotopy, and the lift of a
		path to the universal covering tree of its connected component shows that
		two reduced edge paths with the same endpoints are homotopic only when they
		are equal.  Thus the
		topological edge-path normal form agrees with
		\cref{thm:reduced-word-normal-form}.  Composition and inverses agree on both
		sides.
	\end{proof}
	
	The full fundamental groupoid on all points of $\abs{G}$ is equivalent, not
	usually isomorphic, to the vertex-based groupoid in
	\cref{prop:graph-fundamental-groupoid}.
	
	\subsection{The presentation complex}
	
	Let $C=(G,R,\partial^-,\partial^+)$ be a groupoidal $2$-computad.  For each
	$\rho\in R$, the word
	\[
	w_\rho=(\partial^+\rho)^{-1}\partial^-\rho
	\]
	is a loop in $\freegpd G$.  Use \cref{prop:graph-fundamental-groupoid} to
	represent it by a cellular loop
	$\lambda_\rho\colon S^1\to\abs{G}$.  A reduced word of length $m$ traverses
	one edge on each of $m$ consecutive subintervals; the empty word gives a
	constant loop.  We use a based cellular parametrization of the displayed
	signed word.  Different such parametrizations are based homotopic and do not
	affect the construction below up to cellular homotopy.
	
	\begin{definition}[Presentation complex]
		\label{def:presentation-complex}
		The \emph{presentation complex} of $C$ is the $2$-dimensional CW complex
		$X_C$ defined by the pushout
		\begin{equation}
			\label{eq:presentation-complex-pushout}
			\begin{tikzcd}
				R\copow S^1 \ar[r,"{\coprod\lambda_\rho}"] \ar[d,hook] &
				{\lvert G\rvert} \ar[d] \\
				R\copow D^2 \ar[r] & X_C.
			\end{tikzcd}
		\end{equation}
		For an ordinary $2$-computad, the same definition is applied after its
		positive boundary paths are included in $\freegpd G$.
	\end{definition}
	
	The complex has one $0$-cell for each element of $G_0$, one $1$-cell for each
	element of $G_1$, and one $2$-cell for each element of $R$.
	
	\begin{theorem}[Fundamental groupoid of a presentation complex]
		\label{thm:presentation-complex}
		For every groupoidal $2$-computad $C$, the edge-path map induces a canonical
		isomorphism
		\[
		\Pp_{\mathrm{gpd}}C\ \cong\ \PiOne(X_C,G_0).
		\]
		For an ordinary $2$-computad $C$, this reads
		\[
		\gpcomp\Pp_1C\ \cong\ \PiOne(X_C,C_0).
		\]
	\end{theorem}
	
	\begin{proof}
		By \cref{prop:graph-fundamental-groupoid}, the vertex-based fundamental
		groupoid of the $1$-skeleton is $\freegpd G$.  By
		\cref{cor:attaching-two-cells}, attaching the disk indexed by $\rho$ forces
		the loop $w_\rho$ to be an identity.  This is equivalent to the equation
		$\partial^-\rho=\partial^+\rho$.  Hence the resulting fundamental groupoid is
		the quotient of $\freegpd G$ by exactly the groupoid congruence defining
		$\Pp_{\mathrm{gpd}}C$.  The ordinary case follows from
		$\Pp_{\mathrm{gpd}}C\cong\gpcomp\Pp_1C$.
	\end{proof}
	
	\begin{example}[The triangle, one last time]
		\label{ex:o-triangulo-vira-disco}
		For the relation in \cref{ex:o-triangulo-volta}, the attaching word is
		$c^{-1}ba$.  Hence $X_C$ is a disk, and
		$\PiOne(X_C,\{x,y,z\})$ is the connected thin groupoid on these three
		objects.  The deficiency is
		\[
		\abs{C_1}-\abs{C_0}+1-\abs{C_2}
		=3-3+1-1=0.
		\]
		The normalization used below is therefore already visible in the opening
		example.
	\end{example}
	
	\begin{corollary}[Classical presentation complexes]
		\label{cor:classical-presentation-complex}
		For a group presentation $P=\langle S\mid R\rangle$, the one-vertex case of
		$X_C$ is the classical presentation complex: its $1$-skeleton is
		$\bigvee_{s\in S}S^1$, and the relator disks are attached along the words in
		$R$.  Its fundamental group is the group presented by $P$.
	\end{corollary}
	
	\begin{proof}
		Apply \cref{def:presentation-complex} to the one-object groupoidal computad
		with loop generators $S$ and relators $R$.  The description of the cells is
		immediate, and the assertion about the fundamental group is the one-object
		case of \cref{thm:presentation-complex}.
	\end{proof}
	
	\begin{example}[The torus]
		The computad with one object, loop arrows $a,b$, and one relation $ab=ba$ has
		presentation complex
		\[
		S^1\vee S^1\cup_{aba^{-1}b^{-1}}D^2,
		\]
		where the displayed attaching word is chosen up to cyclic
		reparametrization.  This is the standard CW structure on the torus.  The
		example is why the
		$2$-cells in \eqref{eq:presentation-complex-pushout} cannot generally be
		split off as a wedge of disks or spheres.
	\end{example}
	
	\subsection{Functoriality}
	
	A generator-preserving computad map induces a cellular map of presentation
	complexes once the characteristic maps are chosen compatibly.  A general
	computad map may send a generating arrow to a path and a relation to a
	composite derivation; realizing it requires explicit choices of cellular
	representatives and disk fillings.  Different fillings need not be homotopic
	when the target has nonzero $\pi_2$.  The construction is therefore
	functorial on generator-preserving maps equipped with compatible
	characteristic maps, which is the functoriality used below.
	
	\begin{remark}
		Ordinary paths with equal-subdivision composition are not strictly
		associative, and reduced groupoid words are not natural under arrow
		identifications.  Moore paths provide a strict model when a path algebra is
		needed.  Here the cellular pushout
		\eqref{eq:presentation-complex-pushout} supplies the required universal
		property directly.
	\end{remark}
	
	\begin{proposition}[Euler count]
		\label{prop:presentation-euler-count}
		If $C$ is finite, then
		\[
		\chi(X_C)=\abs{C_0}-\abs{C_1}+\abs{C_2}.
		\]
		If $C$ is finite and connected, then over any field $k$,
		\[
		\abs{C_1}-\abs{C_0}+1-\abs{C_2}
		=b_1(X_C;k)-b_2(X_C;k).
		\]
	\end{proposition}
	
	\begin{proof}
		The first equality is the alternating cellular count.  The second is
		\eqref{eq:euler-two-complex} rewritten using the first.
	\end{proof}
	
	\section{Deficiency of groupoids}
	\label{sec:deficiency}
	
	Maximal-tree reduction forces the normalization by the number of objects.  It
	also turns the presentation score into $1-\chi$.
	
	\subsection{Finite and rank-finite presentations}
	
	\begin{definition}[Presentation deficiency]
		\label{def:presentation-deficiency}
		Let $C$ be a connected groupoidal $2$-computad with finite relation set and
		rank-finite generating graph.  Its \emph{presentation deficiency} is
		\[
		\defi(C)=\beta_1(C_{\leq1})-\abs{C_2}.
		\]
		If the graph is finite, this is
		\begin{equation}
			\label{eq:finite-presentation-deficiency}
			\defi(C)=\abs{C_1}-\abs{C_0}+1-\abs{C_2}
			=1-\chi(X_C).
		\end{equation}
	\end{definition}
	
	\begin{definition}[Groupoid deficiency]
		\label{def:groupoid-deficiency}
		A groupoid has a \emph{finite computadic presentation} if it is presented by
		a groupoidal $2$-computad with finite sets in dimensions $0$, $1$, and $2$.
		It has a \emph{rank-finite computadic presentation} if the generating graph
		is rank-finite and the relation set is finite.  For a connected groupoid
		$\mathcal G$, define
		\begin{align*}
			\defi_{\mathrm{fin}}(\mathcal G)
			&=\sup\{\defi(C)\mid C\text{ is a finite presentation of }\mathcal G\},\\
			\defi_{\mathrm{rk}}(\mathcal G)
			&=\sup\{\defi(C)\mid C\text{ is a rank-finite presentation of }
			\mathcal G\}.
		\end{align*}
		The supremum is taken in $\Z\cup\{+\infty\}$; the invariant is undefined if
		the indicated class of presentations is empty.
	\end{definition}
	
	Using a supremum avoids assuming that an optimal presentation is attained.
	The adjective \emph{finite} is reserved for finite cell sets.  Rank-finite
	presentations may have infinitely many objects and tree arrows, so they are
	not called finitely presented without qualification.
	
	\subsection{Reduction to an isotropy group}
	
	For a group $H$ admitting a finite presentation, write
	$\defi_{\mathrm{grp}}(H)$ for the supremum of
	$\abs S-\abs R$ over its finite group presentations
	$\langle S\mid R\rangle$ \cite{Johnson1990,Rapaport1973}.
	
	\begin{theorem}[Isotropy-group reduction]
		\label{thm:connected-groupoid-isotropy-deficiency}
		Let $\mathcal G$ be a connected groupoid and $x\in\Ob\mathcal G$.
		\begin{enumerate}
			\item If $\mathcal G$ has finitely many objects and admits a finite
			computadic presentation, then
			\[
			\defi_{\mathrm{fin}}(\mathcal G)
			=\defi_{\mathrm{grp}}\bigl(\mathcal G(x,x)\bigr).
			\]
			\item For an arbitrary object set, whenever a rank-finite presentation
			exists,
			\[
			\defi_{\mathrm{rk}}(\mathcal G)
			=\defi_{\mathrm{grp}}\bigl(\mathcal G(x,x)\bigr).
			\]
		\end{enumerate}
	\end{theorem}
	
	\begin{proof}
		Given a connected computadic presentation $C$, choose a maximal tree in its
		generating graph.  By \cref{thm:group-presentation-comparison}, collapsing the
		tree gives a presentation of the isotropy group with
		$\beta_1(C_{\leq1})$ generators and $\abs{C_2}$ relators.  It has the same
		deficiency as $C$.  Taking suprema gives the inequality from left to right.
		
		Conversely, choose one arrow $\tau_y\colon x\to y$ in $\mathcal G$ for each
		object $y$, with $\tau_x=1_x$.  Start with a group presentation of
		$\mathcal G(x,x)$ at $x$ and adjoin a tree arrow $\tau_y$ from $x$ to each
		$y\neq x$.  No new relation is required: every arrow $y\to z$ is uniquely
		represented as
		\[
		\tau_z h\tau_y^{-1}\qquad
		(h\in\mathcal G(x,x)).
		\]
		The resulting groupoidal computad presents $\mathcal G$ and has the same
		cycle rank and relation count as the original group presentation.  It is
		finite when the object set is finite, and rank-finite for an arbitrary object
		set.  Taking suprema gives the reverse inequalities.
	\end{proof}
	
	\begin{corollary}
		\label{cor:group-deficiency-agreement}
		For every finitely presented group $H$,
		\[
		\defi_{\mathrm{fin}}(\Sigma H)
		=\defi_{\mathrm{grp}}(H).
		\]
		Thus groupoid deficiency restricts to classical group deficiency in the
		one-object case.
	\end{corollary}
	
	\begin{proof}
		Apply \cref{thm:connected-groupoid-isotropy-deficiency} to the unique object
		of $\Sigma H$; its isotropy group is $H$.
	\end{proof}
	
	\subsection{The homological bound}
	
	Let $k$ be a field.  For a connected groupoid $\mathcal G$, write
	$b_i(\mathcal G;k)=\dim_kH_i(B\mathcal G;k)$.
	
	\begin{theorem}[Homological deficiency bound]
		\label{thm:homological-deficiency-bound}
		Let $C$ be a connected groupoidal $2$-computad with rank-finite generating
		graph and finite relation set, and suppose that it presents $\mathcal G$.
		Then the Betti numbers below are finite, and
		\begin{equation}
			\label{eq:homological-deficiency-bound}
			\defi(C)\leq b_1(\mathcal G;k)-b_2(\mathcal G;k).
		\end{equation}
		Consequently, whenever the indicated presentations exist,
		\begin{align*}
			\defi_{\mathrm{fin}}(\mathcal G)
			&\leq b_1(\mathcal G;k)-b_2(\mathcal G;k),\\
			\defi_{\mathrm{rk}}(\mathcal G)
			&\leq b_1(\mathcal G;k)-b_2(\mathcal G;k).
		\end{align*}
		Equivalently, every such presentation satisfies the relation lower bound
		\begin{equation}
			\label{eq:relation-lower-bound}
			\abs{C_2}\geq
			\beta_1(C_{\leq1})-b_1(\mathcal G;k)+b_2(\mathcal G;k).
		\end{equation}
	\end{theorem}
	
	\begin{proof}
		Choose a maximal tree and apply
		\cref{thm:group-presentation-comparison}.  Its collapse gives a finite
		one-object presentation $P$ of an isotropy group of $\mathcal G$, with
		$\beta_1(C_{\leq1})$ generators, $\abs{C_2}$ relators, and the same
		deficiency as $C$.  The finite presentation complex $X_P$ satisfies
		\[
		\defi(C)=1-\chi(X_P)=b_1(X_P;k)-b_2(X_P;k).
		\]
		Since $X_P$ is a finite $2$-complex, its first two Betti numbers are finite.
		The isomorphism on $H_1$ and the Hopf surjection on $H_2$ therefore show that
		$b_1(\mathcal G;k)$ and $b_2(\mathcal G;k)$ are finite as well.
		Its classifying map gives
		$b_1(X_P;k)=b_1(\mathcal G;k)$ and
		$b_2(X_P;k)\geq b_2(\mathcal G;k)$ by
		\eqref{eq:betti-inequalities}.  This proves
		\eqref{eq:homological-deficiency-bound}.  Taking the two suprema proves the
		groupoid statements, and rearranging the definition of $\defi(C)$ gives
		\eqref{eq:relation-lower-bound}.
	\end{proof}
	
	The theorem is an upper bound on deficiency and, equivalently, a lower bound
	on the number of relations.
	
	\begin{theorem}[Fixed-skeleton sharpness]
		\label{thm:cada-corda-paga-uma-relacao}
		Let $G$ be a connected rank-finite graph and let $T\subseteq G$ be a maximal
		tree.  Choose a root $r$, and let
		$\tau_x\colon r\to x$ be the unique reduced signed path in $T$.  For every chord
		$e\colon x\to y$, put
		\[
		p_e=\tau_y\tau_x^{-1}\colon x\to y.
		\]
		The groupoidal $2$-computad $C(G,T)$ with generating graph $G$ and one
		relation
		\begin{equation}
			\label{eq:cada-corda-volta-para-a-arvore}
			e=p_e
			\qquad(e\in G_1\setminus T_1)
		\end{equation}
		presents the thin groupoid freely generated by $T$.  It has
		\[
		\abs{C(G,T)_2}
		=\abs{G_1\setminus T_1}
		=\beta_1(G)
		\]
		and therefore $\defi(C(G,T))=0$.
		
		This relation count is minimal with the generating graph fixed: every
		groupoidal $2$-computad with generating graph $G$, a finite relation set, and
		a thin presented groupoid has at least $\beta_1(G)$ relations.  Thus the
		presentation \eqref{eq:cada-corda-volta-para-a-arvore} is sharp.
	\end{theorem}
	
	\begin{proof}
		Define a functor
		\[
		\freegpd G\longrightarrow\freegpd T
		\]
		to be the identity on every tree arrow and to send a chord $e$ to $p_e$.
		It respects the relations
		\eqref{eq:cada-corda-volta-para-a-arvore}, and hence factors through a
		functor
		\[
		\Pp_{\mathrm{gpd}}C(G,T)\longrightarrow\freegpd T.
		\]
		In the other direction, inclusion of the tree induces
		$\freegpd T\to\Pp_{\mathrm{gpd}}C(G,T)$.  The first composite is the identity
		on the tree generators.  The second is the identity on the tree generators
		and, by the defining relation, also on every chord.  The two functors are
		therefore inverse isomorphisms.  Since an unoriented tree has a unique reduced
		signed path between any two objects, $\freegpd T$ is thin.
		
		By the definition of $\beta_1$ for a rank-finite graph,
		$\abs{G_1\setminus T_1}=\beta_1(G)$.  There is one relation for each chord,
		and the asserted value of the deficiency follows.
		
		For minimality, let $C$ be any groupoidal presentation of the indicated kind.
		A connected thin groupoid is equivalent to the terminal groupoid, so its
		first and second Betti numbers vanish over every field.  Applying
		\cref{thm:homological-deficiency-bound} to $C$ gives
		\[
		\abs{C_2}\geq\beta_1(G).
		\]
		The displayed presentation attains this bound.
	\end{proof}
	
	\subsection{Thin groupoids}
	
	\begin{theorem}[Deficiency of a connected thin groupoid]
		\label{thm:thin-groupoid-deficiency}
		Let $\mathcal T$ be a connected thin groupoid.
		\begin{enumerate}
			\item If $\Ob\mathcal T$ is finite, then $\mathcal T$ has a finite computadic
			presentation and
			$\defi_{\mathrm{fin}}(\mathcal T)=0$.
			\item For an arbitrary small object set, $\mathcal T$ has a rank-finite
			computadic presentation and
			$\defi_{\mathrm{rk}}(\mathcal T)=0$.
		\end{enumerate}
		In both cases, a spanning tree with no $2$-generators is an efficient
		presentation: it attains the displayed supremum.
	\end{theorem}
	
	\begin{proof}
		Every isotropy group of a thin groupoid is trivial.  The trivial group has
		the empty presentation of deficiency zero, and the homological bound gives
		that no presentation has positive deficiency.  The result now follows from
		\cref{thm:connected-groupoid-isotropy-deficiency}.
		
		Concretely, take any tree on the object set.  Its free groupoid has exactly
		one arrow between every ordered pair of objects, so it is isomorphic to
		$\mathcal T$.  Its cycle rank and relation count are both zero.
	\end{proof}
	
	For a one-object groupoid, this specializes to the elementary equality
	$\defi_{\mathrm{grp}}(1)=0$.
	
	The definition depends on cycle rank, which is unchanged when a tree is
	adjoined.  Since an arbitrary monad carries no corresponding numerical size
	for its generator and relation objects, we do not attach a deficiency to the
	general monadic presentations of \cref{sec:presentations}.
	
	\section{Thinness and rewriting}
	\label{sec:thin-categories}
	
	Applying groupoid completion to the boundaries of an ordinary $2$-computad
	gives a groupoidal $2$-computad.  This transfers the homological bounds to
	categorical presentations.  Reflecting thinness back to the category is a
	separate problem, governed by faithful localization.
	
	\subsection{The groupoid-completion bound}
	
	For a finite connected ordinary $2$-computad $C$, put
	\[
	\defi_{\Cat}(C)=\abs{C_1}-\abs{C_0}+1-\abs{C_2}.
	\]
	For a category $A$ admitting such presentations, define
	$\defi_{\Cat}(A)$ as the supremum of these integers.  This is a presentation
	invariant of $A$ relative to ordinary computads; it is not asserted to equal
	the groupoid invariant.
	
	\begin{proposition}
		\label{prop:category-groupoid-deficiency-bound}
		Let $C$ be a finite connected ordinary $2$-computad presenting a connected
		category $A$.  Applying groupoid completion to its boundary paths presents
		$\gpcomp A$ and leaves the presentation deficiency unchanged.  Consequently,
		\[
		\defi_{\Cat}(A)\leq
		\defi_{\mathrm{fin}}(\gpcomp A)
		\]
		whenever the left-hand side is defined.  In particular, if $\gpcomp A$ is
		thin, every finite connected computadic presentation of $A$ has deficiency
		at most zero.
	\end{proposition}
	
	\begin{proof}
		Groupoid completion is a left adjoint and preserves the coequalizer that
		defines $\Pp_1C$.  No generators or relations are changed, so the cell count
		is unchanged.  The inequality follows because ordinary presentations form a
		subclass of groupoidal presentations.  If $\gpcomp A$ is thin, apply
		\cref{thm:thin-groupoid-deficiency}.
	\end{proof}
	
	\subsection{Faithful localization}
	
	\begin{theorem}[Faithful localization]
		\label{thm:faithful-localization}
		For a category $A$, the following are equivalent:
		\begin{enumerate}
			\item the unit $\eta_A\colon A\to U\gpcomp A$ is faithful;
			\item $A$ admits a faithful functor into a groupoid.
		\end{enumerate}
		If these conditions hold and $\gpcomp A$ is thin, then $A$ is thin.
	\end{theorem}
	
	\begin{proof}
		The unit itself proves $(1)\Rightarrow(2)$.  Conversely, any functor
		$F\colon A\to H$ into a groupoid factors uniquely through $\eta_A$.  If $F$
		is faithful and $\eta_A(f)=\eta_A(g)$, then $F(f)=F(g)$ and hence $f=g$.
		The final assertion is \cref{prop:localization-reflects-thin-under-faithfulness}.
	\end{proof}
	
	Johnstone gives intrinsic embedding criteria for categories in groupoids
	\cite{Johnstone2008}.  Two-sided cancellation is necessary for an embedding
	but is not sufficient, already for monoids.  Thus a proof that establishes
	only cancellation cannot invoke \cref{thm:faithful-localization}.
	
	The same distinction persists one dimension higher: a locally cancellative
	$2$-category need not have locally faithful groupoid completion.  We return
	to this point in \cref{sec:two-categories}.
	
	\subsection{A rewriting criterion}
	
	Localization is not the only way to prove thinness.  When the equations of
	an ordinary computad are oriented as rewrite rules, rewriting is often more
	informative.
	
	Let $C$ be an ordinary $2$-computad.  Orient each generating relation as a
	rewrite rule $\ell_\alpha\to r_\alpha$ between parallel paths.  A one-step
	rewrite may be applied inside arbitrary left and right path contexts.
	
	\begin{definition}
		The rewrite relation is \emph{terminating} if it has no infinite chain.  It is
		\emph{locally confluent} if every one-step fork is joinable: whenever
		$p\to q$ and $p\to r$, there are finite reduction sequences from $q$ and
		$r$ to a common path $s$,
		\[
		\begin{tikzcd}[column sep=large,row sep=large]
			& p \ar[dl] \ar[dr] & \\
			q \ar[dr,dashed,"{\to^*}"'] && r \ar[dl,dashed,"{\to^*}"] \\
			& s . &
		\end{tikzcd}
		\]
		A path is a \emph{normal form} if no rule applies to it.
	\end{definition}
	
	\begin{theorem}[Rewriting criterion for thinness]
		\label{thm:rewriting-thinness}
		Assume the oriented relations of $C$ are terminating and locally confluent.
		Then every arrow of $\Pp_1C$ has a unique normal-form representative.
		Moreover, $\Pp_1C$ is thin if and only if, for every ordered pair
		$(x,y)$, there is at most one normal path $x\to y$.
	\end{theorem}
	
	\begin{proof}
		Newman's lemma \cite{Newman1942} says that a terminating, locally confluent
		rewrite relation is confluent.  Hence every path reduces to a normal form and
		any two congruent paths reduce to the same normal form.  Conversely, every
		generating rewrite
		is an equation in $\Pp_1C$, so paths with the same normal form represent the
		same arrow.  Thus normal paths are in bijection with arrows of the presented
		category.  The last assertion is the definition of thinness.
	\end{proof}
	
	An overlapping one-step fork is \emph{critical} when it is minimal under the
	removal of common left and right contexts.  For the finite path-rewriting
	systems considered here, the critical-pair lemma reduces local confluence to
	joinability of the finitely many critical overlaps
	\cite[Chapter~2]{BookOtto1993}.  In \cref{sec:descent-example}, the same idea
	one dimension higher turns critical confluences into generating $3$-cells.
	
	\subsection{Constructive thin presentations}
	
	The \emph{preorder reflection} of a category $A$ is the thin category with
	the same objects as $A$ and with an arrow $x\to y$ precisely when
	$A(x,y)$ is nonempty.  A graph $G$ is \emph{acyclic} when $\freecat G$ has
	no nonidentity endomorphism.  We denote the preorder reflection of
	$\freecat G$ by $\operatorname{Reach}(G)$.  Thus
	$\operatorname{Reach}(G)$ has one arrow $x\to y$ precisely when there is a
	path from $x$ to $y$ in $G$.
	
	\begin{proposition}[Finite acyclic reachability]
		\label{prop:acyclic-reachability-presentation}
		Let $G$ be a finite acyclic graph.  For every reachable ordered pair
		$(x,y)$ choose one path $n_{x,y}$.  Add a relation
		$p=n_{x,y}$ for every other path $p\colon x\to y$.  This finite computad
		presents $\operatorname{Reach}(G)$.
	\end{proposition}
	
	\begin{proof}
		A finite acyclic graph has only finitely many paths, so the relation set is
		finite.  Let
		\[
		q\colon\freecat G\longrightarrow\operatorname{Reach}(G)
		\]
		be the canonical functor.  Its kernel congruence identifies precisely the
		paths with the same domain and codomain.  The displayed relations generate
		this congruence, since every such path is equated with the chosen
		$n_{x,y}$.  Hence the quotient of $\freecat G$ by these relations is
		canonically isomorphic to $\operatorname{Reach}(G)$.
	\end{proof}
	
	The construction lists all competing paths.  A convergent orientation can
	often replace it by relations associated with critical branchings.
	
	We next construct thin presentations from tree normal forms.
	
	We call a graph $T$ a \emph{tree} when its underlying unoriented multigraph
	is a tree.  In the next theorem, a \emph{chord} is an arrow adjoined to
	$T$.
	
	For a category $A$ and a set $S$ of arrows, $A[S^{-1}]$ denotes the
	universal category obtained by making every arrow of $S$ invertible.
	
	\begin{theorem}[Tree-path presentations]
		\label{thm:tree-path-presentations}
		Let $T$ be a tree and obtain a finite connected graph $G$ by adding
		chords.  Suppose first that each chord
		$e\colon x\to y$ is parallel to a path $p_e$ in $T$.  The relations
		\[
		e=p_e
		\]
		define a $2$-computad $C$ presenting $\freecat T$, hence a thin category.
		There is one relation per chord and $\defi_{\Cat}(C)=0$.
		
		More generally, suppose that every chord is either of the preceding forward
		kind or reverses a path in $T$.  Write each backward chord as
		$f\colon x_k\to x_0$ against a path
		\[
		x_0\xrightarrow{a_1}x_1\xrightarrow{a_2}\cdots
		\xrightarrow{a_k}x_k.
		\]
		Impose the $k+1$ cyclic relations
		\begin{equation}
			\label{eq:backward-chord-relations}
			(a_i\cdots a_1)f(a_k\cdots a_{i+1})=1_{x_i}
			\qquad(0\leq i\leq k),
		\end{equation}
		where an empty factor is omitted.  With these relations for every backward
		chord and $e=p_e$ for every forward chord, the presented category is thin.
	\end{theorem}
	
	\begin{proof}
		In the forward-only case, every chord is eliminated in favour of its tree
		path.  The quotient is $\freecat T$, and a tree has at most one path between
		two objects.  If there are $q$ chords, then
		$\beta_1(G)=q$ and there are $q$ relations.
		
		Eliminate the forward chords first.  Let $S$ be the set of tree arrows that
		occur in a path reversed by at least one backward chord.  For the
		cycle displayed in the statement and $1\leq i\leq k$, put
		\[
		b_i=(a_{i-1}\cdots a_1)f(a_k\cdots a_{i+1})
		\colon x_i\longrightarrow x_{i-1}.
		\]
		The cyclic relations at $i$ and $i-1$ say, respectively,
		$a_i b_i=1_{x_i}$ and $b_i a_i=1_{x_{i-1}}$.  Thus every $a_i\in S$ is
		invertible in the quotient.  Its inverse is unique, so the expressions
		obtained from different cycles agree.  The endpoint relations also say that
		$f$ is the inverse of the full tree path $a_k\cdots a_1$.
		
		Conversely, all cyclic relations hold in the localization
		$\freecat T[S^{-1}]$ when $f$ is sent to that inverse path.  The preceding
		formulas therefore give mutually inverse functors between the presented
		category and $\freecat T[S^{-1}]$.  A morphism of this localization has a
		reduced signed tree path in which only arrows of $S$ may occur backwards.
		The underlying unoriented tree has a unique simple path between two
		objects, so there is at most one such reduced path with prescribed
		endpoints.  Hence the localization, and therefore the presented category, is
		thin.
	\end{proof}
	
	The $k+1$ relations in \eqref{eq:backward-chord-relations} cannot in general
	be replaced by only the two endpoint relations.  Consider the cyclic graph
	\[
	x_0\xrightarrow{a_1}x_1\xrightarrow{a_2}x_2
	\xrightarrow{f}x_0.
	\]
	The two endpoint relations make
	$e=a_1fa_2\colon x_1\to x_1$ idempotent, since
	\[
	e^2=a_1f(a_2a_1f)a_2=e,
	\]
	but do not force $e$ to be the identity.  Indeed, take
	$X_0=X_2=\{*\}$ and $X_1=\{0,1\}$ in $\Set$, let
	$a_1(*)=0$, and let $a_2$ and $f$ be the unique maps between their
	respective singleton sets.  The endpoint relations hold, whereas $e$ is
	the constant map with value $0$.  Thus the theorem gives a sufficient
	presentation; it makes no claim of optimality.
	
	\section{Two-categorical presentations}
	\label{sec:two-categories}
	
	The top-dimensional construction of Part~I gives free strict
	$2$-categories.  A $3$-computad imposes equations between their $2$-cells.
	Topology enters through identities among relations, which are modeled
	algebraically by crossed modules.
	
	\subsection{Free two-categories}
	
	Let $C$ be a $2$-computad.  Assign weight one to every generating $2$-cell
	and weight zero to identity $2$-cells.  Extend weight additively under both
	vertical and horizontal composition.
	
	\begin{proposition}
		\label{prop:free-two-category-noninvertibility}
		Weight is well defined on the $2$-cells of $\F_2C$.  A vertically invertible
		$2$-cell of $\F_2C$ has weight zero; hence it is built without generating
		$2$-cells.  In particular, a free strict $2$-category with a nonidentity
		generating $2$-cell is not locally groupoidal.
	\end{proposition}
	
	\begin{proof}
		Strict associativity and unit laws only rebracket or delete identity cells,
		and interchange has the same multiset of generating $2$-cells on both sides.
		Thus total weight is invariant under the axioms of a strict $2$-category.  If
		$\alpha$ has a vertical inverse $\beta$, then
		$0=\operatorname{wt}(1)=\operatorname{wt}(\beta\alpha)
		=\operatorname{wt}(\beta)+\operatorname{wt}(\alpha)$, so both weights are
		zero.
	\end{proof}
	
	For a strict $2$-category $A$, applying groupoid completion to each
	hom-category gives a locally groupoidal $2$-category
	$\gpcomp_{\mathrm{loc}}A$ and a strict $2$-functor
	\[
	\eta_A^{\mathrm{loc}}\colon A\longrightarrow
	\gpcomp_{\mathrm{loc}}A.
	\]
	For categories $H,K$, the canonical comparison
	\[
	\gpcomp(H\times K)\longrightarrow\gpcomp H\times\gpcomp K
	\]
	is an isomorphism: both sides are obtained by inverting the arrows
	$(h,1)$ and $(1,k)$, and every arrow $(h,k)$ factors as their composite.
	The homwise composition functors therefore extend uniquely and still satisfy
	the strict associativity and unit equations.
	
	\begin{proposition}
		\label{prop:track-completion-adjunction}\label{otracklivre}
		Homwise groupoid completion is left adjoint to the inclusion of locally
		groupoidal strict $2$-categories into strict $2$-categories.  Consequently,
		$\gpcomp_{\mathrm{loc}}\F_2C$ is the free track $2$-category on an ordinary
		$2$-computad $C$.
	\end{proposition}
	
	\begin{proof}
		The adjunction between categories and groupoids extends independently on
		each hom-category.  The isomorphism
		$\gpcomp(H\times K)\cong\gpcomp H\times\gpcomp K$ makes the extended
		composition strict, and the homwise universal properties assemble into the
		required strict $2$-functor.  Applying the adjunction to the universal
		property of $\F_2C$ proves the final assertion.
	\end{proof}
	
	We say that a strict $2$-functor $F\colon A\to B$ is
	\emph{locally faithful} when every functor
	$F_{x,y}\colon A(x,y)\to B(Fx,Fy)$ is faithful.
	
	\begin{proposition}[Locally faithful localization]
		\label{prop:local-faithful-localization}
		If $\eta_A^{\mathrm{loc}}$ is locally faithful and
		$\gpcomp_{\mathrm{loc}}A$ is locally thin, then $A$ is locally thin.
		Local cancellation does not, by itself, imply that
		$\eta_A^{\mathrm{loc}}$ is locally faithful.
	\end{proposition}
	
	\begin{proof}
		Apply \cref{thm:faithful-localization} to each hom-category.  For the final
		sentence, let $M$ be a cancellative monoid which does not embed in a group;
		examples are given in \cite{EdwardesHeath2025}.  Form a strict
		$2$-category with objects $0,1$, one $1$-cell $f\colon0\to1$, and
		$\operatorname{End}_2(f)=M$; apart from identities, there are no other
		$1$- or $2$-cells.  This $2$-category is locally cancellative, but the
		completion unit on its hom-category $A(0,1)$ is not faithful.
	\end{proof}
	
	\subsection{Strict, locally groupoidal, and groupoidal completions}
	
	The same lower-dimensional generators determine three free constructions.
	If $C$ is an ordinary $2$-computad, let $C^{\mathrm{gpd}}$ denote the
	groupoidal $2$-computad obtained by applying groupoid completion to its
	boundary paths.  Write $\lambda_2$ for local groupoid completion,
	$\lambda_1$ for the subsequent inversion of the $1$-cells, and
	$\lambda_{1,2}$ for their composite.  Thus there is a commutative diagram of
	strict $2$-functors
	\begin{equation}
		\label{eq:tres-leituras-sem-tabela}
		\begin{tikzcd}[column sep=huge,row sep=large]
			\F_2C
			\ar[r,"\lambda_2"]
			\ar[dr,swap,"\lambda_{1,2}"]&
			\gpcomp_{\mathrm{loc}}\F_2C
			\ar[d,"\lambda_1"]\\
			&\F_2^{\mathrm{gpd}}C^{\mathrm{gpd}} .
		\end{tikzcd}
	\end{equation}
	The upper-left object is a strict $2$-category.  The upper-right object is a
	track $2$-category: every $2$-cell is invertible, but a $1$-cell need not be.
	The lower object is a strict $2$-groupoid.  Rewriting below is carried out in
	the track $2$-category, whereas the crossed-module construction is applied
	to the strict $2$-groupoid.
	
	\begin{definition}[$3$-computads of three types]
		\label{def:three-computad}
		For an ordinary $2$-computad $C$, write
		\[
		\F_2^{\mathrm{str}}C=\F_2C,\qquad
		\F_2^\top C=\gpcomp_{\mathrm{loc}}\F_2C.
		\]
		For a groupoidal $2$-computad $C$, let $\F_2^{\mathrm{gpd}}C$ be the free
		strict $2$-groupoid on its graph and generating $2$-cells, with the specified
		signed-word boundaries.  More precisely, strict $2$-functors
		$\F_2^{\mathrm{gpd}}C\to Q$ into a strict $2$-groupoid $Q$ are exactly the
		assignments of the generators of $C$ which respect those signed boundaries.
		This construction inverts both the $1$- and $2$-cells, whereas
		\cref{prop:track-completion-adjunction} inverts only the $2$-cells.  A strict
		or locally groupoidal $3$-computad consists of an ordinary $C$, a set $D_3$,
		and parallel boundaries
		\[
		\partial^-,\partial^+\colon D_3\longrightarrow
		\bigl(\F_2^\epsilon C\bigr)_2
		\qquad(\epsilon=\mathrm{str}\text{ or }\top).
		\]
		A groupoidal $3$-computad instead has a groupoidal $C$ and boundaries in
		$\F_2^{\mathrm{gpd}}C$.  In each case the presented object is
		\[
		\Pp_2^\epsilon D=\F_2^\epsilon C/
		(\partial^-\Gamma=\partial^+\Gamma\mid\Gamma\in D_3).
		\]
		It is, respectively, a strict $2$-category, a track $2$-category, or a
		strict $2$-groupoid.
	\end{definition}
	
	For $\epsilon=\mathrm{str}$, the quotient in
	\cref{def:three-computad} is the coequifier of
	\cref{prop:higher-relations-coequifier}.  By
	\cref{prop:higher-presentation-obstruction}, a strict $2$-category whose
	underlying $1$-category is not free cannot be presented by a strict
	$3$-computad of this form.  It may still have a general monadic
	presentation.
	
	\subsection{Cellular realization in dimension three}
	
	Let $C$ be a connected groupoidal $2$-computad, and let $X_C$ be its
	presentation complex as in \cref{def:presentation-complex}.  Fix a vertex
	$x$ and a maximal tree $T\subseteq X_C^1$, and put
	\[
	M_C=\pi_2(X_C,X_C^1,x),
	\qquad P=\pi_1(X_C^1,x).
	\]
	Whitehead's theorem identifies the fundamental crossed module
	\[
	\partial_C\colon M_C\longrightarrow P
	\]
	with the free crossed module on the generating $2$-cells of $C$, after
	transporting their attaching loops to $x$ along $T$.  If $\alpha,\beta$ are
	parallel $2$-cells of $\F_2^{\mathrm{gpd}}C$, their vertical difference
	$\beta^{-1}\alpha$, transported through the same one-object reduction,
	determines an element
	\[
	r(\alpha,\beta)\in\ker\partial_C.
	\]
	For a $3$-generator $\Gamma$, we denote
	\begin{equation}
		\label{eq:acertandoaesfera}
		r_\Gamma
		\defeq r(\partial^-\Gamma,\partial^+\Gamma).
	\end{equation}
	Since $X_C^1$ is a graph, $\pi_2(X_C^1,x)=0$; hence the exact homotopy
	sequence of the pair gives an isomorphism
	\begin{equation}
		\label{eq:aesferaeonucleo}
		\pi_2(X_C,x)\cong\ker\partial_C.
	\end{equation}
	
	We recall explicitly the one-object strict $2$-groupoid associated with a
	crossed module $\partial\colon M\to P$.  Its $1$-cells are the elements of
	$P$.  An element $m\in M$ labels a $2$-cell $u\Rightarrow v$ when
	$u=v\partial(m)$:
	\[
	\begin{tikzcd}[column sep=huge]
		*
		\ar[r,bend left=38,"u",""{name=upper,below}]
		\ar[r,bend right=38,"v"',""{name=lower,above}]&
		* .
		\arrow[Rightarrow,from=upper,to=lower,"m" description]
	\end{tikzcd}
	\qquad u=v\partial(m).
	\]
	Vertical composition is multiplication in right-to-left order, and
	whiskering is determined by the $P$-action.  Equivariance makes the
	boundaries of the whiskered cells correct, while the Peiffer identity is
	exactly the interchange equation.  For $\partial_C$, this convention gives
	\[
	\partial_C(m)=[v]^{-1}[u].
	\]
	Consequently, $\beta^{-1}\alpha$ corresponds to
	$m_\beta^{-1}m_\alpha$.  All products below are written with this
	convention.
	In particular,
	\[
	r_\Gamma=
	m_{\partial^+\Gamma}^{-1}m_{\partial^-\Gamma}.
	\]
	
	\begin{lemma}[Tree reduction and the free crossed module]
		\label{lem:arvore-cruzada}
		Collapsing $T$ gives a one-object reduction of
		$\F_2^{\mathrm{gpd}}C$ which, under the equivalence between categorical
		groups and crossed modules, corresponds to
		\[
		\partial_C\colon M_C\longrightarrow P.
		\]
		More explicitly, a formal $2$-cell $\alpha\colon u\Rightarrow v$ determines
		an element $m_\alpha\in M_C$ whose boundary is
		$[v]^{-1}[u]\in P$.  If
		$\alpha,\beta\colon u\Rightarrow v$ are parallel, then
		\[
		r(\alpha,\beta)=m_\beta^{-1}m_\alpha
		\in\ker\partial_C .
		\]
	\end{lemma}
	
	\begin{proof}
		After the tree has been collapsed, every generating arrow is a loop at
		$x$.  The universal property of the resulting one-object strict
		$2$-groupoid says that its $2$-cells are freely generated by the
		$2$-generators of $C$, with their prescribed boundary words and with
		whiskering by $P$.  Under the Brown--Spencer equivalence
		\cite{BrownSpencer1976}, this is precisely the universal property of the
		free crossed module on those boundary words.  Whitehead's free
		crossed-module theorem \cite[\S16]{Whitehead1949} identifies that crossed
		module with
		$\pi_2(X_C,X_C^1,x)\to\pi_1(X_C^1,x)$.  The stated formula for
		$m_\alpha$ records the oriented boundary of its representative disk.
		Parallel cells have the same boundary, so
		$m_\beta^{-1}m_\alpha$ lies in the kernel.  Vertical composition,
		whiskering, and interchange correspond, respectively, to multiplication,
		the $P$-action, and the Peiffer identity.
	\end{proof}
	
	\begin{definition}[Cellular groupoidal $3$-computad]
		A \emph{cellular groupoidal $3$-computad} consists of a datum
		$D=(C,D_3,\partial^-,\partial^+,(a_\Gamma)_{\Gamma\in D_3})$, where $C$ is a
		connected groupoidal $2$-computad, $D_3$ is a set,
		$\partial^\pm\colon D_3\to(\F_2^{\mathrm{gpd}}C)_2$ are parallel boundary
		maps, and every element has a specified cellular map
		\[
		a_\Gamma\colon S^2\longrightarrow X_C.
		\]
		Base the sphere at the source vertex of the common $1$-boundary and transport
		it to $x$ along the unique path in $T$.  Under
		\eqref{eq:aesferaeonucleo}, its homotopy class is required to correspond
		precisely to $r_\Gamma$.  Equivalently, one may obtain $a_\Gamma$ by gluing
		chosen cellular disk representatives of the two formal boundaries, with the
		second hemisphere oppositely oriented.  The realization of $D$ is
		\[
		X_D=X_C\cup_{(a_\Gamma)}(D_3\copow D^3).
		\]
		We write $D_i=C_i$ for $0\leq i\leq2$.  Changing $x$ or $T$ transports both
		$r_\Gamma$ and the based class of $a_\Gamma$, so the normalization is
		independent of the auxiliary choice.
	\end{definition}
	
	Thus the attaching map, rather than only its unbased homotopy class, belongs
	to the cellular datum.  The displayed normalization fixes the element of
	$\pi_2(X_C,x)$ represented by that map.
	
	\begin{definition}[$P$-normal closure]
		\label{def:fecho-normal-p}
		If $\partial\colon M\to P$ is a crossed module and
		$S\subseteq\ker\partial$, the \emph{$P$-normal closure}
		$\langle\!\langle S\rangle\!\rangle_P$ is the smallest normal subgroup of
		$M$ which contains $S$ and is stable under the action of $P$.  It is
		contained in $\ker\partial$: the kernel is both normal and $P$-stable.
		Consequently
		\[
		\langle\!\langle S\rangle\!\rangle_P\longrightarrow 1
		\]
		is a crossed submodule and the action and boundary descend to
		$M/\langle\!\langle S\rangle\!\rangle_P\to P$.
	\end{definition}
	
	For a cellular groupoidal $3$-computad $D$, define
	\[
	M_D=M_C/\langle\!\langle r_\Gamma\mid\Gamma\in D_3
	\rangle\!\rangle_P
	\quad\text{and}\quad
	\partial_D\colon M_D\longrightarrow P,
	\]
	where the denominator is the $P$-normal closure.  Finally, let
	$\mathcal Q_D$ be the one-object strict $2$-groupoid corresponding to the
	crossed module $\partial_D$.  A different root or maximal tree gives an
	equivalent connected $2$-groupoid.
	
	The quotient just defined is a morphism of crossed modules:
	\begin{equation}
		\label{eq:o-quociente-do-modulo-cruzado}
		\begin{tikzcd}[column sep=large,row sep=large]
			M_C \ar[r,"\partial_C"] \ar[d,two heads]&
			P \ar[d,equal]\\
			M_D \ar[r,"\partial_D"']&
			P .
		\end{tikzcd}
	\end{equation}
	Its kernel is the $P$-normal closure of the elements represented by the
	attaching spheres of the generating $3$-cells.
	
	\begin{definition}[Homotopical $2$-truncation]
		\label{def:homotopical-two-truncation}
		Let $K$ be a connected reduced, hence one-object, crossed complex with
		degree-one group $P$, degree-two group $M$, and boundary
		$\partial_3\colon K_3\to M$.  Its \emph{homotopical $2$-truncation} is the
		crossed module
		\[
		\tau_{\leq2}K
		\defeq
		\left(
		M/\langle\!\langle\operatorname{im}\partial_3
		\rangle\!\rangle_P\longrightarrow P
		\right).
		\]
		Thus the degree-three boundaries, together with all their $P$-translates,
		are made trivial in degree two.
	\end{definition}
	
	\begin{lemma}[Adjoining the $3$-cells]
		\label{lem:tres-celulas-quociente}
		Under the identification of \cref{lem:arvore-cruzada}, coequifying the two
		boundary $2$-cells of every $\Gamma\in D_3$ gives the crossed module
		\[
		\partial_D\colon M_D\longrightarrow P.
		\]
		The same crossed module is the $2$-truncation of the fundamental crossed
		complex of the skeletal filtration of $X_D$.
	\end{lemma}
	
	\begin{proof}
		In the one-object strict $2$-groupoid associated with
		$\partial_C$, the equation
		$\partial^-\Gamma=\partial^+\Gamma$ is equivalent to
		$r_\Gamma=1$.  Whiskering the equation by a loop applies the corresponding
		$P$-action.  Thus the universal quotient kills exactly
		$\langle\!\langle r_\Gamma\mid\Gamma\in D_3\rangle\!\rangle_P$, proving
		the first assertion.
		
		For the second, let $\Pi(X_D^\ast)$ denote the fundamental crossed complex
		of the skeletal filtration.  Its degree-two term is $M_C$, and the
		cellular boundary of the degree-three generator $\Gamma$ is
		$r_\Gamma$.  Its $2$-truncation is therefore
		\[
		\tau_{\leq2}\Pi(X_D^\ast)=
		\left(
		M_C/\langle\!\langle\operatorname{im}\partial_3\rangle\!\rangle_P
		\longrightarrow P
		\right)
		\cong(\partial_D\colon M_D\to P).
		\]
		Here the crossed-complex boundary and its $P$-translates give precisely the
		crossed submodule just described.  This is also the cellular instance of
		the crossed-complex colimit theorem
		\cite{BrownHiggins1981,BrownHigginsSivera2011}.
	\end{proof}
	
	A map of spaces is a \emph{$2$-equivalence} when it induces a bijection on
	path components and isomorphisms on $\pi_1$ and $\pi_2$ at every choice of
	basepoint.
	
	\begin{theorem}[Crossed-module comparison]
		\label{thm:crossed-module-comparison}
		Let $D$ be a connected cellular groupoidal $3$-computad.  The $1$- and
		$2$-generators define $\partial_C$, and the $3$-generators define the
		quotient $\partial_D$ above.  After choosing a root and maximal tree, the
		presented strict $2$-groupoid $\Pp_2^{\mathrm{gpd}}D$ is biequivalent to the
		one-object strict $2$-groupoid $\mathcal Q_D$ associated with
		$\partial_D$.  The classifying-space comparison induces a map, canonical up
		to the stated choices and homotopy,
		\[
		X_D\longrightarrow B_2\mathcal Q_D,
		\]
		which is a $2$-equivalence: the target has no homotopy above dimension two
		and the map induces isomorphisms
		\[
		\pi_i(B_2\mathcal Q_D,x)\cong\pi_i(X_D,x)
		\qquad(i=1,2).
		\]
		Thus $B_2\mathcal Q_D$ is a model for the homotopy $2$-type of $X_D$.
	\end{theorem}
	
	\begin{proof}
		\Cref{lem:arvore-cruzada,lem:tres-celulas-quociente} identify the
		one-object reduction of the presented strict $2$-groupoid with
		$\partial_D$; hence it is biequivalent to $\mathcal Q_D$.
		
		We spell out the topological comparison.  Collapse the chosen maximal tree
		$T$ and write $\overline X_D=X_D/T$.  Since $T$ is a contractible CW
		subcomplex, the quotient map $X_D\to\overline X_D$ is a homotopy
		equivalence.  The corresponding collapse of the pair
		$(X_D^2,X_D^1)$ identifies its degree-one and degree-two fundamental
		crossed-complex terms with $P$ and $M_C$.  Let
		$\Pi(\overline X_D^\ast)$ be the fundamental crossed
		complex of the skeletal filtration.  Its homotopical $2$-truncation is
		\[
		\tau_{\leq2}\Pi(\overline X_D^\ast)=
		\left(
		M_C/\langle\!\langle\operatorname{im}\partial_3
		\rangle\!\rangle_P\longrightarrow P
		\right)
		\cong(\partial_D\colon M_D\to P)
		\]
		by \cref{lem:tres-celulas-quociente}.
		
		The truncation morphism
		\[
		\Pi(\overline X_D^\ast)\longrightarrow
		\tau_{\leq2}\Pi(\overline X_D^\ast)
		\]
		determines, by the Brown--Higgins classification theorem
		\cite[Theorem~A*]{BrownHiggins1991}, a homotopy class of maps
		\[
		\overline X_D\longrightarrow
		B\tau_{\leq2}\Pi(\overline X_D^\ast)
		\simeq B_2\mathcal Q_D .
		\]
		Composing with $X_D\to\overline X_D$ gives the asserted comparison.  The
		construction is summarized by
		\[
		\begin{tikzcd}[column sep=large]
			X_D \ar[r,"\simeq"]&
			\overline X_D \ar[r,"\simeq_{\leq2}"]&
			B_2\mathcal Q_D .
		\end{tikzcd}
		\]
		The first arrow collapses the maximal tree and is a homotopy equivalence;
		the second is the map induced by crossed-complex truncation and is a
		$2$-equivalence.  Its homotopy
		class is canonical up to the choices already stated.  The proof of the
		$n=2$ case of
		\cite[Theorem~4.1]{BrownHiggins1991} shows that it induces isomorphisms on
		$\pi_1$ and $\pi_2$.  Explicitly, these groups are
		\[
		\pi_1=\operatorname{coker}\partial_D,\qquad
		\pi_2=\ker\partial_D
		\]
		by \cite[Proposition~2.6]{BrownHiggins1991}: attaching the $3$-cells
		quotients $\ker\partial_C$ by their $P$-normal boundary classes and does
		not change the cokernel.  The target is a $2$-type, so the comparison is
		a $2$-equivalence.
	\end{proof}
	
	\begin{corollary}[Homotopical obstruction to local thinness]
		\label{cor:pi2-local-thinness-obstruction}
		The presented strict $2$-groupoid $\Pp_2^{\mathrm{gpd}}D$ is locally thin if
		and only if $\pi_2(X_D,x)=0$.  Equivalently, the same is true of
		$\mathcal Q_D$.
	\end{corollary}
	
	\begin{proof}
		Combine \cref{thm:crossed-module-comparison} with
		\eqref{eq:local-thin-pi2}.  A biequivalence induces equivalences of
		hom-groupoids, and a groupoid equivalent to a thin groupoid is thin.
	\end{proof}
	
	\begin{corollary}[Homological $3$-cell lower bound]
		\label{cor:three-cell-bound}
		Let $D$ be a finite connected cellular groupoidal $3$-computad and suppose
		that the presented strict $2$-groupoid $\Pp_2^{\mathrm{gpd}}D$ is locally
		thin.  Put $G=\pi_1(X_D,x)$.  For every field $k$,
		\begin{equation}
			\label{eq:three-cell-homological-bound}
			\abs{D_3}\geq
			\abs{D_0}-\abs{D_1}+\abs{D_2}-1
			+b_1(G;k)-b_2(G;k).
		\end{equation}
		In particular, if $X_D$ is simply connected, then
		\begin{equation}
			\label{eq:three-cell-bound}
			\abs{D_3}\geq
			\abs{D_0}-\abs{D_1}+\abs{D_2}-1.
		\end{equation}
		Equivalently, in the simply connected case $\chi(X_D)\leq1$.
	\end{corollary}
	
	\begin{proof}
		Local thinness and \cref{cor:pi2-local-thinness-obstruction} give
		$\pi_2(X_D)=0$.  The relevant fragment of the Hopf exact sequence is
		\[
		\pi_2(X_D)\longrightarrow
		H_2(X_D;\Z)\longrightarrow H_2(G;\Z)\longrightarrow0 .
		\]
		It therefore identifies $H_2(X_D;\Z)$ with $H_2(G;\Z)$, while the
		classifying map also identifies integral first homology.  Naturality of the
		universal-coefficient short exact sequence now gives
		$H_i(X_D;k)\cong H_i(G;k)$ for $i=1,2$.  Hence
		\[
		\abs{D_0}-\abs{D_1}+\abs{D_2}-\abs{D_3}
		=1-b_1(G;k)+b_2(G;k)-b_3(X_D;k).
		\]
		Rearranging and using $b_3(X_D;k)\geq0$ gives
		\eqref{eq:three-cell-homological-bound}.  If $X_D$ is simply connected, then
		$G$ is trivial and its positive-degree Betti numbers vanish, giving
		\eqref{eq:three-cell-bound}.
	\end{proof}
	
	Local thinness alone does not imply simple connectedness; the general bound
	keeps the first two group-homology terms precisely for that reason.
	
	\subsection{Contractible collapses}
	
	\begin{proposition}[Contractible-spine collapse]
		\label{prop:contractible-spine-collapse}
		Let $D$ be a connected cellular groupoidal $3$-computad and
		$K\subseteq D$ a subcomputad whose realization $X_K$ is a contractible CW
		subcomplex of $X_D$.  Then
		\[
		X_D\longrightarrow X_D/X_K
		\]
		is a homotopy equivalence.  If $D/K$ is again a connected cellular
		groupoidal $3$-computad realizing the right-hand side, then the strict
		$2$-groupoids associated with $D$ and $D/K$ have equivalent homotopy
		$2$-types.  Local thinness is therefore preserved.
	\end{proposition}
	
	\begin{proof}
		The inclusion $X_K\hookrightarrow X_D$ is a cofibration because $X_K$ is a
		CW subcomplex.  Since $X_K$ is contractible, the quotient map
		$X_D\to X_D/X_K$ is a homotopy equivalence.  Apply
		\cref{thm:crossed-module-comparison} to the quotient realization.  Local
		thinness depends only on the vanishing of $\pi_2$ by
		\eqref{eq:local-thin-pi2}.
	\end{proof}
	
	The proposition requires a contractible realization and yields an equivalence
	of $2$-types.  A locally discrete thin sub-$2$-category need not satisfy the
	topological hypothesis, and even when it does the conclusion is not a strict
	isomorphism of $2$-categories.
	
	\subsection{The Eckmann--Hilton example}
	
	\begin{proposition}[Eckmann--Hilton case]
		\label{prop:eckmann-hilton-computad}
		Let $C$ have one $0$-cell, no nonidentity generating $1$-cells, and a set $S$
		of generating $2$-cells.  In the free strict $2$-groupoid on $C$, the
		$2$-endomorphisms of the identity $1$-cell form the free abelian group
		$\Z^{(S)}$.  The cellular realization is $\bigvee_{S}S^2$, and
		\[
		\pi_2\!\left(\bigvee_{S}S^2\right)\cong\Z^{(S)}.
		\]
	\end{proposition}
	
	\begin{proof}
		Vertical and horizontal composition have the same unit and satisfy
		interchange, so the Eckmann--Hilton argument makes them equal and
		commutative.  Freely inverting the $2$-generators therefore gives the free
		abelian group on $S$.  Topologically, the $1$-skeleton is a point and every
		$2$-cell is attached by the constant map, giving a wedge of $2$-spheres.  The
		Hurewicz theorem identifies its second homotopy group with its free second
		homology group.
	\end{proof}
	
	\section{The descent computad}
	\label{sec:descent-example}
	
	We finish with the finite descent computad that motivated the
	$3$-dimensional part of the paper.  We first apply local groupoid completion
	to its free strict $2$-category.  This inverts the $2$-cells but not the
	$1$-cells, and it is the setting for rewriting and the homotopy-basis
	argument.  We then invert the $1$-cells as well.  The resulting strict
	$2$-groupoid is the one to which the crossed-module comparison and the
	homological lower bound apply.
	
	The notation records the usual strict descent pattern.  The arrows
	$d^0,d^1$ compare the two pullbacks of a datum $d$, and $s^0$ is the
	degeneracy used for the unit condition.  The cell $\vartheta$ compares the
	two faces of $d$, while $n_0$ and $n_1$ are the structural degeneracy cells.
	At the next level, $\partial^0,\partial^1,\partial^2$ represent the three
	faces over a triple overlap, and the $\sigma_{ij}$ witness the corresponding
	cosimplicial face comparisons.  The confluence $I$ below is the unit
	condition on the descent comparison, and $A$ is its cocycle condition.  These
	are the two equations of descent data.  For the classical theory of descent
	morphisms, see Janelidze and Tholen \cite{JanelidzeTholen1997}; Street gives
	the categorical and combinatorial descent shapes used here
	\cite{Street2004Descent}.  Our notation follows our treatments by
	pseudo-Kan extensions, absolute Kan extensions, and semantic factorization
	\cite{Nunes2018PseudoKan,Nunes2021Absolute,Nunes2022Semantic}.
	
	\subsection{Generators}
	
	Let $C_{\mathrm{desc}}$ be the $2$-computad with objects
	$0,1,2,3$, generating $1$-cells
	\[
	\begin{tikzcd}[column sep=large]
		0 \ar[r,"d"] &
		1 \ar[r,,bend left=35, shift left=1.1ex,"d^0"]
		\ar[r,bend right=35,shift right=1.1ex,swap,"d^1"] &
		2 \ar[l,swap, "s^0" description]
		\ar[r,bend left=35,shift left=1.5ex,"\partial^0"]
		\ar[r,"\partial^1" description]
		\ar[r,bend right=35, shift right=1.5ex,swap,"\partial^2"] & 3,
	\end{tikzcd}
	\]
	and generating $2$-cells
	\begin{align*}
		n_0&\colon s^0d^0\Rightarrow1_1,&
		n_1&\colon 1_1\Rightarrow s^0d^1,&
		\vartheta&\colon d^1d\Rightarrow d^0d,\\
		\sigma_{01}&\colon\partial^1d^0\Rightarrow\partial^0d^0,&
		\sigma_{02}&\colon\partial^2d^0\Rightarrow\partial^0d^1,&
		\sigma_{12}&\colon\partial^2d^1\Rightarrow\partial^1d^1.
	\end{align*}
	We work in the free locally groupoidal $2$-category, so the generating
	$2$-cells may be used in either direction.  Orient them as rewrite rules on
	$1$-cells by using
	\[
	\begin{gathered}
		n_0\colon s^0d^0\longrightarrow1_1,\qquad
		\bar n_1=n_1^{-1}\colon s^0d^1\longrightarrow1_1,\qquad
		\vartheta\colon d^1d\longrightarrow d^0d,\\
		\sigma_{01}\colon\partial^1d^0\longrightarrow\partial^0d^0,\quad
		\sigma_{02}\colon\partial^2d^0\longrightarrow\partial^0d^1,\quad
		\sigma_{12}\colon\partial^2d^1\longrightarrow\partial^1d^1.
	\end{gathered}
	\]
	Equivalently, one may replace the generator $n_1$ by an oppositely oriented
	generator $\bar n_1$ before forming the free track $2$-category; the resulting
	free track $2$-category is isomorphic to the original one.  Thus the oriented
	system below is a $2$-polygraph to which the homotopy-basis theorem applies.
	
	\subsection{Termination and critical branchings}
	
	For the oriented $2$-computad above, a \emph{one-step reduction} is a
	generating rewrite rule placed in a left and a right $1$-cell context.  A
	\emph{local branching} is a pair of distinct one-step reductions with the
	same source.  Such a branching is disjoint when the two applications of
	rewrite rules occupy disjoint subwords.  Otherwise it is overlapping, and it
	is \emph{critical} when it is minimal among overlapping branchings under the
	addition of contexts.  A branching is \emph{confluent} when its two targets
	can be further reduced to a common target.
	
	\begin{lemma}
		\label{lem:descent-termination}
		The oriented rewrite system above is terminating.  It has exactly two
		overlapping critical branchings, with source words
		\[
		s^0d^1d
		\qquad\text{and}\qquad
		\partial^2d^1d.
		\]
	\end{lemma}
	
	\begin{proof}
		Order words first by length and then lexicographically from the left with
		\[
		\partial^2>\partial^1>\partial^0,\qquad d^1>d^0,
		\]
		extending this to any total order on the remaining generators.  The rules
		$n_0$ and $\bar n_1$ strictly shorten a word; the other rules preserve length
		and strictly decrease the lexicographic part.  The order is well founded and
		stable under left and right contexts, proving termination.
		
		Every left-hand side has length two.  The only suffix of one left-hand side
		that is the prefix of another is $d^1$: it produces the overlaps
		$s^0d^1d$ and $\partial^2d^1d$.  All other simultaneous rewrites are
		disjoint and commute by interchange.
	\end{proof}
	
	The first critical branching is completed by the two reductions
	\begin{equation}
		\label{eq:identity-branching}
		\begin{tikzcd}[column sep=large]
			s^0d^1d \ar[r,"{\bar n_1\ast1_d}"]
			\ar[d,swap,"{1_{s^0}\ast\vartheta}"] & d \\
			s^0d^0d \ar[ur,swap,"{n_0\ast1_d}"] &
		\end{tikzcd}
	\end{equation}
	Put
	\[
	u_I=\bar n_1\ast1_d,
	\qquad
	\ell_I=(n_0\ast1_d)(1_{s^0}\ast\vartheta).
	\]
	The \emph{identity descent $3$-cell} is the cell between the two displayed
	pasting paths:
	\begin{equation}
		\label{eq:a-celula-tres-da-identidade}
		\begin{tikzcd}[column sep=huge]
			s^0d^1d
			\ar[r,bend left=40,"u_I",""{name=upper,below}]
			\ar[r,bend right=40,"\ell_I"',""{name=lower,above}]&
			d
			\arrow[Rightarrow,from=upper,to=lower,"I" description]
		\end{tikzcd}.
	\end{equation}
	The second is completed by
	\begin{equation}
		\label{eq:associativity-branching}
		\begin{tikzcd}[column sep=huge,row sep=large]
			\partial^2d^1d
			\ar[r,"{\sigma_{12}\ast1_d}"]
			\ar[d,swap,"{1_{\partial^2}\ast\vartheta}"] &
			\partial^1d^1d
			\ar[r,"{1_{\partial^1}\ast\vartheta}"] &
			\partial^1d^0d
			\ar[d,"{\sigma_{01}\ast1_d}"] \\
			\partial^2d^0d
			\ar[r,swap,"{\sigma_{02}\ast1_d}"] &
			\partial^0d^1d
			\ar[r,swap,"{1_{\partial^0}\ast\vartheta}"] &
			\partial^0d^0d .
		\end{tikzcd}
	\end{equation}
	Put
	\[
	\begin{aligned}
		u_A&=(\sigma_{01}\ast1_d)
		(1_{\partial^1}\ast\vartheta)
		(\sigma_{12}\ast1_d),\\
		\ell_A&=(1_{\partial^0}\ast\vartheta)
		(\sigma_{02}\ast1_d)
		(1_{\partial^2}\ast\vartheta).
	\end{aligned}
	\]
	The \emph{associativity descent $3$-cell} is
	\begin{equation}
		\label{eq:a-celula-tres-da-associatividade}
		\begin{tikzcd}[column sep=huge]
			\partial^2d^1d
			\ar[r,bend left=40,"u_A",""{name=upper,below}]
			\ar[r,bend right=40,"\ell_A"',""{name=lower,above}]&
			\partial^0d^0d
			\arrow[Rightarrow,from=upper,to=lower,"A" description]
		\end{tikzcd}.
	\end{equation}
	By the \emph{full subcomputad} on a selected set of objects we
	mean the computad retaining those objects and exactly the displayed
	generators whose entire globular boundary lies on them.
	
	Let $C_{\mathrm{id}}$ be the full subcomputad of $C_{\mathrm{desc}}$ on
	$0,1,2$, with $2$-generators $n_0,n_1,\vartheta$.  Let
	$D_{\mathrm{id}}$ be the locally groupoidal $3$-computad obtained from it by
	adjoining $I$, and let $D_{\mathrm{desc}}$ be obtained from
	$C_{\mathrm{desc}}$ by adjoining $I$ and $A$.
	
	For a track $2$-category $Q$, its \emph{locally thin reflection} is the
	quotient which identifies every pair of parallel $2$-cells; it is universal
	among strict $2$-functors from $Q$ to locally thin track $2$-categories.  We
	denote the locally thin reflections of the free track $2$-categories on
	$C_{\mathrm{id}}$ and $C_{\mathrm{desc}}$ by
	\[
	\dot{\Delta}_{\mathrm{Str}_2}
	\qquad\text{and}\qquad
	\dot{\Delta}_{\mathrm{Str}},
	\]
	respectively.  This notation agrees with the strict descent shapes used in
	\cite{Nunes2018PseudoKan,Nunes2021Absolute}.
	
	For a $2$-computad $C$, a family of $3$-cells between parallel $2$-cells of
	the free track $2$-category on $C$ is a \emph{homotopy basis} when the
	quotient obtained by imposing their boundaries is locally thin.  Thus every
	$2$-sphere becomes filled in the quotient.  An oriented $2$-computad is
	\emph{convergent} when its rewriting relation is terminating and confluent.
	
	\begin{theorem}[Descent homotopy basis]
		\label{thm:descent-homotopy-basis}
		The set $\{I,A\}$ is a homotopy basis: after imposing these two equations,
		every pair of parallel $2$-cells is equal.  Equivalently,
		\[
		\Pp_2^\top D_{\mathrm{desc}}
		\cong\dot{\Delta}_{\mathrm{Str}}.
		\]
	\end{theorem}
	
	\begin{proof}
		By \cref{lem:descent-termination}, the system terminates and has only the two
		displayed critical branchings.  The reduction paths in
		\cref{eq:identity-branching,eq:associativity-branching} show that both are
		confluent, while $I$ and $A$ fill the resulting confluence spheres.  Every
		other local branching is either obtained from one of these by adding
		contexts or is disjoint; in the latter case its two reductions commute by
		interchange.  Hence the system is locally confluent, and termination gives
		confluence.  The polygraphic homotopy-basis theorem says that, for a
		convergent $2$-polygraph, one chosen confluence for every critical branching
		forms a homotopy basis
		\cite[Proposition~4.3.4]{GuiraudMalbos2009}.  This is the polygraphic form
		of Squier's homotopical theorem; its monoid antecedent is
		\cite{SquierOttoKobayashi1994}.
		Applying it to the two displayed confluences proves the claim.
	\end{proof}
	
	\begin{corollary}[The identity descent shape]
		\label{cor:identity-descent-homotopy-basis}
		The singleton $\{I\}$ is a homotopy basis for $C_{\mathrm{id}}$.  In
		particular,
		\[
		\Pp_2^\top D_{\mathrm{id}}
		\cong\dot{\Delta}_{\mathrm{Str}_2}.
		\]
	\end{corollary}
	
	\begin{proof}
		The oriented subsystem on $C_{\mathrm{id}}$ is terminating and has only the
		critical branching with source $s^0d^1d$.  Its chosen confluence is $I$.
		The same homotopy-basis theorem used above proves the result.
	\end{proof}
	
	We also consider the full subcomputad on the objects $1,2,3$, omitting the
	datum $d$ at the object $0$.  Let
	\[
	C_{\mathrm{mid}}
	=C_{\mathrm{desc}}\big|_{\{1,2,3\}}
	\]
	be the full subcomputad.  Thus its generating $1$-cells are
	$d^0,d^1,s^0,\partial^0,\partial^1,\partial^2$, and its generating $2$-cells
	are $n_0,n_1,\sigma_{01},\sigma_{02},\sigma_{12}$.  Let
	$\Delta_{\mathrm{Str}}$ denote the full sub-$2$-category of
	$\dot{\Delta}_{\mathrm{Str}}$ on the same three objects.
	
	\begin{theorem}[The intermediate strict descent shape]
		\label{thm:intermediate-strict-descent-shape}
		The canonical strict $2$-functor
		\[
		\F_2^\top C_{\mathrm{mid}}
		\longrightarrow\Delta_{\mathrm{Str}}
		\]
		is an isomorphism.  In particular, the free track $2$-category on
		$C_{\mathrm{mid}}$ is already locally thin and requires no generating
		$3$-cell.
	\end{theorem}
	
	\begin{proof}
		Restrict the orientation of \cref{lem:descent-termination} to
		$C_{\mathrm{mid}}$.  The resulting left-hand sides are
		\[
		s^0d^0,\qquad s^0d^1,\qquad
		\partial^1d^0,\qquad\partial^2d^0,\qquad\partial^2d^1.
		\]
		The same length-lexicographic order proves termination.  There is no
		overlapping critical branching.  Indeed, every displayed word has length two,
		its leftmost letter belongs to
		$\{s^0,\partial^1,\partial^2\}$, and its rightmost letter belongs to
		$\{d^0,d^1\}$; these two sets are disjoint.  A proper overlap would require
		the rightmost letter of one left-hand side to be the leftmost letter of
		another, and the five left-hand sides are themselves distinct.  Thus two
		distinct one-step reductions with a common source are disjoint, and their two
		composites agree by interchange.  The subsystem is therefore locally
		confluent, and hence convergent.
		
		Apply the polygraphic homotopy-basis theorem used in
		\cref{thm:descent-homotopy-basis}.  Since this subsystem has no critical
		branching, the empty family is a homotopy basis.  By definition, this says
		precisely that $\F_2^\top C_{\mathrm{mid}}$ is locally thin.
		
		It remains to identify the displayed functor.  Every $1$-cell generated by
		$C_{\mathrm{desc}}$ whose domain lies in $\{1,2,3\}$ uses only the generators
		of $C_{\mathrm{mid}}$: there is no generating arrow with codomain $0$.
		Likewise, a $2$-cell between such $1$-cells cannot involve $\vartheta$, whose
		boundary paths have domain $0$.  The same observation applies to the
		relations $I$ and $A$.  Whiskering either of them cannot move its domain away
		from $0$, since there is no path from $1$, $2$, or $3$ to $0$.  Passing from
		the free track $2$-category to
		$\Pp_2^\top D_{\mathrm{desc}}\cong\dot{\Delta}_{\mathrm{Str}}$ therefore adds
		neither cells nor equations in the full sub-$2$-category on $1,2,3$.
		Consequently the canonical functor is bijective on objects, $1$-cells, and
		$2$-cells, and is an isomorphism of strict $2$-categories.
	\end{proof}
	
	\begin{remark}
		\Cref{thm:intermediate-strict-descent-shape} concerns local groupoid
		completion: the $2$-cells, but not the $1$-cells, are inverted.  Full
		groupoid completion is obtained afterwards through the canonical
		localization
		\[
		\F_2^\top C_{\mathrm{mid}}\longrightarrow
		\F_2^{\mathrm{gpd}}C_{\mathrm{mid}}^{\mathrm{gpd}},
		\]
		which also inverts the generating $1$-cells.
	\end{remark}
	
	\subsection{Efficiency}
	
	Apply now full groupoid completion.  Write
	$D_{\mathrm{id}}^{\mathrm{gpd}}$ and
	$D_{\mathrm{desc}}^{\mathrm{gpd}}$ for the resulting groupoidal computads.
	Thus the $1$-cells are freely inverted, and the confluence spheres in
	\eqref{eq:identity-branching} and
	\eqref{eq:associativity-branching} are used as the cellular attaching maps.
	This is different from the local groupoid completion in
	\cref{thm:descent-homotopy-basis}.
	
	\begin{proposition}[Efficient descent presentations]
		\label{prop:efficient-descent-presentations}
		The groupoidal computads $D_{\mathrm{id}}^{\mathrm{gpd}}$ and
		$D_{\mathrm{desc}}^{\mathrm{gpd}}$ present locally thin strict
		$2$-groupoids.  Among finite cellular groupoidal $3$-computads with the same
		respective $2$-skeleta and a locally thin presented $2$-groupoid, the first
		needs at least one $3$-cell and the second at least two.  Thus the displayed
		groupoidal presentations attain the lower bound of
		\cref{cor:three-cell-bound}.
	\end{proposition}
	
	\begin{proof}
		Collapse the maximal tree with arrows $d,s^0$ for the identity presentation,
		and write $\overline X_{\mathrm{id}}$ for the resulting quotient CW complex.
		This collapse is a homotopy equivalence by
		\cref{prop:contractible-spine-collapse}.  Using the ordered cellular bases
		\[
		C_1(\overline X_{\mathrm{id}};\Z)\colon (x=d^0,y=d^1),\qquad
		C_2(\overline X_{\mathrm{id}};\Z)\colon(n_0,n_1,\vartheta),
		\]
		orient a $2$-cell so that its cellular boundary is the exponent vector of
		its source word minus that of its target word.  After the tree collapse the
		three attaching words therefore give
		\[
		\partial_2(n_0)=x,\qquad
		\partial_2(n_1)=-y,\qquad
		\partial_2(\vartheta)=y-x.
		\]
		Thus
		\[
		\partial_2=
		\begin{pmatrix}
			1&0&-1\\
			0&-1&1
		\end{pmatrix}.
		\]
		Indeed,
		\[
		\partial_2(a,b,c)=(a-c,-b+c),
		\]
		so $\ker\partial_2$ is freely generated by
		$v_I=(1,1,1)$.
		
		For clarity, we also calculate the cellular incidence of the $3$-cell.
		Write $e_\gamma$ for the standard basis vector of $C_2$ indexed by a
		$2$-generator $\gamma$.
		Incidence vectors add under the composition of $2$-cells.  Whiskering by a
		$1$-cell has incidence zero of its own and hence does not alter the
		incidence of the whiskered $2$-cell; inversion reverses its sign.  In
		\eqref{eq:identity-branching}, the upper path consequently has incidence
		$-e_{n_1}$, while the lower path has incidence
		$e_{\vartheta}+e_{n_0}$.  Choose the orientation of $I$ so that its
		cellular boundary is the incidence of the lower path minus that of the upper
		path.  Then
		\[
		\partial_3(I)
		=(e_{\vartheta}+e_{n_0})-(-e_{n_1})
		=v_I.
		\]
		The complex has no cells above dimension three.  Since
		$\im\partial_3=\ker\partial_2$ and $\partial_3\colon\Z[I]\to C_2$ is
		injective, it follows separately that
		\[
		H_2(\overline X_{\mathrm{id}};\Z)=0,
		\qquad
		H_3(\overline X_{\mathrm{id}};\Z)=0.
		\]
		The relations $n_0,n_1$ kill the two generators of $\pi_1$, so the
		realization is simply connected.  Hurewicz gives $\pi_2=0$, and
		\cref{cor:pi2-local-thinness-obstruction} shows that the presented strict
		$2$-groupoid is locally thin.
		
		For the full descent presentation, collapse the maximal tree with arrows
		$d,s^0,\partial^0$, and write $\overline X_{\mathrm{desc}}$ for the quotient.
		Use the ordered cellular bases
		\[
		\begin{split}
			C_1(\overline X_{\mathrm{desc}};\Z)&\colon
			(x=d^0,y=d^1,p=\partial^1,q=\partial^2),\\
			C_2(\overline X_{\mathrm{desc}};\Z)&\colon
			(n_0,n_1,\vartheta,\sigma_{01},\sigma_{02},\sigma_{12}).
		\end{split}
		\]
		With the same source-minus-target convention, the cellular boundaries are
		\[
		\begin{aligned}
			\partial_2(n_0)&=x,&
			\partial_2(n_1)&=-y,&
			\partial_2(\vartheta)&=y-x,\\
			\partial_2(\sigma_{01})&=p,&
			\partial_2(\sigma_{02})&=q+x-y,&
			\partial_2(\sigma_{12})&=q-p.
		\end{aligned}
		\]
		For example, the source and target of $\sigma_{02}$ become $qx$ and $y$
		after the tree collapse, while those of $\sigma_{12}$ become $qy$ and
		$py$.  Hence
		\[
		\partial_2=
		\begin{pmatrix}
			1&0&-1&0& 1& 0\\
			0&-1&1&0&-1& 0\\
			0&0&0&1& 0&-1\\
			0&0&0&0& 1& 1
		\end{pmatrix}.
		\]
		If $(a,b,c,r,s,t)$ denotes a vector in the displayed $2$-cell basis, the
		kernel equations are
		\[
		a-c+s=0,\qquad -b+c-s=0,\qquad r-t=0,\qquad s+t=0.
		\]
		Their general solution is
		\[
		(a,b,c,r,s,t)
		=\lambda(1,1,1,0,0,0)
		+\mu(0,0,-1,1,-1,1).
		\]
		Thus $\ker\partial_2$ is freely generated by
		\[
		v_I=(1,1,1,0,0,0),\qquad
		v_A=(0,0,-1,1,-1,1).
		\]
		The preceding incidence calculation for $I$ is unchanged.  For $A$, the
		upper path in \eqref{eq:associativity-branching} has incidence
		\[
		e_{\sigma_{12}}+e_{\vartheta}+e_{\sigma_{01}},
		\]
		whereas the lower path has incidence
		\[
		e_{\vartheta}+e_{\sigma_{02}}+e_{\vartheta}.
		\]
		Here again the whiskerings introduce no further $2$-cell incidence.
		Choose the orientation of $A$ so that its cellular boundary is the incidence
		of the upper path minus that of the lower path.  Then
		\[
		\partial_3(A)
		=(e_{\sigma_{12}}+e_{\vartheta}+e_{\sigma_{01}})
		-(2e_{\vartheta}+e_{\sigma_{02}})
		=v_A.
		\]
		Consequently $\im\partial_3=\ker\partial_2$.  Moreover, $v_I$ and $v_A$ are
		linearly independent, so
		\[
		\partial_3\colon\Z[I,A]\longrightarrow C_2
		\]
		is injective.  Since there are no higher cells, these two statements give,
		independently,
		\[
		H_2(\overline X_{\mathrm{desc}};\Z)=0,
		\qquad
		H_3(\overline X_{\mathrm{desc}};\Z)=0.
		\]
		The relations $n_0,n_1$ kill $x,y$, and
		$\sigma_{01},\sigma_{02}$ then kill $p,q$, so $\pi_1=0$.
		Hurewicz and \cref{cor:pi2-local-thinness-obstruction} prove local thinness.
		
		Finally, the respective cell counts are
		\[
		(|D_0|,|D_1|,|D_2|)=(3,4,3)
		\quad\text{and}\quad
		(|D_0|,|D_1|,|D_2|)=(4,7,6).
		\]
		Every groupoidal $3$-computad with either of these fixed $2$-skeleta has
		simply connected realization: the $2$-cells already kill $\pi_1$, and
		adjoining $3$-cells does not change it.  The bound
		\eqref{eq:three-cell-bound} is therefore one in the first case and two in the
		second.  The cells $I$ and $I,A$, respectively, attain it.
	\end{proof}
	
	The efficiency statement is numerical and relative to the fixed
	$2$-skeleton.  A different lower skeleton may have different numbers of
	generating cells in dimensions $0,1,2$; the proposition makes no global
	minimality claim.
	
	\section{Final remarks}
	\label{sec:conclusion}\label{voltandoaocomeco}
	
	The main categorical result identifies the free extension of a strict
	$(n-1)$-category by $n$-cells as the coinserter of its boundary pair.  More
	precisely, this construction is left adjoint to the underlying
	derivation-scheme functor.  When the lower skeleton is freely generated by a
	computad, the extension is $\F_nC$; an $(n+1)$-computad then imposes
	equations between its $n$-cells by a coequifier.  On generator-preserving
	computad maps, the free functor also has the recursively defined right
	adjoint of \cref{thm:free-underlying-computad-adjunction}.  In dimension one,
	the comparison with monadic presentations shows that object equations can
	be imposed first and the remaining relation datum is computadic.  The path
	normal form gives, in addition, the classification of the nonempty totally
	ordered free categories as the finite ordinals, $\N$, $\N^{\op}$, and $\Z$.
	
	The interval-weighted topological coinserter gives the $1$-skeleton of the
	presentation complex.  After the relation disks are attached, its
	fundamental groupoid on the vertices is canonically the groupoid presented
	by the computad.  Collapsing a maximal tree then gives both the correction
	by the number of objects in groupoid deficiency and the comparison with the
	deficiency of an isotropy group.  On a fixed connected rank-finite graph,
	one relation for each arrow outside a maximal tree presents a thin groupoid,
	and the homological bound shows that fewer relations cannot suffice.
	Crossed modules give the analogous comparison for $3$-cells and identities
	among relations.
	
	For the descent computad, the identity and associativity confluences form a
	homotopy basis in the free track $2$-category.  After full groupoid
	completion, the same confluences attain the homological lower bound relative
	to the two fixed $2$-skeleta.  The full subcomputad on the objects $1,2,3$
	has no critical branching, and its free track $2$-category is already
	locally thin.
	
	\subsection*{Historical position and chronology}
	
	The two-dimensional universal properties used here belong to the calculus of
	weighted colimits developed by Street and Kelly
	\cite{Street1976,Kelly1982,Kelly1989}.  In the same 1976 paper, Street
	introduced computads; his later accounts treat higher computads and free
	strict higher categories \cite{Street1995Higher,Street1996}.  Blackwell,
	Kelly and Power placed the
	relevant algebraic constructions in two-dimensional monad theory
	\cite{BlackwellKellyPower1989}.  Johnson gave a combinatorial proof of the
	$n$-categorical pasting theorem, and Power subsequently generalized the
	class of pasting schemes and expressed Street's free strict
	$n$-categories by left adjoints \cite{Johnson1989,Power1991}.  For a
	retrospective account of this circle of ideas, see Street
	\cite{Street2010Conspectus}.  Within this setting, the contribution of the
	present paper is to identify a single strict cellular extension as a
	coinserter, with the adjoined cells appearing as the components of its
	universal $2$-cell, and to apply this construction to presentations,
	topological realizations, and deficiency.
	
	Burroni introduced polygraphs as a higher-dimensional language for rewriting
	\cite{Burroni1993}.  Batanin defined $A$-computads for finitary monads $A$
	on globular sets, encompassing Street's strict computads and
	higher-operadic examples \cite{Batanin1998}.  M\'etayer developed
	polygraphic resolutions \cite{Metayer2003}.  Lafont subsequently synthesized
	the interaction of polygraphic rewriting with homotopy and homology
	\cite{Lafont2007}, and M\'etayer later gave a direct proof that strict
	$\omega$-categories are monadic over polygraphs \cite{Metayer2016}.  The
	results proved here concern strict $n$-categories.  In dimension two the
	transformations used here are Lack's icons; the $n$-icons are their strict
	top-dimensional counterpart in the iterated-icon setting
	\cite{Lack2010Icons,ChengGurski2014}.
	
	The mathematical backbone of this paper was completed during my PhD in
	mathematics in the joint doctoral programme of the Universities of Coimbra
	and University of Porto.  It appeared as DMUC Preprint 17--20 and as arXiv:1704.04474 on
	14 April 2017 \cite{Nunes2017Preprint}, and became Chapter~2 of the
	dissertation submitted in September 2017 and defended in January 2018
	\cite{Nunes2017Thesis}.  The coinserter constructions, the graph and
	$2$-computad realizations, the ordered-free-category classification, the
	groupoid-deficiency calculations, and the strict descent calculations were
	already present in that version.  The descent notation belongs to the
	contemporaneous work on pseudo-Kan extensions and to its later formulations
	by absolute Kan extensions and semantic factorization
	\cite{Nunes2018PseudoKan,Nunes2021Absolute,Nunes2022Semantic}.
	
	After the 2017 preprint, strict and geometric computads were developed
	further through non-unital polygraphs, constructible realizations, discrete
	Conduch\'e $\omega$-functors, and comparisons of pasting formalisms
	\cite{Henry2019,Hadzihasanovic2019,Guetta2020,Forest2022}; the monograph
	\cite{AraEtAl2025} now gives a broad account of strict polygraphs.  In weak
	and generalized settings, later work develops inductive computads,
	generalized signatures, hom $\omega$-categories, and specified invertible
	generators
	\cite{DeanEtAl2022,Markakis2024,BenjaminMarkakis2025,
		BenjaminChampinMarkakis2026}.
	
	The present revision makes the distinction between finite and rank-finite
	presentations explicit and proves the isotropy, homological, and
	crossed-module comparisons in the stated generality.
	
	\subsection*{Acknowledgements}
	
	I thank Maria Manuel Clementino for her careful supervision of the doctoral
	work from which this paper arose.  I thank Martin Hyland for his careful
	reading and comments as an examiner of the thesis.  The anonymous referee of
	the original submission gave a detailed report which led to substantial
	corrections and the present reorganization.
	
	The preparation of the present version was supported by the Deutsche
	Forschungsgemeinschaft (DFG, German Research Foundation) through the project
	Higher-Order Monad-based Programming and Reasoning (HOMBRe), project number
	501369690, and by the Centre for Mathematics of the University of Coimbra
	(CMUC; doi:10.54499/UID/00324/2025), funded by the Funda\c{c}\~ao para a
	Ci\^encia e a Tecnologia through grants UID/00324/2025 and
	UID/PRR/00324/2025.

	\providecommand{\bysame}{\leavevmode\hbox to3em{\hrulefill}\thinspace}
	\providecommand{\MR}{\relax\ifhmode\unskip\space\fi MR }
	\providecommand{\MRhref}[2]{\href{http://www.ams.org/mathscinet-getitem?mr=#1}{#2}}
	\providecommand{\href}[2]{#2}

\end{document}